\documentclass[11pt]{amsart}

\newif\iffinalrun

\usepackage{amssymb,eucal}
\usepackage[all,cmtip,poly]{xy}
\usepackage{hyperref}
    \usepackage{enumitem}

\newcommand{\WD}{{\operatorname{WD}}}
\newcommand{\rec}{{\operatorname{rec}}}
\newcommand{\Fpbar}{{\overline{\F}_p}}
\newcommand{\Fp}{{\F_p}}
\newcommand{\mc}{\mathcal}
\newcommand{\mf}{\mathfrak}
\newcommand{\m}{\mathfrak{m}}
\newcommand{\n}{\mathfrak{n}}
\newcommand{\epsilonbar    }{\overline{\epsilon}}   
\newcommand{\varepsilonbar    }{\overline{\varepsilon}}   
\newcommand{\rbar}{\bar{r}}
\newcommand{\ang}[1]{\langle #1 \rangle}

\newcommand{\GQp}{{G_\Qp}}

\newcommand{\T}{\mathbb{T}}

\newcommand{\tV}{\widetilde{{V}}}

\newcommand{\tv}{{\widetilde{{v}}}}

\newcommand{\cA}{\mathcal{A}}

\newcommand{\cG}{\mathcal{G}}

\newcommand{\cM}{\mathcal{M}}
\newcommand{\cN}{\mathcal{N}}
\newcommand{\cO}{\mathcal{O}}
\newcommand{\cP}{\mathcal{P}}

\newcommand{\cR}{\mathcal{R}}

\newcommand{\A}{\mathbb A}
\newcommand{\C}{\mathbb C}
\newcommand{\F}{\mathbb F}
\newcommand{\fp}{\F_p}

\newcommand{\fpb}{\overline \F_p} 

\newcommand{\HH}{\mathcal H}
\renewcommand{\k}{{\bar k}}
\renewcommand{\O}{{\mathcal O}}

\newcommand{\Q}{\mathbb Q}

\newcommand{\qp}{{\Q_p}}

\newcommand{\qpb}{\overline\Q_p} 

\renewcommand{\SS}{\mathcal S}
\newcommand{\TT}{\mathbb T}
\newcommand{\Z}{\mathbb Z}
\newcommand{\Zbar}{\overline{\Z}} 
\newcommand{\zp}{\Z_p}
\newcommand{\zpb}{\overline\Z_p}

\newcommand{\Zp}{{\Z_p}}

\newcommand{\into}{\hookrightarrow}

\newcommand{\onto}{\twoheadrightarrow}
\newcommand{\congto}{\xrightarrow{\,\sim\,}}
\newcommand{\isoto}{\congto}
\newcommand{\tocong}{\xleftarrow{\,\sim\,}}

\newcommand{\s}{^\times}
\renewcommand{\ss}{^{\mathrm{ss}}}
\newcommand{\lalg}{^{\mathrm{lalg}}}
\newcommand{\sm}{^{\mathrm{sm}}}
\newcommand{\Fss}{^{\mathrm{F-ss}}}

\newcommand{\ab}{^{\mathrm{ab}}}

\newcommand{\dual}{^\vee}
\newcommand{\gen}{_{\mathrm{gen}}}
\newcommand{\genw}{_{\mathrm{gen},w}}

\newcommand{\rhobar}{\overline{\rho}}

\newcommand{\U}{\operatorname{U}}
\newcommand{\GL}{\operatorname{GL}}
\newcommand{\SL}{\operatorname{SL}}

\newcommand{\HT}{\operatorname{HT}} 
\newcommand{\Fil}{\operatorname{Fil}}

\newcommand{\Qbar}{\overline{\Q}} 

\newcommand{\Qp}{{\Q_p}}
\newcommand{\Ql}{\Q_l} 
 
\newcommand{\Qpbar}{{\overline{\Q}_p}}
\newcommand{\Zpbar}{{\overline{\Z}_p}}

\newcommand{\Fbar}{\overline{\F}} 
\newcommand{\M}{\mathcal{M}} 

\newcommand{\Slift}{\widetilde{S}}

\newcommand{\gr}{\operatorname{gr}} 
\newcommand{\FBrModdd}[1][r]{\text{$\F$-$\operatorname{BrMod}_{\mathrm {dd}}^{#1}$}}
\newcommand{\FpBrMod}[1][r]{\text{$\Fp$-$\operatorname{BrMod}^{#1}$}}
\newcommand{\OEModdd}[1][r]{\text{$\cO_E$-$\Mod_{\mathrm {dd}}^{#1}$}}
\newcommand{\ZpMod}[1][r]{\text{$\Zp$-$\Mod^{#1}$}}
 
\newcommand{\Res}{\operatorname{Res}}

\DeclareMathOperator{\Max}{Max}
\DeclareMathOperator{\Mod}{Mod}
\DeclareMathOperator{\Art}{Art}
\DeclareMathOperator{\val}{val}
\DeclareMathOperator{\Id}{Id}
\DeclareMathOperator{\JH}{JH}

\DeclareMathOperator{\End}{End}
\DeclareMathOperator{\Hom}{Hom}

\DeclareMathOperator{\rank}{rank}

\DeclareMathOperator{\Gal}{Gal}

\DeclareMathOperator{\tr}{tr}

\DeclareMathOperator{\ord}{ord}

\DeclareMathOperator{\Ind}{Ind}
\DeclareMathOperator{\ind}{c-Ind}
\DeclareMathOperator{\nInd}{n-Ind}

\DeclareMathOperator{\im}{im}
\DeclareMathOperator{\Frob}{Frob}

\newcommand{\plim}{\varprojlim}  \newcommand{\ilim}{\varinjlim}   
\newcommand{\bs}{\backslash}
\renewcommand{\o}[1]{\overline{#1}}
\newcommand{\wt}[1]{\widetilde{#1}}
\newcommand{\wh}[1]{\widehat{#1}}
\newcommand{\su}[1][j]{^{(#1)}}

\newcommand{\hg}{\HH_G(V)}
\newcommand{\kind}{\ind_K^{G(F)}}
\newcommand{\vnk}{V^{N(k)}}
\newcommand{\vonk}{V_{\o N(k)}}
\newcommand{\vuk}{V^{U(k)}}
\newcommand{\vouk}{V_{\o U(k)}}

\newcommand{\bigO}{\mathcal{O}}

\newcommand{\N}{\mathbb{N}}

\newcommand{\dR}{\mathrm{dR}}
\newcommand{\cris}{\mathrm{cris}}
\newcommand{\st}{\mathrm{st}}

\newcommand{\textA}{\mathrm{A}}
\newcommand{\textB}{\mathrm{B}}
\newcommand{\textD}{\mathrm{D}}
\newcommand{\textV}{\mathrm{V}}
\newcommand{\textT}{\mathrm{T}}

\newcommand{\hA}{\wh\textA}
\newcommand{\Acris}{\textA_{\cris}}
\newcommand{\hAcris}{\wh\textA_{\cris}}
\newcommand{\hAst}{\wh\textA_{\st}}
\newcommand{\Bcris}{\textB_{\cris}}
\newcommand{\BdR}{\textB_{\dR}}

\newcommand{\Bst}{\textB_{\st}}
\newcommand{\Dcris}{\textD_{\cris}}

\newcommand{\Dst}{\textD_{\st}}
\newcommand{\Vst}{\textV_{\st}}
\newcommand{\Tst}{\textT_{\st}}

\newcommand{\D}{{\mathcal{D}}}
\newcommand{\MF}{{\mathcal{MF}}}

\iffinalrun
  \newcommand{\need}[1]{}
  \newcommand{\mar}[1]{}
\else
  \newcommand{\need}[1]{{\tiny *** #1}}
  \newcommand{\mar}[1]{\marginpar{\raggedright\tiny #1}}
\fi

\renewcommand{\(}{\textup{(}}
\renewcommand{\)}{\textup{)}}

\usepackage{hyperref}
    \usepackage{enumitem}

\theoremstyle{plain} 
 \newtheorem{ithm}{Theorem}
\newtheorem{lm}[equation]{Lemma}
\newtheorem{lem}[equation]{Lemma}

\newtheorem{prop}[equation]{Proposition}
\newtheorem{thm}[equation]{Theorem}
\newtheorem{coroll}[equation]{Corollary}
\newtheorem{corollary}[equation]{Corollary}

\theoremstyle{definition}
\newtheorem{df}[equation]{Definition}
\theoremstyle{remark}
\newtheorem{rk}[equation]{Remark}
\newtheorem{defn}[equation]{Definition}
\theoremstyle{remark} 
\newtheorem{rem}[equation]{Remark}
\newtheorem{remark}[equation]{Remark}

\numberwithin{equation}{subsection}
\numberwithin{figure}{subsection}

\newcommand{\eqninc}{\addtocounter{equation}{1}}
\newcommand{\ssinc}{\addtocounter{subsubsection}{1}}

\def\RCS$#1: #2 ${\expandafter\def\csname RCS#1\endcsname{#2}}
\RCS$Revision: 36f
8 $
\RCS$Date: 2013-04-14 11:48:38 +0800 (Sun, 14 Apr 2013) $

\begin{document}

\title{Weight cycling and Serre-type conjectures for unitary groups}
\author{Matthew Emerton}
\address{Department of Mathematics, Northwestern University, 2033 Sheridan Road, Evanston, IL 60208-2730, USA}
\email{emerton@math.northwestern.edu}
\thanks{The first author was partially supported by NSF grants DMS-0701315 and DMS-1002339}
\author{Toby Gee}
\address{Department of Mathematics, Northwestern University, 2033 Sheridan Road, Evanston, IL 60208-2730, USA}
\email{gee@math.northwestern.edu} 
\thanks{The second author was partially supported by NSF grant DMS-0841491}
\author{Florian Herzig}
\address{Institute for Advanced Study, Einstein Drive, Princeton, NJ 08540, USA}
\email{herzig@math.ias.edu}
\thanks{The third author was partially supported by NSF grant DMS-0902044 and by agreement DMS-0635607}
\maketitle

\begin{abstract}  We prove that for forms
  of $\U(3)$ which are compact at infinity and split at places dividing
  a prime $p$, in generic situations the Serre weights of a mod $p$ modular
  Galois representation which is irreducible when restricted to each
  decomposition group above $p$ are exactly those previously predicted by the
  third author.  We do this by combining explicit computations
  in $p$-adic Hodge theory, based on a formalism of strongly divisible
  modules and Breuil modules with descent data which we develop in the paper,
  with a technique that we call ``weight cycling''.
\end{abstract}

\section{Introduction}
\label{sec:introduction}

The weight part of (generalisations of) Serre's conjecture has 
received considerable attention in recent years. There have been many
new conjectures formulated (cf.~\cite{bib:BDJ},
\cite{bib:herzig-thesis}, \cite{MR2430440}, \cite{bib:Gee-lifts})
 and several cases of these conjectures have
been established (cf.~\cite{bib:Gee-bdj},
\cite{geesavitttotallyramified}). However, there have been essentially
no theoretical results for groups of semisimple rank greater than
one. In this paper we prove the first such results, by showing (under
suitable genericity hypotheses) that for a Galois representation which is
irreducible
when restricted to each decomposition group above~$p$,
its Serre weights on forms of $\U(3)$
which are compact at infinity and split at places dividing~$p$,
and for which it is modular,
are precisely those predicted by~\cite{bib:herzig-thesis}.
(The same set of weights is also predicted to be a subset of the modular weights by Doud's modification \cite{doud-supersingular}
of the conjecture of Ash, Doud, Pollack, and Sinnott \cite{bib:ADP}, \cite{bib:ASinn}.)
We remark that
other results have recently been obtained in 
\cite{gg} (in the ordinary case), and the forthcoming \cite{blggUn}.

We now explain our results in more detail.
Let $F$ be an imaginary CM field in which $p$ splits completely, and $\rbar:G_F\to\GL_3(\Fpbar)$ be a continuous
irreducible representation.  Let $F^+$ be the maximal totally real subfield
of $F$, and let $G$ be a unitary group over $F^+$ which is isomorphic
to $U(3)$ at each infinite place and split at each prime above
$p$. Assume that $F^+\ne \Q$.
In Definition~\ref{defn: modular of some Serre
  weight} we define what it means for $\rbar$ to be modular for $G$ (of some
specific weight); roughly
speaking, it means that the characteristic polynomials of Frobenius
elements at unramified places should correspond to Hecke polynomials
of some automorphic Hecke eigenform on $G$.
This implies
that $\rbar$ is essentially conjugate self-dual. We have a
notion of a (strongly) generic Serre weight (see Definition~\ref{defn:
  generic global Serre weight}; it is a condition saying that a weight
is sufficiently far away from the walls of the alcoves, and implies
for example that $p$ is at least $17$), and we let $W\gen(\rbar)$ be the set
of generic Serre weights for which $\rbar$ is modular. We also define
$W^?(\rbar)$, the set of conjectural Serre weights for $\rbar$, following
the recipe of \cite{bib:herzig-thesis}. (See
Subsection~\ref{subsection:main theorem for unitary groups}.) We
remind the reader that $W^?(\rbar)$ is defined purely locally, and
indeed only depends on the restrictions of $\rbar$ to inertia groups
at places dividing $p$. Then our main
theorem (Theorem~\ref{thm:main unitary global} in the body of the text) is:

\begin{ithm}\label{thm:main intro version}
  Suppose that $\rbar: G_F \to \GL_3(\fpb)$ is a continuous representation. 
  For all places $w | p$ of $F$ suppose that $\rbar|_{G_{F_w}}$ is irreducible.
  If $\rbar$ is modular of some strongly generic Serre weight, then $W\gen(\rbar) = W^?(\rbar)$.
\end{ithm}


This result may be regarded as a generalisation of the results of
\cite{bib:Gee-bdj} and
\cite{geesavitttotallyramified} to 3-dimensional representations, but
the methods of proof are rather different, and in particular make no
use of automorphy lifting theorems. The argument naturally breaks up
into two parts. We first prove that any generic weight for which $\rbar$
is modular is one predicted by \cite{bib:herzig-thesis}.  We then prove
that all these predicted weights actually appear.

The first part of the argument makes
use of calculations in $p$-adic Hodge theory.  While 
in essence these are of a similar nature
to those of \cite{bib:Gee-bdj} and \cite{MR2433612} (which in spirit
go back at least as far as \cite{MR1176206}),
they are rather more involved than those earlier
calculations, for several reasons.  Firstly, the representation
theory of $\GL_3(\Fp)$ is significantly more complicated than that of
$\GL_2(\Fp)$, which makes the combinatorics rather more involved. More
significantly, the calculations of \cite{bib:Gee-bdj} and
\cite{MR2433612} are for potentially crystalline Galois
representations with Hodge--Tate weights in $[0,1]$, whose reductions
mod~$p$ correspond to finite flat group schemes with descent data. In
our setting we must work instead with $3$-dimensional potentially
crystalline Galois representations with Hodge--Tate weights in
$[0,2]$, whose reductions no longer correspond to finite flat
group schemes. However, thanks to the recent work of Tong~Liu
(\cite{MR2388556}) which established a conjecture of Breuil on the
existence of strongly divisible lattices in semistable Galois
representations with Hodge--Tate weights in $[0,p-2]$ (with no
restrictions on the ramification), one knows that (the reductions mod~$p$
of) 
the representations we need to study will correspond to
certain Breuil modules with descent data.

In fact,
we are not aware of any results in the literature covering the
extension of Liu's results to incorporate coefficients and descent
data, so we establish the basic results in this paper. We do not attempt a
comprehensive treatment, but limit ourselves to allowing coefficients
in the rings of integers of finite extensions of $\Qp$, and consider
only tame descent data. Our treatment follows that of
\cite{MR2137952}, which considers the case of Hodge--Tate weights in
$[0,1]$. It will not surprise an expert that this extension is
possible, and our arguments here are of an essentially formal nature. 

Having established these results, we are reduced to performing certain
calculations with Breuil modules with descent data. Again, these
calculations are far more delicate than those of  \cite{bib:Gee-bdj}
and \cite{MR2433612}, for several reasons. Firstly, we no longer have
the underlying crutch of finite flat group schemes, which means that
we cannot immediately make use of scheme-theoretic closure techniques,
or make use of properties of \'etale models. Thankfully, the work of
Caruso (\cite{carusomodp}), and in particular his theory of maximal
Breuil modules, provides an adequate substitute for our purposes. In
addition, we know of no way to make use of the fact that our
representations are potentially crystalline, rather than just
potentially semistable, and we do not know of any way to directly
transfer the information that the Hodge--Tate weights of our
representations are ${0,1,2}$ to the Breuil modules we
consider. Instead, we assume only that our representations are
potentially semistable with Hodge--Tate weights in $[0,2]$, and use
various ad hoc arguments using the assumption that our mod~$p$
representations are (locally) irreducible. We also make significant
use of the determinant of the representation (in particular, this is
how we exploit the fact that our Hodge--Tate weights are
$0,1,2$).

Having in this way established the ``easy'' direction of the weight
part of Serre's conjecture, we now tackle the ``hard'' direction of
proving that a mod~$p$ representation that is modular of some particular
weight is actually modular of all the
expected weights. In contrast to \cite{bib:Gee-bdj} and \cite{blggUn},
we do not make use of any automorphy lifting theorems to do
this. Instead, we make use of a method which we call
``weight cycling''.
In the case of the group  ${\GL_2}_{/\Qp}$, this technique was discovered by
Kevin Buzzard,\footnote{In fact, the computations that underly the technique seem to go
back (at least) to the proofs of Thm.~3.2 of the paper
\cite{MR629471} of Hida and Thm.~(3.4) of the paper \cite{MR688264} of Ribet.}
and was first written up in Section~5 of
\cite{MR2290604}. The argument was extended to  ${\GL_2}_{/F}$ for $F$ an
arbitrary finite extension of $\Qp$ in \cite{MR2407230}, and was also
exploited in \cite{geesavitttotallyramified}. We extend it
to arbitrary split
connected reductive groups over $p$-adic fields. The result we obtain
is of the following form: if a mod~$p$ representation $\rbar$ is
modular of some given weight, and a certain Hecke operator at~$p$
vanishes on the corresponding cohomology class, then $\rbar$ is
modular of some other weight in a short list of possibilities. 
Using the
elimination result already proved, this list can often be reduced to a
single weight. The argument can then be repeated, and in this way we
can ``cycle'' through all the conjectural weights.

Of course, in order to apply the weight cycling result we need to establish
a vanishing result for certain Hecke operators in characteristic~$p$
(Corollary~\ref{cor:supersingular-hecke-evals}). 
Here we again use our assumption that our mod~$p$ Galois
representation is irreducible when restricted to any decomposition
group above~$p$. We establish a rather general comparison between
Hecke operators at~$p$ in characteristic zero and characteristic~$p$
(Proposition~\ref{prop: compatibility of Hecke actions in char 0 and char p}),
and use it to show that if the relevant Hecke operators did not
vanish, certain (normalised) Hecke operators in characteristic 0 would
have to act via $p$-adic units. Using local-global compatibility one can deduce
that the local Galois representation would be reducible, a
contradiction. This argument may be regarded as a generalisation of the
fact that the Galois representations associated to ordinary modular
forms are necessarily reducible at $p$.

We remark that while automorphy lifting theorem techniques such as
those used in \cite{blggUn} are undeniably a powerful method for
proving results about Serre weight conjectures, it does not seem to be
possible to prove Theorem \ref{thm:main intro version} by such methods. The problem is
that these methods rely on producing automorphic Galois
representations in characteristic 0, and deducing information about
Serre weights by reducing modulo~$p$. However, irreducible
representations of the algebraic group $\GL_3$ over $\Fpbar$ cannot always be lifted to characteristic 0,
 making the link between
characteristic 0 and characteristic~$p$ far weaker than it is for
$\GL_2$. Concretely, if we are working over a CM field $F$ of degree
$2d$ over $\Q$ in which $p$ splits completely, then we generically expect to have $9^d$ Serre weights
for a given Galois representation (assumed irreducible when restricted
to decomposition groups at places above $p$), and in this paper we produce all
$9^d$ weights (in those cases to which our results apply).
By contrast, automorphy lifting techniques appear to be
limited to producing $6^d$ weights (the ``obvious'' weights; for
example, in the case that the restriction to the decomposition groups
at $p$ was ordinary, these would be the ordinary weights); see the
forthcoming paper \cite{blggUn} for this.

It would be of interest to
weaken the hypothesis of (strong) genericity in our results. This
would appear to require a significant improvement to our $p$-adic
Hodge theoretic techniques, in order to prove the ``elimination'' part
of the argument; in particular, it seems likely that one would have to
make full use of the precise Hodge--Tate weights, which we do not
currently know how to do, or to understand the reductions modulo $p$
of three-dimensional crystalline representations with Hodge--Tate
weights in $[0,2p]$, which also seems to be a hard problem.
Even if the ``elimination'' part can be completed, the conjectural weight
set need not consist of $9^d$ weights in degenerate situations and it
seems that weight cycling does not necessarily cycle through all these 
weights.

\subsection{The organisation of the paper}
In order to make it easier for the reader to see exactly which
assumptions are used in each argument, and with an eye to future
applications, we have gone to considerable lengths to axiomatise as
many of our assumptions as possible, and to argue within our abstract
axiomatic framework, rather than with specific unitary groups. For
example, where arguments apply for a general $\GL_n$, we have
presented them in this generality.

In Section~\ref{sec:repr-theory} we recall the definitions of 
weights and their associated Hecke algebras, as well as the Satake
isomorphism of \cite{bib:herzig-satake}. We also prove our general
weight cycling result for an arbitrary split connected reductive
group, and establish some special cases of the tame inertial local
Langlands correspondence for $\GL_n$.

In Section~\ref{sec:p-adic-hodge}
we establish the necessary results from $p$-adic Hodge theory; in
particular, we generalise some of the results of \cite{MR2388556} to
allow for coefficients and descent data. We then use these results
to give an upper bound for the set of
possible reductions mod~$p$ of 3-dimensional potentially semistable
representations of $\GQp$ with Hodge--Tate weights $0$, $-1$, $-2$
which are of particular tame types, under the
assumption that these reductions are irreducible. 

In Section~\ref{sec: abstract framework} we introduce our basic
axiomatic framework. We have sets of axioms dealing with
both characteristic $p$ and characteristic zero contexts,
and we show how the characteristic zero axioms
imply the characteristic $p$ ones; this is an abstraction of the
results on local-global compatibility and the comparison of Hecke
operators in characteristics zero and $p$ discussed above.

In Section~\ref{sec:elimination} we use the calculations of
Section~\ref{sec:p-adic-hodge} to establish the elimination result
that we need; that is, to show that in generic situations, if $\rbar$
is modular of some weight then that weight is one of the weights
predicted in \cite{bib:herzig-thesis}. We carry out these arguments in
the setting of Section~\ref{sec: abstract framework}. 

In
Section~\ref{subsec:weight-cycling-gl_3} we combine the results
of Section~\ref{sec:elimination} with our weight
cycling technique and the vanishing of Hecke operators proved in
Subsection~\ref{subsec:niveau-3-implies} to prove an abstract version of
Theorem~\ref{thm:main intro version} (namely Theorem~\ref{thm:main})
in our axiomatic framework.

In Section~\ref{sec:global-applications} we begin by recalling some standard
material on automorphic forms on compact unitary groups and their
associated Galois representations. We then establish that the axioms
of Section~\ref{sec: abstract framework} hold in this setting, and
thus establish Theorem~\ref{thm:main intro version}. Finally, in
Subsection~\ref{examples:automorphic induction} we use automorphic
induction and base change for unitary groups to show that there are
many Galois representations which satisfy the hypotheses of our main
theorem. We remark that the assumption that $F^+\ne\Q$ is only used in
the proof of Theorem~\ref{thm: existence of Galois reps attached to
  algebraic modular forms}, where it is needed due to a limitation on
our knowledge of base change between unitary groups and $\GL_3$.

We refer to the beginning of each section for more details regarding its
contents.

\subsection{Acknowledgements}
\label{subsec:acknowledgements}

We would like to thank Xavier Caruso, David Savitt, and Tong Liu for some helpful conversations related to the material in Section~\ref{sec:p-adic-hodge};
in particular, we are grateful to Caruso for the proof of
Proposition~\ref{prop:descent-data-on-rank-one-subobject}.
We thank Guy Henniart for helpful comments regarding the local
Langlands correspondence and Sug Woo Shin for a helpful discussion
concerning the proof of Theorem~\ref{thm: existence of Galois reps
  attached to algebraic modular forms}. We thank Diego Izquierdo for
his helpful comments on an earlier version of the paper. 
Some of this work was carried
out at the IAS (by F.H.) and during a conference visit at the
University of Kyoto in November 2008 (by T.G. and F.H.). We thank
these institutions for the excellent working conditions they
provided. Finally we thank the referees for many helpful
comments.


\subsection{Notation and terminology}
\label{subsec:notation-terminology}

For any field $K$, we let $G_K$ denote an absolute Galois group of
$K$, the precise choice of $G_K$ depending on the choice of an
algebraic closure $\overline{K}$ of $K$.  In the case of the field
$\mathbb Q_p$, it will be convenient to fix once and for all an
algebraic closure $\Qpbar$, with ring of integers $\Zpbar$ and residue
field $\Fpbar$. If $K$ is an algebraic extension of
$\mathbb Q_p$,
it will always be implicit that $K$ is a subfield
of this fixed algebraic closure $\Qbar_p$, and we may then
unambiguously set $G_K := \Gal(\Qbar_p/K)$.  It is also convenient to
fix once and for all an isomorphism $\imath:\Qpbar\isoto\C$, and we do
so.

Let $F$ be a finite extension of $\Ql$, for $l$ possibly equal to $p$. 
We let $\Art_F$ be the isomorphism $F^\times\isoto W_F\ab$ of local class field theory, normalised so that geometric
Frobenius elements correspond to uniformisers. Let $\rec_{F,\C}$ denote the local Langlands correspondence from isomorphism
classes of irreducible smooth representations of $\GL_n(F)$ over $\C$ to isomorphism classes of $n$-dimensional Frobenius
semisimple Weil--Deligne representations of $W_F$. (See \cite{ht}.) We define
the local Langlands correspondence $\rec_F$ over $\qpb$ by $\imath \circ \rec_F = \rec_{F,\C} \circ \imath$. It depends only 
on $\imath^{-1}(q^{(n-1)/2})$, where $q$ is the cardinality of the residue field.
Let $N(\pi)$ denote the monodromy operator of
$\rec_F(\pi)$.

We let $\Qbar$ denote the algebraic closure of $\mathbb Q$ in
$\Qbar_p$.
All algebraic extensions $F$ of $\mathbb Q$ are implicitly understood to be subfields
of $\Qbar$, so that we may unambiguously define $G_F := \Gal(\Qbar/F)$. If
$v$ is a finite place of $F$ we denote by $\Frob_v \in G_F$ a (choice of) geometric Frobenius element
at $v$.

Where possible, we use $\rho$ and $\rhobar$ to denote representations
of the absolute Galois groups of local fields, and $r$ and $\rbar$ to
denote representations of the absolute Galois groups of number
fields. We let $\varepsilon$ denote the $p$-adic cyclotomic character,
and we let $\varepsilonbar$ denote its reduction modulo $p$.

If $K$ is a finite extension of $\Qp$ and $\rho$ is a continuous de Rham representation of $G_K$ over $\Qpbar$,
then we will write $\WD(\rho)$ for the corresponding
Weil--Deligne representation of $W_K$ (defined, for example, in
Appendix B.1 of \cite{MR1639612}), 
and if $\tau:K \into \Qpbar$ then we will write $\HT_\tau(\rho)$ for the multiset of
Hodge--Tate weights of $\rho$ with respect to $\tau$.  By definition, if $W$ is a de Rham representation of $G_K$ over
$\Qpbar$ and if $\tau:K \into \Qpbar$ then the multiset $\HT_\tau(W)$ contains $i$ with multiplicity $\dim_{\Qpbar} (W
\otimes_{\tau,K} \widehat{\o{K}}(i))^{G_K} $. Thus for example
$\HT_\tau(\varepsilon)=\{ -1\}$.
We will refer to (the
isomorphism class of) $\WD(\rho)|_{I_K}$ as the \emph{inertial type}
of $\rho$. We will let $\WD(\rho)\Fss$ denote the Frobenius
semisimplification of $\WD(\rho)$.

We will refer to the mod $p$ Galois representations that are
associated to automorphic representations as {\em modular}, as this
seems to be standard practice in the literature; we remark that one
could equally well refer to them as {\em automorphic}.

We also signal to the reader that we will frequently use the
word {\em weight} to refer to an irreducible representation (or
an isomorphism class thereof) of the
$k$-valued points of a reductive group over a finite field 
$k$.  (See Definition~\ref{df:serre-wts} below.)
A {\em Serre weight} is then a variant of this notion, defined
in a suitable global context.
(See Definition~\ref{df:serre-weight}.)
Hopefully no confusion will arise with the other standard
use of the word ``weight'' in representation theory (i.e.\ in
the sense of a character of a maximal torus).

\section{Representation theory}
\label{sec:repr-theory}
This section is devoted to various representation-theoretic
preliminaries.  In Subsection~\ref{subsec:intro-to-weights} we recall
some basic terminology and facts related to weights, and in
Subsection~\ref{subsec:hecke-operators} we recall the Satake
isomorphism of~\cite{bib:herzig-satake}.  In
Subsection~\ref{subsec:weight-cycling} we then introduce the technique
of weight cycling in a general setting.  Finally, in
Subsection~\ref{subsec:inertial-ll}, which is in a somewhat different
vein to the preceding three subsections, we establish some simple
instances of the so-called inertial local Langlands correspondence for
the group $\GL_n(F)$, with $F$ a finite extension of $\Q_p$.

\subsection{Weights}
\label{subsec:intro-to-weights}
We introduce notation that will be in force for this subsection and the two
that follow.  Namely,
let~$F$ be a finite extension of~$\qp$ with ring of integers~$\O_F$, uniformiser~$\varpi$, and residue field~$k$.
Let $G_{/\O_F}$ be a split connected reductive group and fix a maximal split torus $T_{/\O_F}$. Let $\Phi \subset X^*(T)$ denote
the set of roots and choose a system of simple roots~$\Delta \subset \Phi$. Denote by $B_{/\O_F}$ denote the associated Borel
subgroup and by $U_{/\O_F}$ its unipotent radical. Let~$W$ be the Weyl group and $K := G(\O_F)$,
a hyperspecial maximal compact subgroup of~$G(F)$.

\begin{df}\label{df:serre-wts}
  A \emph{weight} is an isomorphism class of irreducible representations~$V$ of~$G(k)$ over~$\k$.
\end{df}

Since $G(k)$ is a finite group, there are only finitely many weights.

\begin{remark}
\label{rem:abuse}
We will constantly engage in the following standard abuse of terminology:
namely, we will speak of some particular 
$G(k)$-representation $V$ being a weight,
when we actually mean that $V$ is an isomorphism class representative of 
a weight.
\end{remark}

For a standard parabolic subgroup $P = MN$ we denote by $\o P = M\o N$ the opposite parabolic.
The following result
is Lemma~2.5 in~\cite{bib:herzig-satake}.

\begin{lm}\label{lm:serre-wts-invts}
  Suppose that $V$ is a weight and that $P = MN$ is a standard parabolic. Then $\vnk$ and $\vonk$ are weights for
  $M$ and the natural, $M(k)$-equivariant map $\vnk \into V \onto \vonk$ is an isomorphism. In particular, $\vuk \cong \vouk$ is
  one-dimensional.
\end{lm}

Suppose that $\mu \in X_*(T)$. Let $P_\mu = M_\mu N_\mu$ denote the parabolic subgroup of $G$
defined by $\mu$, where the Levi subgroup $M_\mu$ contains $T$.
(See \cite[Prop.\ 8.4.5]{bib:Springer_LAG}. In particular, if $\mu$ is dominant
then $P_\mu$ is the standard parabolic subgroup defined by the set of simple roots that are orthogonal
to $\mu$.)
The following lemma is a special case of Proposition~3.8 in~\cite{bib:herzig-satake}.
We use the shorthand ${}^t K = t K t^{-1}$.

\begin{lm}\label{lm:buildings-lemma}
  Let $t = \mu(\varpi)$. Then the image of $K \cap {}^t K$ in $G(k)$ is $P_{-\mu}(k)$.
\end{lm}

\subsection{Hecke operators and the Satake isomorphism}
\label{subsec:hecke-operators}

If $V$ is a weight, then we define the
\emph{Hecke algebra of $V$} as $\HH_G(V) := \End_{G(F)}(\kind V)$.  By
Frobenius reciprocity
we can and usually will think of it as $\k$-algebra of compactly supported functions $f : G(F) \to \End_{\k}
V$ satisfying $f(k_1gk_2) = k_1 \circ f(g)\circ k_2$ for all $k_1$, $k_2 \in K$, $g \in G(F)$, where the multiplication is given by
convolution. Note that if $\pi$ is a smooth $G(F)$-representation,
$\HH_G(V)$ naturally acts on the left on $(V \otimes_{\fpb} \pi)^K$. Explicitly, if $\phi \in \HH_G(V)$ and
$x \in (V \otimes_{\fpb} \pi)^K$, then
\begin{equation}\label{eq:20}
\phi \cdot x = \sum_{\gamma \in G(F)/K} (\phi(\gamma) \otimes \gamma)x.
\end{equation}

We now recall some results of~\cite{bib:herzig-satake}.
Let~$T^-$ denote the submonoid of~$T(F)$,
\begin{equation*}
  T^- = \{t \in T(F) : \ord_F(\alpha(t)) \le 0 \quad \forall \alpha \in \Delta\},
\end{equation*}
and let $\HH^-_T(V^{U(k)})$ denote the subalgebra of~$\HH_T(V^{U(k)})$ consisting of those $\varphi : T(F) \to \k$ that are
supported on~$T^-$.

\begin{thm}\label{thm:satake} Suppose that $V$ is a weight.
  Then
  \begin{align*}
    \SS_G : \HH_G(V) &\to \HH_T(V^{U(k)})\\
    f &\mapsto \left(t \mapsto \sum_{u \in U(F)/U(\O_F)} f(tu)\Big|_{V^{U(k)}}\right)
  \end{align*}
  is an injective $\k$-algebra homomorphism with image $\HH^-_T(V^{U(k)})$.
\end{thm}

In particular, $\hg \cong \k[X_*(T)_-]$ is commutative and noetherian (by Gordan's lemma).
Here, $X_*(T)_- = \{ \mu \in X_*(T) : \ang{\mu,\alpha} \le 0\ \forall \alpha \in \Delta \}$.
We recall that $G(F) = \coprod K\mu(\varpi)K$, where
$\mu$ ranges over $X_*(T)_-$ (refined Cartan decomposition). Moreover,
$\hg$ has a basis (as a $\overline{k}$-vector space) given by $\o T_{\mu,\varpi}$ ($\mu
\in X_*(T)_-$), where $\o T_{\mu,\varpi}$ has support $K\mu(\varpi)K$ and sends $\mu(\varpi)$ to the endomorphism $V \onto
V_{N_\mu(k)} \tocong V^{N_{-\mu}(k)} \into V$.  Note that up to scalar multiple, $\o T_{\mu,\varpi}$ is independent of the choice of uniformiser $\varpi$.

\subsection{Weight cycling}
\label{subsec:weight-cycling}
The following proposition and its corollary encode the basic technique of weight cycling.

\begin{prop}\label{prop:weight-cycling}
  Suppose that $V$ is a weight and that $\pi$ is a smooth $G(F)$-representation.  Suppose that $\mu
  \in X_*(T)_-$ and let $t = \mu(\varpi)$. Then we have a commutative
  diagram as follows, where the map on the right is
  induced by Frobenius reciprocity and the map on the top is \emph{injective}.
  \begin{equation*}
    \xymatrix{(V\otimes_{\fpb} \pi)^K \ar@{^{(}-->}[r]^-i \ar[dr]_{\o T_{\mu,\varpi}} & 
      \big((\Ind_{K\cap {}^t K}^K V^{N_{-\mu}(k)})\otimes_{\fpb} \pi\big)^K \ar[d] \\ & (V\otimes_{\fpb} \pi)^K}
  \end{equation*}
\end{prop}

Note that the natural map $V^{N_{-\mu}(k)} \to V$ is $K \cap {}^t K$-equivariant by Lemma~\ref{lm:buildings-lemma}.

\begin{proof}
  Suppose that $x \in (V\otimes_{\fpb} \pi)^K$. By~\eqref{eq:20} and the definition of $\o T_{\mu,\varpi}$ we see that $\o
  T_{\mu,\varpi} \cdot x = \sum_{k \in K/K\cap {}^t K} (k \otimes k) y$, where \[ y := (\o T_{\mu,\varpi}(t) \otimes t) x \in
  (V^{N_{-\mu}(k)} \otimes_{\fpb} \pi)^{K\cap {}^t K}. \]
  We define the map $i$ by composing the map $x \mapsto y$ with the isomorphism
  \begin{equation}\label{eq:10}
    (V^{N_{-\mu}(k)} \otimes_{\fpb} \pi)^{K\cap {}^t K} \congto
    \big((\Ind_{K\cap {}^t K}^K V^{N_{-\mu}(k)})\otimes_{\fpb} \pi\big)^K,
  \end{equation}
  which holds by Frobenius reciprocity.
  Explicitly, the isomorphism in~\eqref{eq:10} is induced by $V^{N_{-\mu}(k)} \otimes_{\fpb} \pi \to
  (\Ind_{K\cap {}^t K}^K V^{N_{-\mu}(k)})\otimes_{\fpb} \pi$ sending $v \otimes z$ to $\sum_{k \in K/(K\cap {}^t K)} [k,v] \otimes kz$,
  where $[k,v]$ denotes the element of the induction that is supported on $(K \cap {}^t K)k^{-1}$ and sends $k^{-1}$ to $v$.
  Also, the natural map $\Ind_{K\cap {}^t K}^K V^{N_{-\mu}(k)} \to V$ sends $[k,v]$ to $kv$. It follows that the above
  diagram commutes.

  Let $\eta : V \onto V_{N_\mu(k)} \onto \vouk \congto \fpb$ in $V\dual$, where 
  the identification with $\fpb$ is arbitrary. Note that $i(x) = 0$ implies that $\eta(x) = 0$ in $\pi$. Since $i$ is $K$-equivariant
  and $\eta$ generates $V\dual$ as $K$-representation, it follows that $V\dual$ kills $\ker(i) \subset V \otimes \pi$, so $i$ is injective.
\end{proof}

\begin{coroll}\label{cor:weight-cycling}
  With the notation as in the proposition, if $\o T_{\mu,\varpi}$ fails to be injective on $(V\otimes_{\fpb} \pi)^K$,
  then $(V'\otimes_{\fpb} \pi)^K \ne 0$ for at least one of the irreducible constituents $V'$ of
  $\ker(\Ind_{K\cap {}^t K}^K V^{N_{-\mu}(k)} \to V)$.
\end{coroll}

\subsection{Inertial local Langlands}
\label{subsec:inertial-ll}

The purpose of this subsection is to establish some simple instances of ``inertial local Langlands'' for
$\GL_n(F)$, where $n \ge 1$ is arbitrary and $F$ is a finite extension
of $\qp$. (We remark that the $n = 2$ case is worked out completely in the appendix of \cite{MR1944572}.)

Let $E/F$ denote the unramified extension of degree $n$. Let $\O_E$ (resp.\ $\O_F$) denote the ring of integers of $E$
(resp.\ $F$), and let $k_E$ (resp.\ $k_F$) denote the residue field of $E$ (resp.\ $F$). 
Fix an $\O_F$-basis of $\O_E$. It gives rise to a rational maximal torus $\TT$ of ${\GL_n}_{/k_F}$ such that $k_E\s =
\TT(k_F)$. A character $\theta : k_E\s \to \qpb\s$ is said to be \emph{primitive} if its $\Gal(k_E/k_F)$-conjugates are all
distinct.  If this is the case, the Deligne--Lusztig representation $(-1)^{n-1} R_{\TT}^\theta$ of $\GL_n(k_F)$ over $\qpb$
is genuine, irreducible, and cuspidal. (See Proposition~7.4 and Theorem~8.3 in~\cite{bib:DL}.)

Let $B_n = T_n U_n$ denote the Borel subgroup of $\GL_n$ of upper triangular matrices. Let $I(1)$ denote the pro-$p$ Iwahori
subgroup. It is the inverse image of $U_n(k_F)$ in $\GL_n(\O_F)$.
We will
sometimes regard $\GL_n(k_F)$-representations as
$\GL_n(\cO_F)$-representations via inflation. If $\chi$ is
a character of $k_F^\times$ (respectively $\theta$ is a character of
$k_E^\times$) we will write $\chi\circ\Art_F^{-1}$ (respectively
$\theta\circ\Art_E^{-1}$) for the character of $I_F$ given by
composition with the homomorphism $I_F \to k_F\s$ induced by $\Art_F^{-1}$ (respectively with
the homomorphism $I_E \to k_E\s$ induced by $\Art_E^{-1}$).

\begin{prop}\label{prop:inertial-llc}
  Suppose that $\pi$ is an irreducible admissible smooth representation of $\GL_n(F)$ over $\qpb$.
  \begin{enumerate}
  \item If $\pi|_{\GL_n(\O_F)}$ contains the cuspidal representation $(-1)^{n-1} R_{\TT}^\theta$ for some primitive character
    $\theta : k_E\s \to \qpb\s$, then \[ \rec_F(\pi)|_{I_F} \cong \bigoplus_{\sigma \in \Gal(k_E/k_F)} \sigma(\theta\circ
    \Art_E^{-1})\] and $N(\pi) = 0$.
  \item If $\pi|_{\GL_n(\O_F)}$ contains the principal series representation
    \begin{equation*}
      \Ind_{B_n(k_F)}^{\GL_n(k_F)} (\chi_1 \otimes \cdots \otimes \chi_n)
    \end{equation*}
    for some \emph{distinct} characters
    $\chi_i : k_F\s \to \qpb\s$, then $\rec_F(\pi)|_{I_F} \cong \bigoplus_{i=1}^n \chi_i \circ \Art_F^{-1}$ and $N(\pi) = 0$.
  \end{enumerate}
\end{prop}

Note that the cuspidal representation in part~(i) and the principal series representation in part~(ii)
are irreducible. (In the latter case this is because the $\chi_i$ are distinct.)

\begin{proof}
  Using our isomorphism $\imath: \qpb \congto \C$ we will think of $\pi$, $\theta$, $\chi_i$, etc., as complex
  representations for the remainder of this proof.  Recall from \S\ref{subsec:notation-terminology} that we have $\imath \circ
  \rec_F = \rec_{F,\C} \circ \imath$.

  For the first part, let $\tau$ denote $(-1)^{n-1} R_{\TT}^\theta$. It has central character $\theta|_{k_F\s}$. Thus there
  is a unique extension of the inflation of $\theta$ to $\O_E\s$ to a character $\wt\theta$ of $E\s$ such that $\wt\theta|_{F\s}$ is
  the central character of $\pi$.  Moreover $\tau$ extends to a $F\s \GL_n(\O_F)$-representation by letting $F\s$ act via
  $\wt\theta$, and Frobenius reciprocity gives a $\GL_n(F)$-equivariant map $\ind_{F\s \GL_n(\O_F)}^{\GL_n(F)} \tau \to
  \pi$. From
  Proposition~7.3 in~\cite{bib:DL} and \cite[(2.1.1)]{bib:BH3} it follows that $\tau$ is
  uniquely characterised among cuspidal representations by the identity $\tr \tau(x) = (-1)^{n-1} \sum_{\sigma \in \Gal(k_E/k_F)}
  \theta({}^\sigma x)$ for all $x \in k_E\s$ such that all $\Gal(k_E/k_F)$-conjugates of $x$ are distinct. 
  From Theorem~2 in \cite[\S 2.4]{bib:BH3} we see that $\pi \cong \ind_{F\s \GL_n(\O_F)}^{\GL_n(F)}\tau$ is the supercuspidal
  representation that is automorphically induced from $\mu_0 \wt\theta : E\s \to \C\s$, where $\mu_0$ is the unramified 
  character sending uniformisers to $(-1)^{n-1}$.
  It follows by induction on~$n$
  from Lemma VII.2.6(6) in \cite{ht} that $\rec_{F,\C}(\pi) \cong \Ind_{W_E}^{W_F} (\mu_0\wt\theta \circ
  \Art_E^{-1})$. Since $\mu_0$ is unramified and this representation is irreducible, the claim follows.
  
  For the second part, let $\chi : T_n(k_F) \to \C\s$ denote the character $\chi_1 \otimes \cdots \otimes \chi_n$.  If $I
  \subset \GL_n(\O_F)$ denotes the inverse image of $B_n(k_F)$, then we can define a smooth character $\rho : I \to \C\s$ as
  the composite of the natural map $I \onto T_n(k_F)$ with $\chi$.  Theorem~7.7 and Remark~7.8 in \cite{bib:roche} show that we have
  an isomorphism $\pi^\rho \congto \pi_{U_n}^\chi$ where $\pi_{U_n}$ is the Jacquet module and where the superscripts on either side
  denote isotypic components. Our assumption, together with Frobenius reciprocity, shows that $\pi^\rho \ne 0$, so
  $\pi_{U_n}^\chi \ne 0$. It follows that there is a $T_n(F)$-equivariant map $\pi_{U_n} \onto \wt\chi$ for some character $\wt\chi :
  T_n(F) \to \C\s$ such that $\wt\chi|_{T_n(\O_F)} = \chi$. By Frobenius reciprocity, $\pi$ is a subrepresentation of
  $\nInd_{B_n(F)}^{G(F)} (\wt\chi \delta_{B_n}^{-1/2})$ (normalised induction, where $\delta_{B_n}$ is the modulus character of $B_n$).
  The latter representation is irreducible by our assumption that the $\chi_i|_{k_F\s}$ are distinct. (See
  Theorem~1.2.1 in~\cite{bib:Kudla}.) Thus $\pi \cong \wt\chi_1 |.|^{(1-n)/2} \boxplus \cdots \boxplus \wt\chi_n |.|^{(n-1)/2}$, so
  $\rec_{F,\C}(\pi) \cong \bigl(\wt\chi_1 |.|^{(1-n)/2}\oplus \cdots \oplus \wt\chi_n |.|^{(n-1)/2}\bigr)\circ \Art_F^{-1}$ and $N(\pi) =
  0$ \cite[p.\ 252]{ht}. As $|.|$ is unramified, this proves the proposition.
\end{proof}

\section{$p$-adic Hodge theory}
\label{sec:p-adic-hodge}
In this section we develop the various results in $p$-adic Hodge theory that
will be required in the sequel.  In Subsections~\ref{subsec:strongly-divis}
and~\ref{subsec:breuil-modules-with-dd} we introduce basic definitions and results
related respectively to strongly divisible modules and to Breuil modules;
the novelty (to the extent that there is any) is that we allow coefficients and 
descent data.  In fact, as noted in the introduction, we do not allow the 
most general form of either: rather we restrict to coefficients in 
the ring of integers $\O_E$ of a finite extension $E$ of $\mathbb Q_p$
(or its residue field $\mathbb F$ in the case of Breuil modules),
and we restrict to descent data for tamely ramified extensions.  

In Subsection~\ref{subsec:breuil-modules-two} we pass to a detailed study
of the particular Breuil modules that are relevant to our applications.
This subsection culminates in Theorem~\ref{thm: explicit list of possible characters for
specific descent data},
which is the key result that we will apply in the sequel.

\subsection{Strongly divisible modules with coefficients and descent data}
\label{subsec:strongly-divis}

In this subsection we extend certain results of
\cite{MR2137952} on the categories of strongly divisible modules with
coefficients and descent data to the case of lattices in potentially
semistable Galois representations with Hodge--Tate weights in $[0,p-2]$
(\cite{MR2137952} works only with representations with Hodge--Tate
weights in $[0,1]$). We do not attempt to work in the same level of
generality as \cite{MR2137952}, and in particular we only allow
coefficients in the ring of integers of a finite extension of $\Qp$.

We begin by recalling some results from Section~2.2 of
\cite{MR2137952}. Let $p$ be an odd prime
and let $E$ and $K$ be finite extensions of $\Qp$ inside
$\Qpbar$. (Throughout this section $K$ will be a field, rather than a
maximal compact subgroup as in section \ref{sec:repr-theory}.) Let
$K_0$ be the maximal unramified subfield of $K$, so that
$K_0=W(k)[1/p]$, where $k$ is the residue field of $K$. Let $\phi \in \Gal(K_0/\qp)$
denote the arithmetic Frobenius. Let $K/K'$ be
a Galois extension, with $K'$ a field lying between $\Qp$ and $K$. Fix
the uniformiser $p\in\Qp$, so we have a fixed embedding
$\Bst\into\BdR$.

\begin{defn}
  A \emph{filtered}-$(\varphi,N,K/K',E)$-\emph{module} of rank $n$ is a
  free $K_0\otimes_{\Qp}E$-module $D$ of rank $n$ together with
  \begin{itemize}
  \item a $\phi \otimes 1$-semilinear automorphism $\varphi$,
  \item a nilpotent $K_0\otimes_{\Qp} E$-linear endomorphism $N$
    such that $N\varphi=p\varphi N$,
  \item a decreasing filtration $(\Fil^i D_K)_{i \in \Z}$ on $D_K=K\otimes_{K_0}D$ consisting of $K\otimes_{\qp} E$-submodules, 
    which is exhaustive and separated, and
  \item a $K_0$-semilinear, $E$-linear action of $\Gal(K/K')$ which
    commutes with $\varphi$ and $N$ and preserves the filtration on $D_K$.
  \end{itemize}
We say that such a module $D$ is \emph{weakly admissible} if the underlying filtered
$(\varphi,N,K/K,\Qp)$-module is weakly admissible.
\end{defn}

Given a potentially semistable representation $\rho:G_{K'}\to\GL(V)$
on an $n$-dimensional $E$-vector space $V$, such that $\rho|_{G_K}$ is
semistable, we set \[\Dst^{K'}(V)=\Hom_{G_K}(V,\Bst),\]a weakly
admissible filtered $(\varphi,N,K/K',E)$-module of rank $n$. In the other
direction, given a weakly admissible filtered $(\varphi,N,K/K',E)$-module
$D$, we
define \[\Vst^{K'}(D)=\Hom_{\varphi,N}(D,\Bst)\cap\Hom_{\Fil}(D_K,\BdR),\]
a finite-dimensional $E$-vector space with an action of $G_{K'}$
via \[(gf)(x)=g(f(\bar{g}^{-1}x))\] where $g\in G_{K'}$ and
$\bar{g}$ is the image of $g$ in $\Gal(K/K')$.
\begin{prop}\label{prop: D and V are equivs, savitt}
  The functors $\Dst^{K'}$ and $\Vst^{K'}$ are mutually quasi-inverse, and provide equivalences
  of categories between the category of weakly admissible filtered
  $(\varphi,N,K/K',E)$-modules and the category of $E$-representations of
  $G_{K'}$ which become semistable upon restriction to $G_K$.
\end{prop}
\begin{proof}
  This follows from Proposition 2.9 of \cite{MR2137952}.
\end{proof}

We now wish to consider $\bigO_E$-lattices inside potentially
semistable Galois representations. This entails the consideration of
strongly divisible modules. We begin with some basic definitions.

We suppose from now on that $K/K'$ is a tamely ramified Galois
extension with ramification index $e(K/K')$, and that we have fixed a
uniformiser $\varpi\in K$ with $\varpi^{e(K/K')}\in K'$. Write
$g(\varpi)=h_g\varpi$ for each $g\in\Gal(K/K')$. Note that $h_g \in W(k)$. Let $e$ be the absolute
ramification index of $K$. Let $E(u) \in W(k)[u]$ be the (monic) minimal polynomial of
$\varpi$ over $K_0$. Let $S$ be the $p$-adic
completion of $W(k)[u,\frac{u^{ie}}{i!}]_{i\in\N}$. We let $\varphi : S \to S$
denote the unique continuous, $\phi$-semilinear map with $\varphi(u)=u^p$ and
$\varphi(u^{ie}/i!)=u^{iep}/i!$, and we let $N$ be the unique continuous,
$W(k)$-linear derivation of $S$ with $N(u)=-u$ and
$N(u^{ie}/i!)=-ieu^{ie}/i!$, so that $N\varphi=p\varphi N$. There is a descending
filtration $(\Fil^{i} S)_{i \ge 0}$ on $S$ given by letting $\Fil^{i} S$ be the $p$-adic
completion of the ideal of $S$ generated by $E(u)^j/j!$ for $j\geq
i$ for each $i\ge 0$. For $i\le p-1$ we have $\varphi(\Fil^{i}S)\subset
p^{i}S$, and we let $\varphi_{i}$ denote $\varphi/p^{i}$ on
$\Fil^{i}S$. We let the group $\Gal(K/K')$ act on $S$ by defining
for each $g\in\Gal(K/K')$ the continuous ring isomorphism $\hat{g}:S\to S$
with $\hat{g}(w_i\frac{u^i}{\lfloor i/e\rfloor
  !})=g(w_i)h_g^i\frac{u^i}{\lfloor i/e\rfloor !}$ (where $w_i\in W(k)$).
This action commutes with $\varphi$ and $N$ and preserves the filtration. (Note that
$\hat g(E(u)) = E(u)$ due to our assumption that $\varpi^{e(K/K')}\in K'$.)

Let $S_{\bigO_E}= S \otimes_\Zp \bigO_E$, and extend the definitions of
$\Fil$, $\varphi$, $\varphi_k$, $N$, and $\hat{g}$ to $S_{\bigO_E}$ in the
obvious $\bigO_E$-linear fashion. Let $S_E=S_{\bigO_E}\otimes_\Zp\Qp$,
and again extend the definitions of
$\Fil$, $\varphi$, $N$, and $\hat{g}$ to $S_{E}$ in the
obvious fashion.

Let $\MF(\varphi,N,K/K',E)$ be the category whose objects are finite free
$S_E$-modules $\D$ with
\begin{itemize}
\item a $\varphi$-semilinear, $E$-linear morphism $\varphi_\D:\D\to\D$ such
  that the determinant of $\varphi_D$ with respect to some choice of
  $S_\qp$ basis is invertible in $S_\qp$ (this does not depend on the choice of basis);
\item a decreasing filtration of $\D$ by $S_E$-submodules $\Fil^i\D$,
  $i\in\Z$, with $\Fil^i\D=\D$ for $i\ll 0$ and
  $(\Fil^iS_E)(\Fil^j\D)\subset\Fil^{i+j}\D$ for all $j$ and all $i\ge 0$;
\item a $K_0\otimes_{\Qp}E$-linear map $N:\D\to\D$ such that
  \begin{itemize}
  \item for all $s\in S_E$, $x\in\D$ we have $N(sx)=N(s)x+sN(x)$,
  \item $N\varphi=p\varphi N$,
  \item $N(\Fil^i\D)\subset\Fil^{i-1}\D$ for all $i$;
  \end{itemize}
\item an $S_E$-semilinear action of $\Gal(K/K')$ on $\D$ which
  commutes with $\varphi$ and $N$ and preserves each $\Fil^i\D$.
\end{itemize}

Fix a positive
integer $r\le p-2$. Then the category $\OEModdd$ of strongly
divisible $\cO_E$-modules with descent data is defined to be the
category of finitely generated free $S_{\cO_E}$-modules $\widehat{\mathcal{M}}$ with a
sub-$S_{\cO_E}$-module $\Fil^r\widehat{\mathcal{M}}$, additive maps $\varphi$, $N:\widehat{\mathcal{M}}\to\widehat{\mathcal{M}}$, and
$S_{\cO_E}$-semilinear bijections $\hat{g}:\widehat{\mathcal{M}}\to\widehat{\mathcal{M}}$ for each
$g\in\Gal(K/K')$ such that the following conditions hold.
\begin{itemize}
\item $\Fil^r\widehat{\mathcal{M}}$ contains $(\Fil^rS_{\cO_E})\widehat{\mathcal{M}}$,
\item $\Fil^r\widehat{\mathcal{M}}\cap I\widehat{\mathcal{M}}=I\Fil^r\widehat{\mathcal{M}}$ for all ideals $I$ of $\cO_E$,
\item $\varphi(sx)=\varphi(s)\varphi(x)$ for all $s\in S_{\cO_E}$, $x\in\widehat{\mathcal{M}}$,
\item $\varphi(\Fil^r\widehat{\mathcal{M}})$ is contained in $p^r\widehat{\mathcal{M}}$ and generates it over $S_{\cO_E}$,
\item $N(sx)=N(s)x+sN(x)$ for all $s\in S_{\cO_E}$ and $x\in\widehat{\mathcal{M}}$,
\item $N\varphi=p\varphi N$,
\item $E(u)N(\Fil^r\widehat{\mathcal{M}})\subset\Fil^r\widehat{\mathcal{M}}$,
\item for all $g\in\Gal(K/K')$, $\hat{g}$ commutes with $\varphi$ and
  $N$, and preserves $\Fil^r\widehat{\mathcal{M}}$,
\item $\hat{g}_1\circ\hat{g}_2=\widehat{g_1\circ g_2}$ for all $g_1$, $g_2\in\Gal(K/K')$.
\end{itemize}
For any $\widehat{\mathcal{M}}\in\OEModdd$, we define an $\cO_E$-module
$\Tst^{K'}(\widehat{\mathcal{M}})$ with an action of $G_{K'}$ as follows. 
Let $[\underline\varpi]$ be
the element of $\Acris$ corresponding to the Teichm\"uller representatives of a compatible system of $p$-th power roots of
$\varpi$. We let
$\hAst$ be the $S$-algebra with a filtration
$\Fil^i\hAst$, a Frobenius $\varphi$, and a monodromy operator $N$
defined in Section~2.2.2 of \cite{MR1681105}. This ring has a natural
action of $G_{K'}$ (see the second paragraph of Section~3.4 of
\cite{MR2137952} and the paragraph before Proposition~3.11 of
\cite{MR2137952}), and the natural map  $\hAst[1/p]\to\BdR$ is
$G_{K'}$-equivariant.  
[From Section~2.2.2 of \cite{MR1681105} we have $\hAst = \{\sum_{n=0}^\infty a_n \frac{X^n}{n!} : a_n \in \Acris,\ a_n
\to 0 \}$ and we define an action of $G_{K'}$ on $\hAst$ by setting $g(\sum a_n \frac{X^n}{n!}) = \sum g(a_n)
\frac{g(X)^n}{n!}$, where
\begin{equation}
  g(X) = \frac{g([\underline\varpi])}{[\underline\varpi]}\frac{\varpi}{g(\varpi)}(X+1)-1.\label{eq:15}
\end{equation}
(Since $g(\varpi) = \varpi$ for $g \in G_K$, this extends the natural action of $G_K$ on $\hAst$.) From the second
paragraph of Section~3.4 of \cite{MR2137952} and the paragraph before Proposition~3.11 of \emph{loc.\ cit.}\ we see that this
is the unique $G_{K'}$-action on $\hAst$ such that the natural map $\hAcris \into \hAst$ is $G_{K'}$-equivariant.
The second paragraph of Section~3.4 of \cite{MR2137952} shows that the map $f_\varpi:\hAcris\to\BdR^+$ is
$G_{K'}$-equivariant, and as explained in Section~2.2.2 of \cite{MR1681105} the map $\hAst[1/p]\to\BdR$ is determined
by the restriction of $f_\varpi$ to $\Acris$ and the fact that it sends the element $1+X\in\hAst$ to
$[\underline\varpi]/\varpi$. It follows that the map $\hAst[1/p]\to\BdR$ is $G_{K'}$-equivariant, as required.]
Then we
put \[\Tst^{K'}(\widehat{\mathcal{M}}):=\Hom_{\Fil^r,\varphi,N}(\widehat{\mathcal{M}},\hAst)\](that is, the
homomorphisms of $S$-modules which preserve $\Fil^r$ and commute with $\varphi$ and $N$). This inherits an $\cO_E$-module
structure from the $\cO_E$-module structure on $\widehat{\mathcal{M}}$, and $G_{K'}$ acts on
$\Tst^{K'}(\widehat{\mathcal{M}})$ by \[(gf)(x)=g(f(\hat{\bar{g}}^{-1}x))\] where $g\in G_{K'}$ and $\bar{g}$ is the image
of $g$ in $\Gal(K/K')$. Note that this action is well defined since the map $S \to \hAst$ is $G_{K'}$-equivariant
and since the $G_{K'}$-action on $\hAst$ commutes with
$\varphi$, $N$ and preserves $\Fil^i \hAst$ for all $i$. (The first is clear by the definition of the $G_{K'}$-actions,
since the map sends $u$ to $[\underline\varpi] (1+X)^{-1}$. The second follows from the definitions of these actions in
Section 2.2.2 of \cite{MR1681105}, and the fact that the $G_{K'}$-action commutes with $\varphi$ and preserves all $\Fil^i$ on 
$\hAcris$ by the second paragraph of Section~3.4 of \cite{MR2137952}.)

The main result of this subsection is the following.

\begin{prop}\label{prop: existence of weakly divisible with descent data}
  The functor $\Tst^{K'}$ provides an anti-equivalence of categories between the category $\OEModdd$ of strongly divisible
  $\cO_E$-modules with descent data and the category of $G_{K'}$-stable $\cO_E$-lattices in finite-dimensional
  $E$-representations of $G_{K'}$ which become semistable over $K$ with Hodge--Tate weights lying
in $[-r,0]$.
\end{prop}

When $K = K'$ and $E = \qp$ this is the main result of \cite{MR2388556}. 

\begin{proof}
  We write $\ZpMod$ for $\OEModdd$ and $\Tst$ for $\Tst^{K'}$ in the case when $K = K'$ and $E = \qp$.
  
  Suppose that $\wh\M$ is an object of $\OEModdd$. Let $\D := \wh\M[\frac 1p]$. We extend $\varphi$ and $N$ by linearity
  and we define a filtration on $\D$ as follows. We let $\Fil^r \D := (\Fil^r \wh\M)[\frac 1p]$ and
  \begin{equation*}
    \Fil^i \D := 
    \begin{cases}
      \D & \text{for $i \le 0$},\\
      \{ x \in \D : E(u)^{r-i} x \in \Fil^r \D \} & \text{for $0 \le i \le r$}, \\
      \sum_{j=0}^{i-1} (\Fil^{i-j} S_{\qp})(\Fil^j \D) & \text{for $i > r$, by induction.}
    \end{cases}
\end{equation*}
  It is not hard to check that $\D$ is an object of $\MF(\varphi,N,K/K',E)$. 
  (To verify the condition about $\varphi$, pick elements $f_i \in \Fil^r \M$ such that $\varphi_r(f_i)$ form an $S$-basis
  of $\M$, and let $e_i := \varphi_r(f_i)$. Write $f_i = \sum a_{ij} \varphi_r(f_j)$ and $\varphi(e_j) = \sum b_{jl} e_l$ for some
  matrices $A = (a_{ij})$, $B = (b_{jl})$ with entries in $S$. We have $p^r I = \varphi(A) B$, so $\det(B) \in S_\qp\s$.)
  By \cite[\S6]{MR1428871} and Corollary~2.1.4 in \cite{MR1681105} we can write $\D = D \otimes_{K_0} S_{\qp}$ for some
  weakly admissible filtered $(\varphi,N,K/K,\qp)$-module $D$, and we have $D \cong \D \otimes_{S_{\qp},s} K_0$ (where $s(h(u))
  = h(0)$) and $D_K \cong \D \otimes_{S_{\qp},s_\varpi} K$ (where $s_\varpi(h(u)) = h(\varpi)$), which induce $\varphi$, $N$ on $D$ and the
  filtration on $D_K$. It follows that $D$ inherits an $E$-action and a semilinear $\Gal(K/K')$-action and it is
  straightforward to verify that it becomes a weakly admissible filtered $(\varphi,N,K/K',E)$-module. Note that $\Fil^{r+1} D_K
  = 0$ since $\Fil^{r+1} \D \subset (\Fil^1 S_{\qp}) \D$.

  We associate an $E$-representation $\Vst^{K'}(\D)$ of
  $G_{K'}$ to $\D$ by defining \[\Vst^{K'}(\D):=\Hom_{\Fil^r,\varphi,N}(\D,\hAst[1/p])\](that
  is, the homomorphisms of $S$-modules which preserve $\Fil^r$ and
  commute with $\varphi$ and $N$). This inherits an $E$-module
  structure from the $E$-module structure on $\D$, and $G_{K'}$ acts
  on $\Vst^{K'}(\D)$ by \[(gf)(x)=g(f(\hat{\bar{g}}^{-1}x)),\] where
  $g\in G_{K'}$ and $\bar{g}$ is the image of $g$ in $\Gal(K/K')$.
  
  We have a natural $\cO_E[G_{K'}]$-linear map
  $\Tst^{K'}(\wh\M) \to \Vst^{K'}(\D)$.  We claim that there is a natural isomorphism of $E[G_{K'}]$-representations
  $\Vst^{K'}(\D)\to \Vst^{K'}(D)$ given by $f \mapsto f|_D$.
  In the case that $E=\Qp$ and $K'=K$, this is
  established in the proof of Proposition 2.2.5 of \cite{MR1971512}
  (which also shows that $f(D) \subset \Bst^+ \into \hAst[1/p]$), so
  in general this map gives a natural isomorphism $\Vst^{K'}(\D)\to
  \Vst^{K'}(D)$ of the underlying $\Qp$-vector spaces with an action
  of $G_{K}$. 
  It is clear that this
  isomorphism is $E$-linear, and it remains to check that it is
  compatible with the action of $G_{K'}$. 
  In the proof we just cited it is also shown that we have a commutative diagram
  \begin{equation*}
    \xymatrix@C=2.9pc{\D \ar[rr]^f\ar@{>>}[d]_{f_\varpi} && \hAst[1/p] \ar[d] \\ 
      D_K \ar[r]^<<<<<<{f|_D \otimes 1} & \Bst \otimes_{K_0}K \ar@{^{(}->}[r] & \BdR}
  \end{equation*}
  Since the maps
  $f_\varpi$, $\hAst[1/p]\to\BdR$, and $\Bst\otimes_{K_0}K\to\BdR$ are $G_{K'}$-linear, so is $\Vst^{K'}(\D)\to \Vst^{K'}(D)$.
  By Proposition 2.2.5 of \cite{MR1971512} we know that $\Tst^{K'}(\wh\M) \to \Vst^{K'}(\D)$ is a $G_K$-stable $\zp$-lattice. It
  follows that $\Tst^{K'}(\wh\M)$ is a $G_{K'}$-stable $\cO_E$-lattice in $\Vst^{K'}(\D)$, and $\Vst^{K'}(\D) \cong
  \Vst^{K'}(D)$ becomes semistable over $K$ with Hodge--Tate weights in $[-r,0]$ by Proposition~\ref{prop: D and V are equivs, savitt}.
  
  The faithfulness of $\Tst^{K'}$ follows immediately from the case when $E=\Qp$ and $K'=K$. Suppose now that $\wh\M_1$, $\wh\M_2$
  are objects of $\OEModdd$ and that we are given an $\cO_E[G_{K'}]$-linear map $\theta : \Tst^{K'}(\wh\M_1) \to
  \Tst^{K'}(\wh\M_2)$. From fullness in the case when $E=\Qp$ and $K'=K$ we get a map $\eta : \wh\M_1 \to \wh\M_2$ in $\ZpMod$
  such that $\Tst(\eta) = \theta$. We need to show that $\eta$ is $\cO_E$-linear and compatible with the
  $\Gal(K/K')$-actions.  If $\lambda \in \cO_E$, then $\Tst(\lambda\eta) = \Tst(\eta) \lambda = \lambda \Tst(\eta) =
  \Tst(\eta\lambda)$, so $\lambda\eta = \eta\lambda$, as required.  If $g \in G_{K'}$, then $\hat{\bar{g}}^{-1} \eta
  \hat{\bar{g}}$ is a map in $\ZpMod$ and it is straightforward to show that $\Tst(\eta)g = g\Tst(\eta)$
  is equivalent to $\Tst(\hat{\bar{g}}^{-1} \eta \hat{\bar{g}}) = \Tst(\eta)$, so $\hat{\bar{g}}^{-1} \eta
  \hat{\bar{g}} = \eta$, as required.

  It remains to check that $\Tst^{K'}$ is essentially surjective.
  Suppose that $V$ is a finite-dimensional $E$-representation of $G_{K'}$, which becomes semistable
  over $K$ with Hodge--Tate weights in $[-r,0]$, and suppose that $V_{\cO_E}$ is a $G_{K'}$-stable $\cO_E$-lattice in $V$.
  By Proposition \ref{prop: D and V are equivs, savitt} there is a
  weakly admissible
  filtered $(\varphi,N,K/K',E)$-module $D$ with $\Fil^0D_K=D_K$,
  $\Fil^{r+1}D_K=0$, and $\Vst^{K'}(D)\cong V$ as $E$-representations
  of $G_{K'}$. 

Let $\D:=S_E\otimes_{E}D$, a free $S_E$-module of finite rank. We
endow $\D$ with the structure of an object of $\MF(\varphi,N,K/K',E)$ as follows:
\begin{itemize}
\item $\varphi_\D$ is given by the tensor product of $\varphi$ on $S_E$ and
  $\varphi$ on $D$.
\item $N:=N\otimes\Id+\Id\otimes N:\D\to\D$.
\item $\Fil^i\D:=\D$ for $i\le 0$, and by
  induction for $i\ge 0$ \[\Fil^{i+1}\D:=\{x\in\D|N(x)\in\Fil^i\D\text{ and
  }f_\varpi(x)\in\Fil^{i+1}D_K\}\]where $f_\varpi:\D\to D_K$ is defined by
  $\lambda\otimes x\mapsto s(\lambda)x$, where $s:S_E\to K\otimes_\Qp E$ is the
  unique continuous $K_0\otimes_{\qp} E$-linear map sending $u^{ie}/i!$ to $\varpi^{ie}/i!$.
\item The action of $\Gal(K/K')$ is given by the tensor products of
  the actions on $S_E$ and $D$. 
\end{itemize}

By Theorem 2.3.5 of \cite{MR2388556},
there is a strongly divisible lattice $\widehat{\mathcal{M}}\subset\D$ of weight $r$ with
$\Tst(\widehat{\mathcal{M}})$ corresponding to $V_{\cO_E}$ under the isomorphism
$\Vst^{K'}(\D) \cong \Vst^{K'}(D) \cong V$. Note that $\wh\M \in \ZpMod$.
We now show that $\widehat{\mathcal{M}}$ is stable
under the actions of $\cO_E$ and $\Gal(K/K')$ on $\D$. Since $\cO_E$
is generated by $\cO_E^\times$ as a $\Z_p$-module, to check that $\widehat{\mathcal{M}}$
is $\cO_E$-stable it is enough to check that $x\widehat{\mathcal{M}}=\widehat{\mathcal{M}}$ for all
$x\in\O_E^\times$. Certainly $x^{-1}\widehat{\mathcal{M}}$ is a strongly divisible lattice in
$\D$. We also
have \[\Tst(\widehat{\mathcal{M}})=x\Tst(\widehat{\mathcal{M}})=\Tst(x^{-1}\widehat{\mathcal{M}}),\]so it
suffices to show that if $\widehat{\mathcal{M}}$ and $\widehat{\mathcal{M}}'$ are strongly divisible
lattices in $\D$ with $\Tst(\widehat{\mathcal{M}})=\Tst(\widehat{\mathcal{M}}')$, then
$\widehat{\mathcal{M}}=\widehat{\mathcal{M}}'$.
But this is a formal consequence of Theorem 2.3.5 of
\cite{MR2388556}.
The proof that $\widehat{\mathcal{M}}$ is
$\Gal(K/K')$-stable is very similar: if $g\in G_{K'}$ then
$\hat{\bar g}(\widehat{\mathcal{M}})$ is also a strongly divisible lattice in $\D$,
and \[\Tst(\widehat{\mathcal{M}})=g\Tst(\widehat{\mathcal{M}})=\Tst(\hat{\bar g}(\widehat{\mathcal{M}})),\]so
$\hat{\bar g}(\widehat{\mathcal{M}})=\widehat{\mathcal{M}}$, as required.

It remains to check that if $I$ is an ideal of $\cO_E$ then
$I\Fil^r\widehat{\mathcal{M}}=\Fil^r\widehat{\mathcal{M}}\cap I\widehat{\mathcal{M}}$, and that $\widehat{\mathcal{M}}$ is free as an
$S_{\cO_E}$-module. The first follows since $\Fil^r \wh \M = \wh \M \cap \Fil^r \D$, by definition. The second
follows as in the proof of Proposition
3.2.3.2 of \cite{MR1944572}: let $J$ denote the kernel of the ring homomorphism $S \onto W(k)$, $h(u)\mapsto h(0)$.
Then $\wh\M/J\wh\M \into D$ is a finitely generated $W(k) \otimes_\Zp \cO_E$-module which is a $\zp$-lattice. Since
$D$ is a finite free $K_0 \otimes_\Qp E$-module and $K_0/\qp$ is unramified one can check that $\wh\M/J\wh\M$ is
a finite free $W(k) \otimes_\Zp \cO_E$-module. (The point is that $W(k) \otimes_\Zp \cO_E$ is a product of discrete
valuation rings.) As in \emph{loc.\ cit.}\ we see that any lift of a basis to $\wh\M$ is a basis of $\wh\M$ as
$S_{\cO_E}$-module.
\end{proof}
\begin{remark}
  One of the referees has pointed out to us that it is possible to
  prove Proposition~\ref{prop: existence of weakly divisible with
    descent data} without using the formalism of $S_\Qp$-modules and
  how they relate to weakly admissible modules, as an essentially
formal consequence of the results of
\cite{MR2388556}.
However, it is
  convenient to have the relationship with $S_\Qp$-modules
 available in calculations, such as those of Lemma~\ref{lem: form of rank one objects
    and their generic fibers} below, and so we have let the proof
of the proposition stand in its current form,
as an illustration of that relationship.
\end{remark}

\subsection{Breuil modules with descent data: part one}
\label{subsec:breuil-modules-with-dd}
As its title indicated, in this subsection we introduce Breuil modules with descent data,
which morally speaking are the reductions mod $p$ of strongly divisible modules with
descent data, and which are the basic tool for studying the reductions mod $p$ 
of lattices in potentially semistable Galois representations.

We carry over the running hypotheses of the
preceding subsection, namely we consider a tamely ramified Galois extension $K/K'$,
with $K$ finite over $\mathbb Q_p$, and with ramification index $e(K/K')$,
and we furthermore assume given a uniformiser $\varpi \in K$ such that $\varpi^{e(K/K')} \in K'$.
We let $e$ denote the absolute ramification index of $K$.

Let $\F$ be a finite extension of $\Fp$.
Recall that the category $\FBrModdd$  of \emph{Breuil
modules of weight $r$ with descent data} from $K$ to $K'$ and
coefficients in $\F$  consists of quintuples
$(\mathcal{M},\mathcal{M}_{r},\varphi_{r},\hat{g},N)$ where:
\begin{itemize}\item $\mathcal{M}$ is a finitely generated
  $(k\otimes_{\F_p}\F)[u]/u^{ep}$-module, free over $k[u]/u^{ep}$.
\item $\M_{r}$ is a $(k\otimes_{\F_p}\F)[u]/u^{ep}$-submodule of $\M$
  containing $u^{er}\M$.
\item $\varphi_{r}:\M_{r}\to\M$ is $\F$-linear and $\varphi$-semilinear
  (where $\varphi:k[u]/u^{ep}\to k[u]/u^{ep}$ is the $p$-th power map)
  with image generating $\M$ as a
  $(k\otimes_{\F_p}\F)[u]/u^{ep}$-module.
\item $N:\M\to \M$ is $k\otimes_{\F_{p}}\F$-linear and satisfies
  $N(ux)=uN(x)-ux$ for all $x\in\M$,
  $u^{e}N(\M_{r})\subset\M_{r}$, and
  $\varphi_{r}(u^{e}N(x))=cN(\varphi_{r}(x))$ for all $x\in\M_{r}$. Here,
  $c \in (k[u]/u^{ep})^\times$ is the image
  of $\varphi_1(E(u))$ under the natural map $S \to k[u]/u^{ep}$.
\item $\hat{g}:\M\to\M$ are additive bijections for each
  $g\in\Gal(K/K')$, preserving $\M_{r}$, commuting with the $\varphi_{r}$-
  and $N$-actions, and satisfying $\hat{g}_1\circ
  \hat{g}_2=\widehat{g_1\circ g_2}$ for all $g_1, 
  g_2\in\Gal(K/K')$. Furthermore, if $a\in
  k\otimes_{\F_{p}}\F$, $m\in\M$ then
  $\hat{g}(au^{i}m)=g(a)((g(\varpi)/\varpi)^{i}\otimes
  1)u^{i}\hat{g}(m)$.
\end{itemize}

\begin{rem}
It follows from
  the assumption that  $\mathcal{M}$ is a finitely generated
  $(k\otimes_{\F_p}\F)[u]/u^{ep}$-module, which is free over
  $k[u]/u^{ep}$, that $\mc{M}$ is in fact a free
  $(k\otimes_{\F_p}\F)[u]/u^{ep}$-module. (The proof of Proposition 1.2 of \cite{MR2457845} goes
  through, with $\M_r$, $\varphi_r$ replacing $\M_1$, $\varphi_1$.)
\end{rem}

Suppose that $\F$ is the residue field of a finite extension $E$ of $\mathbb Q_p$.
If $\widehat{\mathcal{M}}$ is an object of $\OEModdd$, then
$\M:=(\widehat{\mathcal{M}}/\m_E\widehat{\mathcal{M}})\otimes_{S_{\cO_E}/\m_ES_{\cO_E}}(k\otimes_{\Fp}\F)[u]/u^{ep}$ is
naturally an object of $\FBrModdd$,
where we define $\M_r$ to be the image of $\Fil^r\widehat{\mathcal{M}}$ in $\M$,
the map $\varphi_r$ is induced by $(1/p^r)\varphi|_{\Fil^r\widehat{\mathcal{M}}}$, and $N$ and
$\hat{g}$ are those coming from $\widehat{\mathcal{M}}$. 
(To see that $\varphi_r$ is well defined on $\M$ comes down to the fact that $\ker(S/pS \to k[u]/u^{ep}) = \Fil^p(S/pS)$.)

We define a functor $\Tst$ from the category
$\FBrModdd$ to the category of finite-dimensional $\F$-representations of $G_{K'}$ as
follows: we let $\hA$ denote the $k[u]/u^{ep}$-algebra 
defined by the pushout (i.e.\ tensor product) square
\begin{equation*}
  \xymatrix{S\ar[r]\ar@{>>}[d]&\hAst \ar@{>>}[d]\\k[u]/u^{ep}\ar[r]&\hA}
\end{equation*}
The vertical arrows induce a filtration $\Fil^r$, maps $\varphi_r$ and $N$, and
an action of $G_{K'}$ on the bottom row. (See \S2.3 of \cite{MR2248152} or \S2.1 of \cite{carusomodp}.)
Explicitly, for $g \in G_{K'}$ we have $g(u) = \left(\frac{g(\varpi)}{\varpi}\right)^q u$ and
$g(X) = \left(\frac{g(\varpi_1)}{\varpi_1}\right)^{f} (1+X) - 1$, where $\varpi_1 \in \cO_{\o K}$ denotes
the $p$-th root of $\varpi$ used in \eqref{eq:15}, where $pq \equiv 1 \pmod {e(K/K')}$, and where $f \equiv 1 \pmod p$
is a multiple of $e(K/K')$. (See also \S4.2 in \cite{carusomodp}.)
For $\M \in \FBrModdd$ we define
\[\Tst(\M)=\Hom_{\operatorname{BrMod}}(\M,\hA),\] 
(that is, the $k[u]/u^{ep}$-homomorphisms which are compatible with the filtration, $\varphi_r$ and $N$). We let
$G_{K'}$ act by \[(gf)(x)=gf(\hat{\bar{g}}^{-1}x).\]

\begin{lm}\label{lm:basic-properties-of-Tst}
  The functor $\Tst$ is faithful. For any $\M \in \FBrModdd$ we have \[\dim_\F
  \Tst(\M)=\rank_{(k\otimes_{\F_p}\F)[u]/u^{ep}}\M.\] If $\wh\M \in \OEModdd$ and
  $\M := (\widehat{\mathcal{M}}/\m_E\widehat{\mathcal{M}})\otimes_{S_{\cO_E}/\m_ES_{\cO_E}}(k\otimes_{\Fp}\F)[u]/u^{ep}$
  denotes the Breuil module corresponding to the reduction of $\wh\M$, then
  \[\Tst^{K'}(\widehat{\mathcal{M}})\otimes_{\cO_E}\F\cong \Tst(\M).\]
\end{lm}

\begin{proof}
  The first two parts follow immediately from the corresponding results without coefficients and descent data
  (use Lemme 2.3.4 in \cite{MR2248152} together with Corollary 2.3.3 of \cite{MR2543474}, respectively
  Lemme 2.3.1.2 in~\cite{MR1681105}).

  For the last part, we first define $\varphi_r : \Fil^r \wh\M \to \M$ by $(1/p^r)\varphi$.
  Then for any ideal $I \subset \cO_E$, the $S/I$-module $\wh\M/I$ inherits a filtration
  (the image of $\Fil^r \wh\M$) and maps $\varphi_r$, $N$, $\hat g$ such that the natural map $\wh\M \to \wh\M/I$ is compatible
  with all structures. Similarly, $\hAst/p^n$ inherits from $\hAst$ a filtration, maps $\varphi_r$, $N$, and a $G_{K'}$-action.
  The natural map $\Tst^{K'}(\wh\M) \to \plim \Hom(\wh\M/p^n,\hAst/p^n)$ is an isomorphism, as explained
  just after D\'efinition~4.1.1.1 in \cite{MR1621389}.
  (Here and in the following, we consider $S$-linear homomorphisms that respect filtration, $\varphi_r$, $N$.) 
  We consider the following three $\cO_E[G_{K'}]$-linear maps:
  \begin{align}
    \label{eq:16}
    \Tst^{K'}(\wh\M)/p &\to \Hom(\wh\M/p,\hAst/p) \\
    \label{eq:17}
    \Hom(\wh\M/p,\hAst/p)/\mf{m}_E &\to \Hom(\wh\M/\m_E, \hAst/p) \\
    \label{eq:18}
    \Hom(\wh\M/\mf{m}_E, \hAst/p) &\to \Hom(\M, \hA) = \Tst(\M)
  \end{align}
  The first two are isomorphisms by the argument used in the proof of Proposition 2.3.2.4 in \cite{MR1681105}.
  The last is an isomorphism by Lemme 2.3.4 of \cite{MR2248152}.
\end{proof}

\begin{prop}\label{prop: unique quotient}
  Suppose that $\cN\in\FBrModdd$, and that $T'$ is a
  $G_{K'}$-subrepresentation of $\Tst(\cN)$ \(so that in particular
  $T'$ has the structure of an $\F$-vector space\). Then there is a
  unique quotient $\cN'$ of $\cN$ in $\FBrModdd$ such that if $f:\cN\to\cN'$ is the
  quotient map, then $\Tst(f)$ is identified with the inclusion
  $T'\into \Tst(\cN)$.
\end{prop}
\begin{proof}
  We write $\FpBrMod$ for the category $\FBrModdd$ in case $K'=K$ and $\F=\F_p$.
  In this case the proposition is proved in Proposition 2.2.5 of
  \cite{carusomodp}. Thus we see that in general there is a unique
  quotient $\cN'$ of $\cN$ in $\FpBrMod$
  such that if $f:\cN\onto\cN'$ is the
  quotient map, then $\Tst(f)$ is the inclusion
  $T'\into \Tst(\cN)$.

  It remains to verify that the $\Gal(K/K')$- and $\F$-actions descend to $\cN'$ via $f$.  (It is then automatic that $\cN'$
  becomes an object of $\FBrModdd$ and the result follows.)  Suppose for the moment that $\M \in \FpBrMod$ and that $g
  \in \Gal(K/K')$. We define the \emph{$g$-twist of $\M$} as follows: we let $\M\su[g] := \M \otimes_{k[u]/u^{ep}, g}
  k[u]/u^{ep}$. Then we have a $g$-semilinear bijection $i: \M \to \M\su[g]$ sending $x$ to $x \otimes 1$. There is a unique way
  to make $\M\su[g]$ into an object of $\FpBrMod$ such that the map $i$ induces a bijection on filtrations and commutes with
  $\varphi_r$ and $N$. (Note that $c$ is fixed by $\Gal(K/K')$.)
  It follows that for our $\cN\in\FBrModdd$ the $g$-semilinear bijection $\hat g : \cN \to \cN$ induces
  an isomorphism $\hat g: \cN\su[g] \congto \cN$ in $\FpBrMod$. The assumption that $T'$ is $G_{K'}$-stable implies that
  there is a unique $\Fp[G_K]$-linear isomorphism $\theta(g) : \Tst(\cN') \congto \Tst({\cN'}\su[g])$ such that
  $\Tst(\hat g) \circ \Tst(f) = \Tst(f\su[g]) \circ \theta(g)$. By the uniqueness of $\cN'$ we get a
  unique isomorphism $\hat g_{\cN'} : {\cN'}\su[g] \congto \cN'$ in $\FpBrMod$ such that $f \circ \hat g = \hat g_{\cN'} \circ f$.

  Finally, a similar but easier argument shows that $\cN'$ inherits a natural action of $\F^\times$ and thus (by
  $\Fp$-linearity) of $\F$, as required.
\end{proof}
\begin{rem}\label{rem: quotient is naive quotient}
 As already remarked in the proof, $\cN$ is actually a quotient in the naive sense
  that $\cN'=\cN/\cN''$ for some $\cN''\in\FBrModdd$ (this follows
  easily from the proof of Proposition 2.2.5 of
  \cite{carusomodp} and the exactness of the functor
  $M_{\tilde{\mathfrak{S}}}$, for which see Theorem 2.1.2 of
  \cite{carusomodp}; we thank Xavier Caruso for explaining this to us).
\end{rem}
It will be convenient for us to be able to employ duality in the
calculations that follow. To this end, we briefly recall the duality
theory for $\FBrModdd$ developed in chapter V of \cite{caruso-thesis} (see also Section~2.1 of \cite{carusomodp}).
\begin{defn}\label{defn: dual Breuil module}If $\M$ is
an object of $\FBrModdd$, then there is an object $\M^*$ of
$\FBrModdd$ such
that \[\Tst(\M^*)=\Tst(\M)^\vee(r),\](i.e.\ $\Tst(\M^*)$ is the
$r$-th Tate twist of the contragredient of $\Tst(\M)$). Explicitly,
$\M^*$ is defined as follows.
\begin{enumerate}
\item $\M^*:=\Hom_{k[u]/u^{ep}-\Mod}(\M,k[u]/u^{ep})$, with $\F$-structure
  inherited from that of $\M$.
\item $\M_r^*:=\{f\in\M^*|f(\M_r)\subset u^{er}k[u]/u^{ep}\}$.
\item Let $\varphi_r:u^{er}k[u]/u^{ep}\to k[u]/u^{ep}$ be the unique
  semilinear map sending $u^{er}$ to $c^r$. Then we define $\varphi_r(f)$
  by $\varphi_r(f)(\varphi_r(x))=\varphi_r(f(x))$ for all $x\in\M_r$,
  $f\in\M_r^*$.
\item $N(f):=N\circ f-f\circ N$, where $N : k[u]/u^{ep} \to k[u]/u^{ep}$ is the unique
  $k$-linear derivation such that $N(u) = -u$.
\item Let $g\in\Gal(K/K')$ act on $\M^*$ by
  $(\hat{g}f)(x)=g(f(\hat{g}^{-1}x))$, where $\Gal(K/K')$ acts on
  $k[u]/u^{ep}$ by $g(au^i)=g(a)(g(\varpi)/\varpi)^iu^i$ for $a \in k$.
\end{enumerate}
\end{defn}
(A priori we have a perfect, $\F_p$-bilinear pairing $\ang{\cdot,\cdot} : \Tst(\M) \times \Tst(\M^*) \to \fp(r)$ that is
$G_K$-equivariant. Explicitly, $\ang{f,f'} = \sum f(e_i) f'(e_i^*)$, where $(e_i)$ is any $k[u]/u^{ep}$-basis of $\M$ 
and $(e_i^*)$ is the dual basis of $\M^*$, by the proof of Theorem V.4.3.1 of \cite{caruso-thesis}. It follows that
$\ang{\cdot,\cdot}$ is $G_{K'}$-equivariant and that $\ang{\lambda f,f'} = \ang{f,\lambda f'}$ for all $\lambda \in \F$.)
We have $\M^{**}\cong\M$. We define the covariant functor $\Tst^{*,r}$ from $\FBrModdd$ to
the category of $\F$-representations of $G_{K'}$
by \[\Tst^{*,r}(\M):=\Tst(\M)^\vee(r),\]so
that \[\Tst^{*,r}(\M)=\Tst(\M^*).\] We then have the following
immediate corollary of Proposition \ref{prop: unique quotient}.
\begin{coroll}
  \label{cor:unique-sub} Suppose that $\cM\in\FBrModdd$, and that $T'$ is a
  $G_{K'}$-subrepresentation of $\Tst^{*,r}(\cM)$ \(so that in particular
  $T'$ has the structure of an $\F$-vector space\). Then there is a
  unique subobject $\cM'$ of $\cM$ in $\FBrModdd$ such that if $f:\cM'\to\cM$ is the
  inclusion map, then $\Tst^{*,r}(f)$ is identified with the inclusion
  $T'\into \Tst^{*,r}(\cM)$.
\end{coroll}
\begin{proof}Let $\cN:=\cM^*$, so that $T'$ is a
  sub-$G_{K'}$-representation of $\Tst(\cN)=\Tst^{*,r}(\cM)$. Then
  by Proposition \ref{prop: unique quotient}, there is a unique
  quotient $\cN'$ of $\cN$ such that applying $\Tst$ to the natural
  projection map $\cN\to\cN'$ gives the inclusion $T'\into
  \Tst(\cN)$. Then $\M':=(\cN')^*$ is the required subobject of $\M$,
  and uniqueness is clear.
  \end{proof}
  \begin{rem}
    Just as in Remark \ref{rem: quotient is naive quotient}, $\M'$ is
    a subobject of $\M$ in the naive sense that it is a
    sub-$(k[u]/u^{ep}\otimes\F)$-module of $\M$, which inherits the
    structure of an object of $\FBrModdd$ from $\M$ in the obvious
    way.
  \end{rem}
Finally, we recall something of the theory of maximal objects from
\cite{carusomodp} (see Sections~4.1 and~4.2 of \cite{carusomodp} for a
proof of the existence of maximal Breuil modules with descent data and
coefficients). 
If $\M$ is an object of $\FBrModdd$, the
corresponding maximal object $\Max(\M)$ is an object of
$\FBrModdd$ together with a morphism
$\iota_{\Max}^\M:\M\to\Max(\M)$ such that
$\Tst^{*,r}(\iota_{\Max}^\M)$ is an isomorphism. Furthermore, it
satisfies the following universal property: if $f:\M\to\M'$ is such
that $\Tst^{*,r}(f)$ is an isomorphism, then there exists a unique
morphism $g:\M'\to\Max(\M)$ such that $g\circ f=\iota_{\Max}^\M$.

\subsection{Breuil modules with descent data: part two}
\label{subsec:breuil-modules-two}
We now specialise to the particular situation of interest to us in
this paper. Let $K_0$ be the unique unramified extension of $\Qp$ of
degree $d$, and let $K=K_0((-p)^{1/(p^d-1)})$. Let $K'=K_0$. Fix
$\varpi=(-p)^{1/(p^d-1)}$. We write $\wt\omega_d:\Gal(K/K_0)\to K_0^\times$ for
the character $g\mapsto g(\varpi)/\varpi$, and we let $\omega_d$ be the
reduction of $\wt\omega_d$ modulo $\varpi$. (By inflation we can also think of
$\wt\omega_d$ and $\omega_d$ as characters of $I_{K_0} = I_\qp$.
Note that $\omega_d$ is a tame fundamental character
of niveau~$d$ and that $\wt\omega_d$ is the Teichm\"uller lift of $\omega_d$.)

Let $E/\Qp$ be a finite extension such that $E$ contains the images of
all embeddings $K\into\Qpbar$.

If $\rho:G_{K_0}\to\GL_n(E)$ is a
potentially semistable representation which becomes semistable over
$K$, then the associated inertial type is a representation of
$I_{K_0}$ which becomes trivial when restricted to $I_K$, so we may
think of it as a representation of $\Gal(K/K_0)\cong I_{K_0}/I_K$. Since this group is abelian, this representation is
isomorphic to a direct sum $\chi_1\oplus\dots\oplus\chi_n$ where each
$\chi_i:\Gal(K/K_0)\to \cO_E^\times$ is a character.

Let $\F$ be the residue field of $\cO_E$.
Let $\phi$ be the arithmetic Frobenius on $k$, and let
$\sigma_0:k\into\F$ be a fixed embedding. Inductively define
$\sigma_1,\dots,\sigma_{d-1}$ by
$\sigma_{i+1}=\sigma_i\circ\phi^{-1}$; we will often consider the
numbering to be cyclic, so that $\sigma_d=\sigma_0$.
There are idempotents $e_i\in
k\otimes_\Fp\F$ such that if $M$ is any $k\otimes_\Fp\F$-module,
then $M=\bigoplus_ie_iM$, and $e_iM$ is the subset of $M$ consisting of
elements $m$ for which $(x\otimes 1)m=(1\otimes\sigma_i(x))m$ for all
$x\in k$. Note that $(\phi\otimes 1)(e_i) = e_{i+1}$ for all $i$.

\begin{prop}\label{prop: form of descent data on breuil module}
  Maintain our current assumptions on $K_0$, $K$ and $K'$, so that
  $K'=K_0$ and $K=K_0((-p)^{1/(p^d-1)}$. Suppose
  that $\rho:G_{K_0}\to\GL_n(E)$ is a continuous representation such
  that $\rho|_{G_K}$ is semistable with Hodge--Tate weights contained
  in $[-r,0]$, with $r\le p-2$, and that $\rho$ has inertial type
  $\chi_1\oplus\dots\oplus\chi_n$ as above. Let $\rhobar$ be the
  reduction mod $\m_E$ of a $G_{K_0}$-stable $\cO_E$-lattice in
  $\rho$. Then there is an element $\M$ of $\FBrModdd$
  with \[\Tst^{*,r}(\M)\cong\rhobar,\] and $\M$ has a
  $(k[u]/u^{ep}\otimes_\Fp\F)$-basis $v_1,\dots,v_n$ such that
  $\hat{g}(v_i)=(1\otimes\o\chi_i(g))v_i$ for all
  $g\in\Gal(K/K_0)$.
\end{prop}
\begin{proof}
  By Proposition
  \ref{prop: existence of weakly divisible with descent data}, there
  is an element $\widehat{\mathcal{M}}$ of $\OEModdd$ such that
  $\Tst^{K'}(\widehat{\mathcal{M}})$ is the given lattice in $\rho$. Then we let
  $\o{\cN}:=(\widehat{\mathcal{M}}/\m_E\widehat{\mathcal{M}})\otimes_{S/pS}k[u]/u^{ep}$, and
  $\M:=\o{\cN}^*$, so that $\Tst^{*,r}(\M)\cong\rhobar$. In the
  case $r=1$, the claim about the form of the descent data is proved
  in Section~5 of \cite{gee-savitt-quaternion}; but the proof given there makes no use of the fact
  that $r=1$, and extends immediately to this more general situation.
\end{proof}

\begin{lem}\label{lem: form of rank one objects and their generic fibers}
  Maintain our current assumptions on $K_0$, $K$ and $K'$, so $e = p^d-1$. Then every
  rank one object of $\FBrModdd$ may be written in the form:
  \begin{itemize}
  \item $\M=((k\otimes_\Fp\F)[u]/u^{ep})\cdot m$,
  \item $e_i\M_r=u^{r_i}e_i\M$,
  \item $\varphi_r(\sum_{i=0}^{d-1}u^{r_i}e_im)=\lambda m$ for some
    $\lambda\in (k\otimes_\Fp\F)^\times$,
  \item $\hat{g}(m)=(\sum_{i=0}^{d-1}(\omega_d(g)^{k_i}\otimes
    1)e_i)m$ for all $g\in\Gal(K/K_0)$, and
  \item $N(m) = 0$.
  \end{itemize}
Here the integers $0\le r_i\le er$ and $k_i$ satisfy $k_i\equiv p(k_{i-1}+r_{i-1})\pmod{e}$ for all
$i$. Conversely, any module $\M$ of this form is a rank one object of $\FBrModdd$. Furthermore, 
\[\Tst^{*,r}(\M)|_{I_{K_0}}\cong\sigma_0\circ\omega_d^{\kappa_0},\]
where $\kappa_0\equiv k_0+s_0 \pmod e$, and \[s_0:=p(r_0p^{d-1}+r_{1}p^{d-2}+\dots+r_{d-1})/(p^d-1).\]
\end{lem}

\begin{proof}
  If $\M$ is a rank one object of $\FBrModdd$, then it is of the claimed form
  by Theorem 3.5 of \cite{MR2457845}, except for the statement that $N(m) = 0$.
  (Note that the cited proof only uses that $u^e\M \subset \M_1$ at one point, once all the work is done:
  to know that $r_i \le e$ in the statement of the theorem. So we
  can just replace $\M_1$ by $\M_r$ throughout.) Now it is easy to see that $N(\M_r) \subset \M_r$.
  Then $N(m) = 0$ by using that $\varphi_{r}(u^{e}N(x))=cN(\varphi_{r}(x))$ for all $x\in\M_{r}$.

  Conversely it is straightforward to verify that any module $\M$ of the above form
  is a rank one object of $\FBrModdd$.
  The calculation of 
  $\Tst^{*,r}(\M)|_{I_{K_0}}$ may be performed in the same way as in Example 3.7 of
  \cite{MR2457845} or Corollaire 4.4.2 of \cite{carusomodp}; we sketch
  the details of an approach following Example 3.7 of
  \cite{MR2457845}. 

Given $\M$ as in the statement of the lemma, define another element
$\M'$ of $\FBrModdd$ by  \begin{itemize}
  \item $\M'=((k\otimes_\Fp\F)[u]/u^{ep})\cdot m'$,
  \item $e_i\M'_r=e_i\M'$,
  \item $\varphi_r(\sum_{i=0}^{d-1}e_im')=\lambda m'$,
  \item $\hat{g}(m')=(\sum_{i=0}^{d-1}(\omega_d(g)^{p^i\kappa_0}\otimes
    1)e_i)m'$ for all $g\in\Gal(K/K_0)$, and
  \item $N(m') = 0$.

  \end{itemize}
We can define a morphism $\M'\to\M$ by mapping $e_im'\mapsto
u^{s_i}e_im$,  where \[s_i:=p(r_ip^{d-1}+r_{i+1}p^{d-2}+\dots+r_{i+d-1})/(p^d-1).\] One easily checks that this is a morphism in the
category $\FBrModdd$; for example, $s_i\ge r_i$, so the map is
compatible with the filtration, and since $s_{i+1}=p(s_i-r_i)$ the
morphism commutes with $\varphi_r$. Since $\Tst^{*,r}$ is faithful
we have a nonzero map $\Tst^{*,r}(\M')\to \Tst^{*,r}(\M)$, which must
be an isomorphism because both source and target are one-dimensional
$\F$-vector spaces.

The calculation of $\Tst^{*,r}(\M')$ may now be carried out in the same
fashion as in Example 3.7 of \cite{MR2457845}. We sketch the
argument. Denote a Teichm\"uller lift by a tilde. One checks without difficulty that $\M'$ may be lifted to a
strongly divisible module with descent data, contained inside the object of $\MF(\varphi,N,K/K',W(\F)[1/p])$ associated to the weakly
admissible filtered $(\varphi,N,K/K',W(\F)[1/p])$-module
$D=(K_0\otimes_\Qp W(\F)[1/p])\cdot v$, where $N=0$, $\varphi(v)=\wt{\lambda}p^rv$,
 ${g}(v)=(1 \otimes \wt{\sigma}_0\circ{\wt\omega_d}(g)^{\kappa_0})v$, and
  $\Fil^i(K\otimes_{K'} D)$ is $0$ for $i\ge r+1$ and $K\otimes_{K'}
  D$ for $i\le r$. The corresponding Galois representation can then
  be computed as in Example 2.13 of \cite{MR2137952}: we find $\Vst^{K'}(D)\dual(r)|_{I_{K_0}} 
  \cong \wt{\sigma}_0\circ{\wt\omega_d}^{\kappa_0}$,
  which indeed reduces to $\sigma_0\circ\omega_d^{\kappa_0}$.
\end{proof}

\begin{lm}\label{lm:maximal-rank-one-breuil-modules} Maintain our assumptions on $K_0$, $K$ and $K'$, so $e = p^d -1$.
  Suppose that $\M$ is the rank one object of $\FBrModdd$ of Lemma~\ref{lem: form of rank one objects and their generic
    fibers}. Then there exists a basis $\wt m$ of $\Max(\M)$ such
  that
  \begin{itemize}
  \item $e_i\Max(\M)_r=u^{er}e_i\Max(\M)$,
  \item $\varphi_r(\sum_{i=0}^{d-1} u^{er}e_i \wt m)=\lambda \wt m$,
  \item $\hat{g}(\wt m)=\sum_{i=0}^{d-1} (\omega_d(g)^{p^i s} \otimes 1)e_i \wt m$ for all $g\in\Gal(K/K_0)$, and
  \item $N(\wt m) = 0$,
  \end{itemize}
  with $s \equiv \kappa_0 - \frac{er}{p-1} \pmod e$, where we keep the notation of Lemma~\ref{lem: form of rank one objects
    and their generic fibers}.
\end{lm}

In particular,
the Breuil module in Lemma~\ref{lem: form of rank one objects and their
  generic fibers} is maximal if and only if $r_i = er$ for all $i$.

\begin{proof}
  Let $\M'$ denote the Breuil module defined in the statement of the lemma. Let
  \begin{equation*}
    t_i = \frac{erp}{p-1} - \frac{r_i p^d + r_{i+1} p^{d-1} + \cdots + r_{i+d-1} p}{e}.
  \end{equation*}
  It is an integer, and since $0 \le r_i \le er$, it lies in the interval
  $\big[0,\frac{erp}{p-1}\big]$. Define a map $f : \M \to \M'$ by $f(e_i m) := u^{t_i} e_i \wt m$. Then $f \ne 0$ since $r <
  p-1$.  To see that $f$ commutes with $\varphi_r$, use that $t_{i+1} = p(r_i+t_i-er)$. Since $p$ divides $t_i$, it
  follows that $f$ commutes with $N$. It is straightforward to check that $f$ commutes with descent data. Since
  $\Tst^{*,r}$ is faithful, it follows that $\Tst^{*,r}(f)$ is an isomorphism.
 
  It remains to verify that $\M'$ is maximal. If we have a map of Breuil modules $h: \M' \to \cN$ such that
  $\Tst^{*,r}(h)$ is an isomorphism, $\cN$ has to be of rank one and we can find a basis $n$ such that
  \begin{itemize}
  \item $e_i \cN_r=u^{r_i'}e_i \cN$,
  \item $\varphi_r(\sum_{i=0}^{d-1} u^{r_i'}e_i n)=\mu n$,
  \item $\hat{g}(n)=\sum_{i=0}^{d-1} (\omega_d(g)^{k_i'}\otimes 1) e_i n$ for all $g\in\Gal(K/K_0)$, and
  \item $N(n) = 0$,
  \end{itemize}
  for some $\mu \in (k\otimes\F)\s$ and integers $r_i' \in [0,er]$ and $k_i'$ such that
  $k_i' \equiv p(k_{i-1}'+r_{i-1}')$ modulo $e$. We can write $h(e_i \wt m) = \nu_i (u^{m_i} + {\mathrm O}(u^{m_i+1})) e_i n$
  for some $\nu_i \in \F\s$ and $m_i \ge 0$. Since $h$ commutes with $\varphi_r$ we find
  that
  \begin{equation}\label{eq:5a}
    m_{i+1} = p(m_i+er-r_i'),
  \end{equation}
  unless both sides are at least $ep$. But $h \ne 0$ implies that
  $m_i < ep$ for some $i$ and if $m_i < C$ for some constant $C \le ep$, then~\eqref{eq:5a} shows that
  $m_{i-1} < C/p$. Inductively we thus find that $m_i = 0$ for all $i$, so $r_i' = er$ for all $i$ and $h$ is an isomorphism.
\end{proof}

\begin{rk}
  Similarly (or by duality) we find that $\M$ in Lemma~\ref{lem: form of rank one objects and their generic fibers} is minimal if
  and only if $r_i = 0$ for all $i$. 
\end{rk}

The following lemma works more generally than for our particular choice of extension $K/K'$ above,
and so just for the duration of its statement and proof, we suppress our running definitions
of $K$ and $K'$, and return to the more general context
of Subsection~\ref{subsec:breuil-modules-with-dd}. 

\begin{lm}\label{lm:special-basis-of-breuil-module}
  Suppose that $K/K'$ is totally ramified and that $\F$ contains an
  embedding of $k$.
  Suppose that $\M$ is a rank $d$ object
  of $\FBrModdd$. Then $\M$ has a basis $m\su$ \($1 \le j \le d$\) such that for all $j$ we have $m\su \in \varphi_r(\M_r)$
  and $(k \otimes_\Fp \F)\cdot m\su$ is $\Gal(K/K')$-stable.
\end{lm}

\begin{proof}
We let $d = [k:\fp]$. Since we are assuming that $\F$ contains an embedding
of $k$, we may define the idempotents $e_i$ as above.
To prove the lemma,
it suffices to show that $e_0 \M$ has an $\F[u]/u^{ep}$-basis $n\su$ such that $\F\cdot n\su$ is $\Gal(K/K')$-stable.
  Thus it is enough to show that there are generators $\alpha\su$ ($1 \le j \le d$) of $e_{d-1} \M_r$
  as an $\F[u]/u^{ep}$-module such that $\F \cdot \alpha\su$ is $\Gal(K/K')$-stable. (Note that $(\varphi_r(\alpha\su))_{j=1}^d$
  is then a basis of $e_0 \M$ as $\F[u]/u^{ep}$-module.)
  
  Note that $\Gal(K/K')$ is abelian of order $e(K/K')$, which is prime to $p$. It acts trivially on $k$. Moreover $k$
  contains the $e(K/K')$-th roots of unity.  Thus $e_{d-1}\M_r$ decomposes as a direct sum of one-dimensional
  $\F$-subspaces that are $\Gal(K/K')$-stable. Pick an $\F$-basis $\beta_l$ of $e_{d-1}\M_r$
  such that $\F\cdot \beta_l$ is $\Gal(K/K')$-stable for all $l$. Pick a subset $\alpha\su$ ($1 \le j \le d$) of that basis
  that forms a basis inside $e_{d-1}(\M_r/u\M_r) \cong \F^d$. (The latter isomorphism follows from the theorem on elementary
  divisors, since $r < p$.) By Nakayama's lemma the $\alpha\su$ generate $e_{d-1} \M_r$, as required.
\end{proof}

\begin{lm}\label{lm:elem-divisors}
  Suppose that $M$ and $M'$ are free $\F[u]/u^{ep}$-modules of rank~$d$ and that $f : M \to M'$ is an $\F[u]/u^{ep}$-linear
  map. Let  $(e_j)_{j=1}^d$ be a basis of $M$. If $f(e_1) \in u^{n+1}M'$ for some positive integer $n+1 \le ep$, then
  $\im(f) \not\supset u^n M'$.
\end{lm}

\begin{proof}
  Suppose to the contrary that we have the containment $\im(f) \supset u^n M'$. 
Let $I$ denote the $\F$-vector space $\im(f)/u\im(f)$.  On the one hand, the images
  of the $f(e_j)$ span $I$, and by assumption $f(e_1)$ is zero in $I$.  Thus $\dim_\F I < d$. On the other hand, by the
  theorem on elementary divisors we can pick a basis $(e_j')_{j=1}^d$ of $M'$ and non-negative integers $a_j$ (possibly
  greater than or equal to $ep$) such that the $(u^{a_j}
  e_j')_{j=1}^{d}$ span $\im(f)$ as an $\F[u]/u^{ep}$-module. Since $u^n M'$ is
  spanned by the $(u^{n} e_j')_{j=1}^d$, we see that $a_j \le n < ep$ for all $j$. Therefore $\dim_\F I = d$ and we get a
  contradiction, as required.
\end{proof}

We owe the following argument to Xavier Caruso.  We would like to thank him for letting us include it here. Without it, we
would need a slightly stronger condition on $(a,b,c)$ in Theorem~\ref{thm: explicit list of possible characters for specific
  descent data} below.  (Namely we would need to demand that $a-b > 3$, $b-c > 3$, $a-c < p-4$ and that none of $2a-b-c$,
$2b-a-c$, $2c-a-b$ are congruent to an element of $[-3,3]$ modulo $p-1$.) In the statement, note that, as $\gcd(p,e) = 1$, we have
a natural ring homomorphism $\Z[1/p] \to \Z/e$ (``reduction modulo $e$'').

\begin{prop}\label{prop:descent-data-on-rank-one-subobject}
  Maintain our assumptions on $K_0$, $K$ and $K'$, so $e = p^d -1$. Assume that $r < \frac{p-1}2$.
  Suppose that $\M$ is a rank $d$ object of
  $\FBrModdd$. Suppose that $\Tst^{*,r}(\M)$ extends to an irreducible representation of $G_\Qp$. Suppose that $\cN$ is
  a rank one object of $\FBrModdd$ and that we are given a map
  $\cN \to \M$ of Breuil modules that is injective \(as an
  $\F$-linear map\).

  Suppose that $\M$ has a basis $m\su$ \($0 \le j \le d-1$\) with $\hat g(m\su) = (\omega_d(g)^{a_j}\otimes 1) m\su$ for all $j$ and all
  $g \in \Gal(K/K_0)$, where the integers $a_j$ satisfy the condition that $(a_j-a_l)/p$ is congruent to an element of
  the interval $\big(\frac{er}{p-1}, \frac{e(p-1-r)}{p-1}\big)$ modulo $e$ whenever $j \ne l$.
  Then there exists a basis $n$ of $\cN$ and a permutation $\pi$ of $\{0,1,\dots,d-1\}$ such that
  \[ \hat g(n) = \sum_{i=0}^{d-1} (\omega_d(g)^{a_{\pi(i)}}\otimes 1) e_i n. \]
\end{prop}

\begin{proof}
  Since $\cN \to \M$ is nonzero and $\Tst^{*,r}$ is faithful we see that $\Tst^{*,r}(\cN)$ is a one-dimensional
  subrepresentation $\F(\chi)$ of $\Tst^{*,r}(\M)$.  As $\Tst^{*,r}(\M)$ extends to an irreducible representation of
  $G_\Qp$ we have $\Tst^{*,r}(\M) \cong \bigoplus_{j=0}^{d-1} \F(\chi^{p^j})$ and the characters $\chi^{p^j}$ ($0 \le j \le
  d-1$) are pairwise distinct. By Lemma~\ref{cor:unique-sub} there are rank one Breuil modules $\cN\su$ for $0 \le j \le d-1$
  together with maps $\cN\su \to \M$ that are injective as $\F$-linear maps and such that $\Tst^{*,r}(\cN\su) \cong
  \F(\chi^{p^j})$.  We can take $\cN\su[0] = \cN$. It follows that we have a map $\alpha : \bigoplus_j \cN\su \to \M$ such
  that $\Tst^{*,r}(\alpha)$ is an isomorphism. By the universal
  property of maximal objects
  and Lemma~1.5.1 of \cite{carusomodp} we have a map $\beta : \M \to
  \bigoplus_j \Max(\cN\su)$ such that $\beta \circ \alpha$ is a direct sum of the natural maps $\cN\su \to \Max(\cN\su)$.

  By Lemma~\ref{lem: form of rank one objects and their generic fibers} there is a basis $n\su$ of $\cN\su$, integers $r_i\su
  \in [0,er]$ and $x_i\su \in \{a_1,\dots,a_d\}$, and $\lambda\su \in (k \otimes_\Fp \F)\s$ such that
  \begin{itemize}
  \item $e_i\cN\su_r=u^{r\su_i}e_i\cN\su$,
  \item $\varphi_r(\sum_{i=0}^{d-1} u^{r\su_i}e_i n\su)=\lambda\su n\su$,
  \item $\hat{g}(n\su)=\sum_{i=0}^{d-1} (\omega_d(g)^{x\su_i}\otimes 1) e_i n\su$ for all $g\in\Gal(K/K_0)$, and
  \item $N(n\su) = 0$.
  \end{itemize}
  (For the descent data note that $e_i(\cN/u\cN)$ injects into $e_i(\M/u\M)$.) By
  Lemma~\ref{lm:maximal-rank-one-breuil-modules} we know that there is a basis $\wt n\su$ of $\Max(\cN\su)$ and an integer
  $s$ such that
  \begin{itemize}
  \item $e_i\Max(\cN\su)_r=u^{er}e_i\Max(\cN\su)$,
  \item $\varphi_r(\sum_{i=0}^{d-1} u^{er}e_i \wt n\su)=\lambda\su \wt n\su$,
  \item $\hat{g}(\wt n\su)=\sum_{i=0}^{d-1} (\omega_d(g)^{p^{i+j}s}\otimes 1) e_i \wt n\su$ for all $g\in\Gal(K/K_0)$, and
  \item $N(\wt n\su) = 0$.
  \end{itemize}
  (Here we use that $\Tst^{*,r}(\cN\su) \cong \F(\chi^{p^j})$.) The same lemma shows that
  \begin{equation}\label{eq:1a}
    p^{i+j}s \equiv x_i\su - \frac{er}{p-1} + \frac{r_i\su p^d + r_{i+1}\su p^{d-1} + \cdots + r_{i+d-1}\su p}{e} \pmod e.
  \end{equation}
  Recall from the proof of Lemma~\ref{lm:maximal-rank-one-breuil-modules} that the natural map 
$$\cN\su \to \Max(\cN\su)$$
  sends $n\su$ to $\sum_{i} u^{t_i\su} e_i \wt n\su$ for certain $t_i\su \in [0,\frac{erp}{p-1}]$. In particular, it follows
  that
  \begin{equation}\label{eq:2a}
    \im(\beta\circ \alpha) \supset \bigoplus_{j=0}^{d-1} u^{erp/(p-1)} \Max(\cN\su).
  \end{equation}

  By Lemma~\ref{lm:special-basis-of-breuil-module} we may assume that $m\su \in \varphi_r(\M_r)$. (Note that we may arrange that
  $\hat g(m\su) = (\omega_d(g)^{a_j}\otimes 1) m\su$: for each $i$ we have $e_i(\M/u\M) \cong \bigoplus_j \F(\omega_d(g)^{a_j})$ as
  $\F[\Gal(K/K_0)]$-module by our assumption on the descent data of $\M$, so we may permute the $e_i m\su$ such that $\hat
  g(e_i m\su) = \omega_d(g)^{a_j} e_i m\su$ for all $j$.)
  It follows that
$$\beta(m\su) \in \bigoplus_l \varphi_r((\Max \cN\su[l])_r) =
  \bigoplus_l (k\otimes_\Fp \F)[u^p]/u^{ep} \cdot \wt n\su[l].$$
  Since $\beta$ commutes with descent data, it furthermore follows that 
$$\beta(e_0
  m\su) = \sum_{l=0}^{d-1} \gamma_{jl} u^{b_{jl}} e_0 \wt n\su[l],$$
where $\gamma_{jl} \in \F$ and $0 \le b_{jl} < ep$ is
  uniquely determined by the congruences
  \begin{equation}\label{eq:4a}
    \begin{split}
      b_{jl} &\equiv 0 \pmod p,  \\
      b_{jl} &\equiv a_j - p^l s \pmod e. 
    \end{split}
  \end{equation}

  It will be convenient to assume in the following that the $a_j$, and hence the $x_i\su$, are divisible by $p$. (Note that so far that
  only their values modulo $e$ mattered and that $(e,p) = 1$.)  Since $r_i\su \in [0,er]$ it follows from \eqref{eq:1a} that
  \begin{equation}
    \label{eq:3a}
    \text{$p^{i+j-1}s$ is congruent to an element of $\bigg[\frac{x_i\su}{p}-\frac{er}{p-1},\frac{x_i\su}{p}\bigg]$ modulo $e$.}
  \end{equation}
  Our assumption on the $a_j$ implies that two such intervals are either equal or disjoint.
  In particular, $x_i\su$ only depends on $i+j$ modulo $d$.

  Suppose now that the $(x_i\su[0])_{i=0}^{d-1}$ are not all
  distinct. (This is the descent data on $\cN$.) Then the
  $(x_0\su[l])_{l=0}^{d-1}$ are not all distinct. So there is $j$ such that $x_0\su[l] \ne a_j$ for all $l$.
  From equations \eqref{eq:4a} and \eqref{eq:3a} we find that $\frac{b_{jl}}p$ is congruent to an element of 
  \begin{equation*}
    \bigg[\frac{a_j-x_0\su[l-1]}{p}, \frac{a_j-x_0\su[l-1]}{p}+\frac{er}{p-1}\bigg]
  \end{equation*} 
  modulo $e$.
  Our assumption on the $a_j$ then implies that $b_{jl} > \frac{erp}{p-1}$ for all~$l$. 
  But, by applying Lemma~\ref{lm:elem-divisors} to the map $\beta : e_0 \M \to \bigoplus_j e_0 \Max(\cN\su)$ 
  and using that $r < p-1$, we see that this contradicts \eqref{eq:2a}.
\end{proof}

\begin{thm}\label{thm: explicit list of possible characters for
    specific descent data}
  Let $\rho:G_{\Qp}\to\GL_3(\Qpbar)$ be a potentially semistable representation with Hodge--Tate weights $-2$, $-1$,
  and $0$ such that $\rhobar$ is irreducible.

  In the following we assume that $a$, $b$, $c$ are integers satisfying $a-b > 2$, $b-c > 2$, $a-c<p-3$.  
  \begin{enumerate}
  \item Suppose that
    $\WD(\rho)|_{I_\Qp}\cong\wt{\omega}^{a}\oplus\wt{\omega}^{b}\oplus\wt{\omega}^{c}$.
    Then
    $\rhobar|_{I_\qp}\cong\psi\oplus\psi^p\oplus\psi^{p^2},$ where
    \[\psi=\omega_3^{(a+a_0)+p(c+a_2)+p^2(b+a_1)}\]or  \[\psi=\omega_3^{(a+2-a_2)+p(b+2-a_1)+p^2(c+2-a_0)}\]with
    $(a_0,a_1,a_2) \in \{(1,1,1),(1,2,0),(2,1,0)\}$.
  \item  Suppose that
    $\WD(\rho)|_{I_\Qp}\cong\wt{\omega}_3^{a+pb+p^2c}\oplus\wt{\omega}_3^{b+pc+p^2a}\oplus\wt{\omega}_3^{c+pa+p^2b}$.
    Then
    $\rhobar|_{I_\qp}\cong\psi\oplus\psi^p\oplus\psi^{p^2}$, where
    \[\psi=\omega_3^{(a+a_0)+p(c+a_2)+p^2(b+a_1)}\]
    with $(a_0,a_1,a_2) \in \{(0,2,1),(1,1,1),(1,2,0)\}$, or
    \[\psi=\omega_3^{(a+a_0)+p(b+a_2)+p^2(c+a_1)}\] with $a_0$,
    $a_1$, $a_2\in[0,2]$ such that $a_1+a_2+a_3=3$.
  \item  Suppose that
    $\WD(\rho)|_{I_\Qp}\cong\wt{\omega}_3^{c+pb+p^2a}\oplus\wt{\omega}_3^{b+pa+p^2c}\oplus\wt{\omega}_3^{a+pc+p^2b}$.
    Then
    $\rhobar|_{I_\qp}\cong\psi\oplus\psi^p\oplus\psi^{p^2}$, where
    \[\psi=\omega_3^{(c+2-a_0)+p(a+2-a_2)+p^2(b+2-a_1)}\]
    with $(a_0,a_1,a_2) \in \{(0,2,1),(1,1,1),(1,2,0)\}$, or
    \[\psi=\omega_3^{(c+2-a_0)+p(b+2-a_2)+p^2(a+2-a_1)}\] with $a_0$,
    $a_1$, $a_2\in[0,2]$ such that $a_1+a_2+a_3=3$.
  \end{enumerate}
\end{thm}

Here we write $\omega$ for $\omega_1$, the mod $p$ cyclotomic character, and $\wt\omega$ for its Teichm\"uller lift.
Note that we consider $\wt{\omega}_3$ as $\qpb$-valued character and $\omega_3$ as $\fpb$-valued character.
This depends on a choice of embedding $K_0 \to \qpb$, but all the statements in the theorem are independent of that choice.  

\begin{proof}
  As before we let $K' = K_0$ be the unique unramified cubic extension of $\Qp$ and let $K=K_0((-p)^{1/(p^3-1)})$.
  We may assume (by a standard Baire category argument) that $\rho:G_{\Qp}\to\GL_3(\cO_E)$ for some finite
  extension $E/\Qp$. (Take $E$ large enough so as to contain the images of all embeddings $K \to \qpb$.)
  The assumptions on $a$, $b$, $c$ show that $p>7$, so we may apply Proposition \ref{prop: form of
    descent data on breuil module} with $n=d=3$ and $r=2 <
  p-1$
  to $\rho|_{G_{K_0}}$,
  we see that there is an object $\M$ of $\FBrModdd[2]$ with
  $\Tst^{*,2}(\M)\cong\rhobar$, and that $\Gal(K/K_0)$ acts on
  some basis of $\M$ via the characters $\omega^a \otimes 1$, $\omega^b \otimes 1$,
  $\omega^c \otimes 1$ in the first case, or $\omega_3^{a+pb+p^2c} \otimes 1$,
  $\omega_3^{b+pc+p^2a} \otimes 1$, $\omega_3^{c+pa+p^2b} \otimes 1$ in the second case,
  or $\omega_3^{c+pb+p^2a} \otimes 1$, $\omega_3^{b+pa+p^2c} \otimes 1$, $\omega_3^{a+pc+p^2b} \otimes 1$ in the third case.
  (Note that, for example in the first case, we at first get a basis on which $\Gal(K/K_0)$ acts via
  $1 \otimes \sigma_0 \omega^a$, $1 \otimes \sigma_0 \omega^b$, $1 \otimes \sigma_0 \omega^c$, but then
  we can just permute the induced basis on each $e_i\M$ to get the desired basis for $\M$.)

Since $\rhobar$ is assumed to be irreducible, $\rhobar|_{G_{K_0}}\cong
\chi\oplus\chi^p\oplus\chi^{p^2}$ for some character
$\chi:G_{K_0}\to\F^\times$ (e.g.\ by Lemma~2.16 of \cite{bib:ADP}). By Corollary
\ref{cor:unique-sub}, there is a rank one submodule $\cN$ of $\M$ with
$\Tst^{*,2}(\cN)\cong\chi$. Then by Lemma \ref{lem: form of rank one
  objects and their generic fibers} $\cN$ has the form  \begin{itemize}
  \item $\cN=((k\otimes_\Fp\F)[u]/u^{ep})\cdot n$,
  \item $e_i\cN_2=u^{r_i}e_i\cN$,
  \item $\varphi_2(\sum_{i=0}^{2}u^{r_i}e_i n)=\lambda n$ for some
    $\lambda\in (k\otimes_\Fp\F)^\times$, and
  \item $\hat{g}(n)=(\sum_{i=0}^{2}(\omega_3(g)^{k_i}\otimes
    1)e_i)n$ for all $g\in\Gal(K/K_0)$,
  \end{itemize}where each $k_i\in \Z$, $r_i \in [0,2e]$, and $r_i\equiv p^2k_{i+1}-k_i\pmod{e}$.
  Note that $e = p^3-1$.
  We now consider the three cases in the statement of the theorem separately.

\emph{First case:} if
    $\WD(\rho)|_{I_\Qp}\cong\wt{\omega}^{a}\oplus\wt{\omega}^{b}\oplus\wt{\omega}^{c}$,
    then $\WD(\det\rho)|_{I_{K_0}}=\wt{\omega}^{a+b+c}$, so that
    $(\WD(\det\rho)\wt{\omega}^{-(a+b+c)})|_{I_{K_0}}=1$. Thus
    $(\det\rho)|_{G_{K_0}}\wt{\omega}^{-(a+b+c)}$ is a crystalline
    character with all Hodge--Tate weights equal to $-3$ (since $\rho$
    is assumed to have Hodge--Tate weights $-2$, $-1$ and $0$), so we see
    that
    $(\det\rho)\wt{\omega}^{-(a+b+c)}|_{I_{K_0}}=\varepsilon^3$, and
    in
    particular \[\det\rhobar|_{I_{K_0}}=\chi^{1+p+p^2}|_{I_{K_0}}=\omega^{a+b+c+3}.\]

By Proposition~\ref{prop:descent-data-on-rank-one-subobject} and the
assumptions on $a$, $b$, $c$ we may assume that 
the $k_i$ are a permutation of $(1+p+p^2)a$, $(1+p+p^2)b$,
$(1+p+p^2)c$.
(Consider $\cN/u\cN$ to compare the $k_i$ with the $a_{\pi(i)}$ in the proposition.)
Write $k_i = (1+p+p^2)x_i$ with $x_i \in \{a,b,c\}$, so that $r_i\equiv
p^2k_{i+1}-k_i\equiv (1+p+p^2)(x_{i+1}-x_i)\pmod{e}$. Since
$r_i\in[0,2e]$ and by our conditions on $(a,b,c)$, we can write $r_i=(1+p+p^2)(x_{i+1}-x_i)+a_ie$ with
$a_i\in[0,2]$. By Lemma~\ref{lem: form of rank one objects and their
  generic fibers}, we see that one of $\chi|_{I_{K_0}}$, $\chi^p|_{I_{K_0}}$ and
$\chi^{p^2}|_{I_{K_0}}$ is equal
to \[\omega_3^{(x_0+a_0)+p(x_2+a_2)+p^2(x_1+a_1)}.\] Since
$\chi^{1+p+p^2}|_{I_{K_0}}=\omega^{a+b+c+3}$, we see
that \[x_0+x_1+x_2+a_0+a_1+a_2\equiv a+b+c+3\pmod{p-1},\] so
$a_0+a_1+a_2=3$. By cyclic symmetry (our choice of $\sigma_0$) we may assume that either $x_0=a$,
$x_1=b$, $x_2=c$ or $x_0=a$, $x_1=c$, $x_2=b$. Our conditions on $(a,b,c)$ and that
each $r_i\in[0,2e]$ then give the claimed result.

\emph{Second case:} if
  $\WD(\rho)|_{I_\Qp}\cong\wt{\omega}_3^{a+pb+p^2c}\oplus\wt{\omega}_3^{b+pc+p^2a}\oplus\wt{\omega}_3^{c+pa+p^2b}$,
  then arguing as in the previous case we see that
  $\chi^{1+p+p^2}|_{I_{K_0}}=\omega^{a+b+c+3}$, and by
  Proposition~\ref{prop:descent-data-on-rank-one-subobject} and the
  assumptions on $a$, $b$, $c$ we may assume that the $k_i$ are
  a permutation of $a+pb+p^2c$, $b+pc+p^2a$, $c+pa+p^2b$.
  (To check the condition on the $a_j$ in the proposition, it suffices to show that $2(p^2+p+1) < p^2x + py + z <
  (p-3)(p^2+p+1)$ whenever $x$, $y$, $z$ are in $\{\pm(a-b),\pm(b-c),\pm(c-a)\}$ and $x > 0$. The second inequality is
  obvious.  For the first, note that $x \ge 3$ and that $y$ and $z$ are bounded below by $-(p-4)$.)
  Write $k_i = x_i + p x_i' + p^2 x_i''$, where $(x_i,x_i',x_i'')$ is a cyclic permutation
  of $(a,b,c)$. So by our conditions on $(a,b,c)$ and since $r_i \in [0,2e]$, we can write
  $r_i = (x_{i+1}' + px_{i+1}'' + p^2x_{i+1})-(x_i + p x_i' + p^2 x_i'')+a_i e$ with $a_i\in[0,2]$.
  By Lemma \ref{lem: form of rank one objects and their
  generic fibers}, we see that one of $\chi|_{I_{K_0}}$, $\chi^p|_{I_{K_0}}$ and $\chi^{p^2}|_{I_{K_0}}$ is equal
  to \[\omega_3^{(x_0+a_0)+p(x_2+a_2)+p^2(x_1+a_1)}.\] Just as above we see that
  $a_0+a_1+a_2=3$. Again we may assume that either $x_0=a$,
  $x_1=b$, $x_2=c$ or $x_0=a$, $x_1=c$, $x_2=b$. The possibilities for the $a_i$ are determined by
  the signs of $x_{i+1}-x_i''$. So in the first case we find that $a_0 \ne 2$, $a_1 \ne 0$, and $a_2 \ne 2$,
  whereas in the second case we get no restriction on the $a_i$
  (because we have $r_i=a_ie$ in this case).

\emph{Third case:} this is analogous to the second case, or can be deduced from the second case by duality.
\end{proof}

\section{Abstract framework}\label{sec: abstract framework}

In order to avoid too much notational complexity, and with
an eye to future applications and possible generalisations of our work,
in this section we set up
an axiomatic framework for the general setting in which our results
will apply. In Section \ref{sec:global-applications} we will establish that our
axioms hold in the particular case of spaces of automorphic forms on
definite unitary groups.

It will be convenient to introduce our axiomatic framework
in two stages.  Since our weight cycling arguments take place
in the context of representations over $\Fbar_p$
of the group $\GL_n(\Q_p)$,
it is most natural to place our framework in the setting of such
$\Fbar_p$-representations.  This is what we do in Subsection~\ref{subsec:mod p framework}
below.  On the other hand, in practice,
our axioms
will be established by a comparison between 
a mod $p$ and a characteristic zero setting,
and it will be convenient to axiomatise this argument as well.  The relevant axiomatic
framework for this argument is the subject of
Subsection~\ref{subsec:char zero framework},
and the argument itself is the subject of Subsections~\ref{subsec:hecke-action-at-p}
and~\ref{subsec:niveau-3-implies}.
In Subsection~\ref{subsec:some-repr-gl_n} we recall some basic facts and notation
related to representations of $\GL_n$.  

\subsection{Some representations of $\GL_n$}
\label{subsec:some-repr-gl_n}

Let $F$ be a number field.
From now on until Section~\ref{subsec:weight-cycling-gl_3} we fix a place
 $w$ of $F$ lying
over $p$.
Fix a uniformiser $\varpi$ of $F_w$ and write $k_w$ for the residue field of~$F_w$.

\subsubsection{Dual Weyl modules}
\label{subsubsec:dual-weyl-modules}
\eqninc

Let $\Z^n_+$ denote the set of tuples $(\lambda_1,\dots,\lambda_n)$ of
integers with $\lambda_1\ge \lambda_2\ge\dots\ge \lambda_n$. Let $B_n \subset \GL_n$ be the Borel subgroup of upper-triangular
matrices and let $T_n \subset B_n$ be the diagonal maximal torus. We can then view any $\lambda\in\Z^n_+$
as a dominant element of $X^*(T_n) = \Hom(T_n,\mathbb G_m)$, and
let $M_\lambda$ be the algebraic $\Z$-representation of $\GL_n$ given
by
$$M_{\lambda}:=\Ind_{B_n}^{\GL_n}(w_0\lambda)_{/\Z},$$
where
$w_0$ is the longest element
of the Weyl group (see \cite{bib:Jan-reps} for more details of these
notions). Then for any commutative ring $A$ we have that $M_\lambda(A) \cong M_\lambda(\Z) \otimes_\Z A$ is a finite
free $A$-module with a natural action of $\GL_n(A)$ that is functorial in $A$.
For any $\lambda\in(\Z^n_+)^{\Hom(F_w,\Qpbar)}$, let $W_\lambda$
be the finite free $\Zpbar$-module with an action of $\GL_n(\cO_{F_w})$ given
by \[W_\lambda:=\bigotimes_{\tau : F_w \to \Qpbar}\big(M_{\lambda_\tau}(\cO_{F_w})\otimes_{\cO_{F_w},\tau}\Zpbar\big).\]
Note that $W_\lambda \otimes_{\zpb} \Qpbar$ has a natural $\GL_n(F_w)$-action, as $M_{\lambda_\tau}(\O_{F_w})
\otimes_{\O_{F_w}} F_w \cong M_{\lambda_\tau} (F_w)$. 

\subsubsection{Weights of $\GL_n(k_w)$}
\label{subsubsec:weights-gln}
\eqninc

We now explicitly describe all weights (i.e.\ irreducible $\fpb$-representations) of $\GL_n(k_w)$.

Given any $a\in\Z^n_+$
we define the algebraic $\fp$-representation $N_a$ of $\GL_n$ to be the
subrepresentation of ${M_a}_{/\fp}$ that is generated by the highest weight vector.
(It is in fact the unique irreducible subrepresentation of
${M_a}_{/\fp}$; see \S II.2 in \cite{bib:Jan-reps}.)

We say that an element $a\in(\Z^n_+)^{\Hom(k_w,\Fpbar)}$ is a \emph{restricted weight} if for each $\sigma\in\Hom(k_w,\Fpbar)$
and each $1\le i\le n-1$ we have $a_{\sigma,i}-a_{\sigma,i+1} \le p-1$.
For a restricted weight $a\in(\Z^n_+)^{\Hom(k_w,\Fpbar)}$, we define
an $\Fpbar$-representation $F_a$ of $\GL_n(k_w)$ by
\begin{equation}
F_a:=\bigotimes_{\sigma : k_w \to \Fpbar}\big(N_{a_\sigma}(k_w)\otimes_{k_w,\sigma}\Fpbar\big).\label{eq:14}
\end{equation}
It is irreducible, and every irreducible $\Fpbar$-representation of $\GL_n(k_w)$ is
of the form $F_a$ for some $a$. 
For restricted weights $a$ and $b$ we have $F_a\cong F_b$ as representations of $\GL_n(k_w)$ if and only if
for each $\sigma\in\Hom(k_w,\Fpbar)$, we
have \[a_{\sigma,i}-a_{\sigma,i+1}=b_{\sigma,i}-b_{\sigma,i+1}\] for each
$1\le i\le n-1$ and the
character $k_w^\times\to\fpb^\times$ given
by \[x\mapsto\prod_{\sigma\in\Hom(k_w,\Fpbar)}\sigma(x)^{a_{\sigma,n}-b_{\sigma,n}}\]is
trivial. 
(For all this, see for example Theorems~3.9 and~3.10 of and
the appendix to~\cite{bib:herzig-thesis}.)

If $F_w/\qp$ is unramified, the sets $\Hom(F_w,\Qpbar)$ and $\Hom(k_w,\Fpbar)$ are
in natural bijection. We can thus
define $\lambda \in (\Z^n_+)^{\Hom(F_w,\Qpbar)}$ by demanding that $\lambda_\tau = a_{\o \tau}$ for all
$\tau$ and we call it the \emph{lift} of $a$. Since $N_a$ is by definition a subrepresentation of ${M_a}_{/\fp}$,
it follows that we have a natural $\GL_n(k_w)$-linear map $N_a(k_w) \into M_a(\O_{F_w}) \otimes_{\O_{F_w}} k_w$
and hence a natural $\GL_n(k_w)$-linear map $F_a \into W_\lambda \otimes_{\zpb} \fpb$.

\ssinc
\subsubsection{Locally algebraic modules}
\label{subsubsec:locally algebraic modules}
\eqninc

If $G$ is an open subgroup of $\GL_n(F_w)$,
we say that a function $f:G \to \Qbar_p$ is \emph{locally algebraic}
if each element of $G$ has a neighbourhood $U$ such
that the restriction $f_{|U}$ of $f$ to $U$ is of the form
$f_{| U} :g \mapsto \phi(g),$ for some element $\phi$
of the affine coordinate ring of the affine algebraic
group  $(\Res_{F_w/\qp} \GL_n)_{/\qpb}$.

If $V$ is a finite free
$\Zbar_p$-module equipped with an action of an open subgroup
$G$ of $\GL_n(F_w)$, we say that the $G$-action
on $V$ is {\em locally algebraic} if each of the matrix coefficients
of $V$ is a locally algebraic function on $G$.
(Recall that the matrix coefficients of $V$ are the functions on
$G$ of the form $g \mapsto \langle g v, v^{\vee} \rangle,$
where $v \in V$ and $v^{\vee} \in V\dual$.)

As an example, note that the $\GL_n(\cO_{F_w})$-action on each of the
dual Weyl modules $W_{\lambda}$ introduced
in~\S\ref{subsubsec:dual-weyl-modules} is algebraic and thus locally
algebraic.

\subsection{Axioms for mod $p$ representations}
\label{subsec:mod p framework}

We keep our notation of the preceding subsection. In particular, recall that we assume that $w$ is a place
of the number field $F$ that lies over $p$. From now on until Section~\ref{subsec:weight-cycling-gl_3} we will assume
moreover that $F_w/\Q_p$ is unramified. (The only exception is Proposition~\ref{prop: newton = hodge implies reducible}.)

Let $\rbar:G_F\to\GL_n(\fpb)$ be a
continuous representation such that $\rbar|_{G_{F_w}}$ is irreducible. 
Let $\cP$ be a set of finite places of $F$ which
lie above split places of $F^+$ and do not lie over $p$, and let $\T^{\cP}$
denote the commutative $\zpb$-polynomial algebra generated by formal
variables $T_v^{(j)}$ where $1\le j\le n$ and $v\in \cP$. We will
write $\T$ for $\T^{\cP}$ when the choice of $\cP$ is clear.

Assume that $\rbar$ is unramified at each place $v\in \cP$. We define a
maximal ideal $\m$ of $\T$ with residue field $\Fpbar$ by demanding
that for each $v\in \cP$, the characteristic polynomial of
$\rbar^\vee(\Frob_v)$ is equal to the reduction modulo $\m$ of\begin{equation*}
X^n
+\dots+(-1)^j(\mathbf{N}v)^{j(j-1)/2}T_v^{(j)}X^{n-j}+\dots
 +(-1)^n(\mathbf{N}v)^{n(n-1)/2}T_v^{(n)}.
\end{equation*}
(By definition, $\m$ is determined by $\rbar$, but in the applications
$\cP$ will be chosen so that conversely $\rbar$ is determined by $\m$.)

We let $V$ denote a weight of the group $\GL_n(k_w)$. (Recall that when $k_w = \fp$ such
a weight is determined by an $n$-tuple of integers $a_1\ge a_2 \ge \dots \ge a_n$
satisfying $a_1-a_2$, \dots, $a_{n-1}-a_n\le p-1$ and we will write $V=F(a_1,\dots,a_n)$. 
For more details
of this, see~\S\ref{subsubsec:weights-gln}.)

Write $G=\GL_n$, and consider the Hecke algebra $\HH_G(V)$ defined in
Section \ref{subsec:hecke-operators}. In Section
\ref{subsec:hecke-operators} we defined the elements $\o{T}_{\mu,\varpi} \in \HH_G(V)$ for $\mu \in X_*(T_n)_-$,
which we can think of as a non-decreasing sequence of $n$ integers.
For $1 \le j \le n$ we let $\o T_j=\o{T}_{\mu\su,\varpi}$,
where $\mu\su=(0,\dots,0,1,\dots,1)$, with $n-j$ zeroes followed by
$j$ ones.

Now let $S$ be an $\Fpbar$-vector space with an action of $\T$, and suppose
that $S$ additionally carries a smooth action of $\GL_n(F_w)$,
and that the actions of $\T$ and $\GL_n(F_w)$ commute. If $V$ is
a weight of $\GL_n(k_w)$, we write $S(V)$ for
$(V \otimes_{\fpb} S)^{\GL_n(\cO_{F_w})}$, so that $S(V)$ carries commuting
actions of $\T$ and $\HH_G(V)$.
We additionally make the following assumptions, letting
$F_{w,n}$ denote the unramified extension of $F_w$ of degree $n$ and $k_{w,n}$ its residue field.
\begin{enumerate}[label=A\arabic*]
\item\label{axiom: Hecke is zero} \textsl{Hecke operators at $p$ act
    by zero.}
  For all weights $V$ of $\GL_n(k_w)$ the vector space $S(V)$ is finite dimensional, and 
  each of the Hecke operators $\o T_1$, \dots, $\o T_{n-1}$ acts nilpotently on the localisation $S(V)_\m$.
\item\label{axiom: crystalline lifts of specified type} \textsl{Lifts
    of specified weight and type.}
  Suppose that $F_w = \qp$.
  Suppose that $S(V)_\m\ne 0$ for some $V=F(a_1,\dots,a_n)$. Then
  \begin{itemize}
  \item $\rbar|_{G_{F_w}}$ has a crystalline lift with Hodge--Tate
    weights $-(a_1+n-1)$, \dots, $-(a_{n-1}+1)$, $-a_n$.
  \item Suppose that $F(a_1,\dots,a_n)$ is a Jordan--H\"older factor of the reduction mod $p$ of the cuspidal representation
    $R_{\TT}^\theta$ of $\GL_n(k_w)$ as in Subsection~\ref{subsec:inertial-ll}, for some primitive character
    $\theta : k_{w,n}\s \to \qpb\s$. Then $\rbar|_{G_{F_w}}$ has a potentially
    semistable lift with Hodge--Tate weights $-(n-1),\dots,-1,0$ and
    inertial type \[\bigoplus_{\sigma \in \Gal(k_{w,n}/k_w)} \sigma(\theta\circ
    \Art_{F_{w,n}}^{-1}).\]
  \item Suppose that $F(a_1,\dots,a_n)$ is a Jordan--H\"older factor of the reduction mod $p$ of the principal series
    representation \[\Ind_{B_n(k_w)}^{\GL_n(k_w)} (\chi_1 \otimes \chi_2 \otimes \dots \otimes \chi_n)\] of $\GL_n(k_w)$ as in
    Subsection~\ref{subsec:inertial-ll}, for some distinct characters $\chi_i : k_w\s \to \qpb\s$.
    Then $\rbar|_{G_{F_w}}$ has a potentially semistable lift with Hodge--Tate
    weights $-(n-1),\dots,-1,0$ and inertial type \[\bigoplus_{i=1}^n \chi_i \circ \Art_{F_w}^{-1}.\]
  \end{itemize}
\end{enumerate}
\begin{remark}
  We could have stated Axiom~\ref{axiom: crystalline lifts of specified type} in greater generality, but
  we decided to give the simplified statement in the case that $F_w=\Qp$ because this is all we need in our main theorems.
\end{remark}
\begin{defn}\label{df:modular}
  We will say that $\rbar$ is \emph{modular of weight} $V$ if
  $S(V)_\m\ne 0$. We let $W_w(\rbar)$ denote the set of weights $V$ for
  which $\rbar$ is modular of weight~$V$.
\end{defn}
In the eventual application to automorphic forms on definite unitary
groups, 
Axioms~\ref{axiom: Hecke is zero} and~\ref{axiom: crystalline lifts of
  specified type} will be consequences of local-global compatibility
at places dividing $p$ for the $p$-adic Galois representations
associated to automorphic forms.

The following lemma will be useful later.

\begin{lm}\label{lm:commuting-lemma}
  Suppose that Axiom~\emph{\ref{axiom: Hecke is zero}} holds. Then the smooth $\GL_n(F_w)$-representation $S$ is admissible, and
  whenever $V$ is a finite free $\fpb$-module with a smooth action of $\GL_n(\O_{F_w})$, we have a natural isomorphism
  \begin{equation*}
    S(V)_{\mf m} \cong (V\otimes_{\fpb} S_{\mf m})^{\GL_n(\O_{F_w})}.
  \end{equation*}
\end{lm}

\begin{proof}
  Suppose that $U$ is any compact open subgroup of $\GL_n(\O_{F_w})$. Then
  \begin{equation*}
    S^U \cong \Hom_U(1, S) \cong \Hom_{\GL_n(\O_{F_w})} (\Ind_U^{\GL_n(\O_{F_w})} 1, S),
  \end{equation*}
  so $\dim S^U$ is bounded above by the sum of the dimensions $\Hom_{\GL_n(\O_{F_w})}(V, S) \cong S(V\dual)$, where $V$ runs
  over all irreducible constituents of $\Ind_U^{\GL_n(\O_{F_w})} 1$. By Axiom~\ref{axiom: Hecke is zero} it follows that
  $\dim S^U < \infty$, so $S$ is admissible.

  Now for any $V$ as in the statement of the lemma, we have that $V \otimes_{\fpb} S \cong \ilim_{U} (V \otimes_{\fpb} S^U)$,
  where $U$ runs through compact open subgroups that are normal in $\GL_n(\O_{F_w})$. Note that the transition maps of the
  inductive limit are injective. It follows that
  \begin{align*}
    S(V)_\m &\cong \ilim_U ((V \otimes_{\fpb} S^U)^{\GL_n(\O_{F_w})})_\m,\\ 
    (V\otimes_{\fpb} S_{\mf m})^{\GL_n(\O_{F_w})} &\cong \ilim_U ((V \otimes_{\fpb} S^U)_\m)^{\GL_n(\O_{F_w})}.
  \end{align*}
  As $S$ is admissible, we have that $V \otimes_{\fpb} S^U$ is finite-dimensional for all $U$.
  For any finite free $\fpb$-module $M$, $M_\m$ is naturally a direct summand of $M$: it is the largest subspace
  on which all $T \not\in \m$ act invertibly. It follows that 
  $((V \otimes_{\fpb} S^U)^{\GL_n(\O_{F_w})})_\m \congto ((V \otimes_{\fpb} S^U)_\m)^{\GL_n(\O_{F_w})}$.
\end{proof}


\subsection{Axioms for characteristic zero representations}
\label{subsec:char zero framework}
We maintain the notation of the preceding subsection.
In this subsection we introduce an axiomatic framework that
relates certain characteristic $p$ and characteristic $0$ representations
of $\GL_n(F_w)$.

To this end, we suppose given a $\Zpbar$-module $\Slift$ and
an $\Fpbar$-vector space $S$, each equipped
with commuting actions of $\GL_n(F_w)$ and $\T$, together with
an embedding
\begin{equation}
\label{eqn:char p to char 0}
S \hookrightarrow \Slift \otimes_{\Zpbar} \Fpbar.
\end{equation}
that is equivariant for the two actions. Suppose that the $\GL_n(F_w)$-action on $S$ is smooth, so $S$ is as in the preceding section.
If $\wt V$ is a $\Zpbar$-module equipped
with a $\GL_n(\O_{F_w})$-action,
then we write $\Slift(\wt V) = (\wt V \otimes_{\zpb} \Slift)^{\GL_n(\O_{F_w})}$. Hence if $\wt V$ is an $A$-module
for a $\Zpbar$-algebra $A$ equipped with a $\GL_n(\O_{F_w})$-action, then $\Slift(\wt V)$ is an $A$-module as well.

For $\lambda\in (\Z^n_+)^{\Hom(F_w,\Qpbar)}$, we defined in~\S\ref{subsubsec:dual-weyl-modules} a dual
Weyl module $W_\lambda$. It is a finite free $\Zpbar$-module with a
(locally) algebraic $\GL_n(\O_{F_w})$-action, and $W_\lambda \otimes_\Zpbar \Qpbar$ has a natural $\GL_n(F_w)$-action.
Thus the $\Qpbar$-vector space
\begin{equation*}
  \Slift(W_\lambda \otimes_\Zpbar \Qpbar) \cong ((W_\lambda \otimes_\Zpbar \Qpbar) \otimes_{\zpb} \Slift)^{\GL_n(\O_{F_w})}
\end{equation*}
carries a natural action of the (commutative) double coset algebra \[ \Zpbar[\GL_n(\O_{F_w}) \bs \GL_n(F_w) /
\GL_n(\O_{F_w})]\] that commutes with the $\T$-action. Explicitly, suppose that $g \in \GL_n(F_w)$ and $f \in ((W_\lambda
\otimes_\Zpbar \Qpbar) \otimes_{\zpb} \Slift)^{\GL_n(\O_{F_w})}$. Write $\GL_n(\mc{O}_{F_w}) g \GL_n(\mc{O}_{F_w})
= \coprod_i k_i g \GL_n(\mc{O}_{F_w})$ as a finite disjoint union with $k_i \in \GL_n(\mc{O}_{F_w})$. Then
\begin{equation}\label{eq:12}
  \left[\GL_n(\O_{F_w}) g \GL_n(\O_{F_w}) \right](f) = \sum_i k_i g f,
\end{equation}
where $k_i g$ acts diagonally.

\begin{df}\label{df:hecke-ops-at-p-in-char-0}
  For $\mu \in X_*(T_n)_-$ denote by $T_{\mu,w}$ the following double coset operator on $\Slift(W_\lambda \otimes_\Zpbar \Qpbar)$:
  \[ T_{\mu,w} := \left[\GL_n(\O_{F_w}) \mu(\varpi) \GL_n(\O_{F_w}) \right]. \] Moreover, let $T_j = T_{\mu\su, w}$, where
  $\mu\su=(0,\dots,0,1,\dots,1)$, with $n-j$ zeroes followed by $j$ ones.
\end{df}

We introduce the following axioms for $\Slift$ related to the above set-up.

\begin{enumerate}[label=\~A\arabic*]
\item
\label{axiom: finite and ext of scalars}
\textsl{Finiteness and extension of scalars.} If $\wt V$ is a free finite-rank $\zpb$-module
equipped with a locally algebraic action of $\GL_n(\O_{F_w})$,
then $\Slift(\wt V)$ is a finite free $\Zpbar$-module and for $A = \qpb$ and $A = \fpb$
the natural map
$$ \Slift(\tV)\otimes_{\Zpbar}A \to
\Slift(\tV\otimes_{\Zpbar} A)$$
is an isomorphism.
\item
\label{axiom: crystalline lifts}
\textsl{Crystalline lifts.}
Suppose that $\lambda\in (\Z^n_+)^{\Hom(F_w,\Qpbar)}$. If $\Slift(W_\lambda \otimes_{\zpb} \Qpbar)_\m \ne 0$, 
then $\o r|_{G_{F_w}}$ has a crystalline lift
  $\rho : G_{F_w} \to \GL_n(\qpb)$ such that for all $\tau : F_w \to \Qpbar$ we have $\HT_\tau(\rho) = -(\lambda_{\tau,1} + n-1,\dots, \lambda_{\tau,n})$.

  If moreover $T_j$ has eigenvalue $t_j \in \Qpbar$ on $\Slift(W_\lambda \otimes_{\zpb} \Qpbar)_\m$
  for $1 \le j \le n$, then we can demand
  that $\varphi^{[F_w:\qp]}$ acting on the $\qpb$-vector space $\Dcris(\rho)$ has characteristic polynomial
  \[(X^n+\dots+(-1)^jq_w^{j(j-1)/2}t_j X^{n-j}+\dots+(-1)^nq_w^{n(n-1)/2}t_n)^{[F_w:\qp]},\] where $q_w = \# k_w = p^{[F_w:\qp]}$.
\item
\label{axiom: tame lifts}
\textsl{Lifts of tame type.}
Suppose that $R$ is a cuspidal representation $(-1)^{n-1} R^\theta_\TT$ of $\GL_n(k_w)$ for some primitive character
$\theta : k_{w,n}\s \to \qpb\s$, respectively a principal series representation $\Ind_{B_n(k_w)}^{\GL_n(k_w)} (\chi_1 \otimes
\cdots \otimes \chi_n)$ of $\GL_n(k_w)$ for some distinct characters $\chi_i : k_w\s \to \qpb\s$. (See Subsection~\ref{subsec:inertial-ll}
for definitions.)

If $\Slift(R)_\m \ne 0$, then $\o r|_{G_{F_w}}$ has a potentially semistable lift
  $\rho : G_{F_w} \to \GL_n(\qpb)$ such that for all $\tau : F_w \to \Qpbar$ we have $\HT_\tau(\rho) = -(n-1,\dots, 1, 0)$ and whose inertial type
  is \[\bigoplus_{\sigma \in \Gal(k_{w,n}/k_w)} \sigma(\theta\circ
    \Art_{F_{w,n}}^{-1}),\]
    respectively
    \[\bigoplus_{i=1}^n \chi_i \circ \Art_{F_w}^{-1}.\]
\end{enumerate}
Theorem \ref{thm:char-0-axioms-imply-char-p-axioms} below shows that
(in an obvious sense) Axioms~{\ref{axiom: finite and ext of
    scalars}--\ref{axiom: tame lifts}} imply Axioms~{\ref{axiom:
    Hecke is zero}--\ref{axiom: crystalline lifts of specified type}}.

\subsection{Hecke action at $p$}
\label{subsec:hecke-action-at-p}

We suppose that $S$ and $\Slift$ are as in the preceding subsection. In this subsection we will only require Axiom~\ref{axiom: finite and ext of scalars}.

Let $a\in(\Z^n_+)^{\Hom(k_w,\Fpbar)}$ be a restricted weight and let $\lambda \in (\Z^n_+)^{\Hom(F_w,\Qpbar)}$ denote its
lift as in \S\ref{subsubsec:weights-gln}.  From the natural injections $F_a \into W_\lambda \otimes_\Zpbar \Fpbar$ and $S
\into \Slift \otimes_\Zpbar \Fpbar$ we get a $\T$-equivariant embedding $S(F_a) \into \Slift(W_\lambda \otimes_\Zpbar
\Fpbar)$.  By Axiom~\ref{axiom: finite and ext of scalars} we can think of it as $\T$-equivariant embedding
\begin{equation}\label{eq:11}
  S(F_a) \into \Slift(W_\lambda) \otimes_\Zpbar \Fpbar.
\end{equation}

The goal of this subsection is to compare the action of the Hecke operators $\o T_{\mu,\varpi}$ on $S(F_a)$ with the action of suitably
normalised Hecke operators on $\Slift(W_{\lambda}) \otimes_\Zpbar \Qpbar$. Note that this last space
is isomorphic to $\Slift(W_\lambda \otimes_\Zpbar \Qpbar)$ by Axiom~\ref{axiom: finite and ext of scalars}, so that
the Hecke operators $T_{\mu,w}$ of Definition~\ref{df:hecke-ops-at-p-in-char-0} act on it.

\begin{prop}\label{prop: compatibility of Hecke actions in char 0 and
    char p}
  Suppose that $\lambda\in (\Z^n_+)^{\Hom(F_w,\Qpbar)}$ is the lift of the restricted weight $a$, and suppose that Axiom~{\em \ref{axiom: finite and ext of scalars}} holds.
  The action of the Hecke operator $(\prod_{\tau : F_w \to \Qpbar} \tau(\varpi)^{-\ang{\mu,\lambda_\tau}}) T_{\mu,w}$ on $\Slift(W_\lambda) \otimes_\Zpbar \Qpbar$
  has the following properties:
  \begin{enumerate}
  \item It stabilises the $\Zpbar$-lattice $\Slift(W_\lambda)$.
  \item The induced action on $\Slift(W_\lambda) \otimes_\Zpbar \Fpbar$ stabilises the subspace $S(F_a)$.
  \item The induced action on $S(F_a)$ coincides with the action of the operator $\o T_{\mu,\varpi}$ in $\HH(F_a)$.
  \end{enumerate}
\end{prop}

\begin{proof}
  Consider equation~\eqref{eq:12} for $g = \mu(\varpi)$.
  To establish part~(i) it suffices to show that $(\prod_{\tau : F_w
    \to \Qpbar} \tau(\varpi)^{-\ang{\mu,\lambda_\tau}}) \mu(\varpi)$ stabilises $W_\lambda \subset W_\lambda \otimes \Qpbar$.
  Recall that $W_{\lambda} \cong \bigotimes_{\tau : F_w \to \Qpbar} \big(M_{\lambda_\tau} (\O_{F_w}) \otimes_{\O_{F_w},\tau} \Zpbar\big)$.
  For $\nu_\tau \in X^*(T_n)$ the element $\mu(\varpi) \in T_n(F_w)$ acts on the $\nu_\tau$-weight space of
  $M_{\lambda_\tau} (\O_{F_w}) \otimes_{\O_{F_w},\tau} \Qpbar$ as the scalar $\tau(\varpi)^{\ang{\mu,\nu_\tau}}$.  Since
  $\nu_\tau \le \lambda_\tau$ (as $\lambda_\tau$ is the highest weight of $M_{\lambda_\tau}$)
  and $\mu \in X_*(T_n)_-$, we see that $\ang{\mu,\nu_\tau} \ge \ang{\mu,\lambda_\tau}$, so
  $\tau(\varpi)^{-\ang{\mu,\lambda_\tau}} \mu(\varpi)$ stabilises $M_{\lambda_\tau} (\O_{F_w}) \otimes_{\O_{F_w},\tau}
  \Zpbar$ and part~(i) follows.

  To establish part~(ii), since $S \into \Slift \otimes \Fpbar$ is
  $\GL_n(F_w)$-equivariant, it suffices to show that $(\prod_{\tau : F_w \to \Qpbar} \tau(\varpi)^{-\ang{\mu,\lambda_\tau}})
  \mu(\varpi)$ stabilises $F_a \subset W_\lambda \otimes \Fpbar$.  Note that the induced action of
  $\tau(\varpi)^{-\ang{\mu,\lambda_\tau}} \mu(\varpi)$ on the reduction $M_{\lambda_\tau} (\O_{F_w})
  \otimes_{\O_{F_w},\tau} \Fpbar$ is the linear projection onto the $\nu_\tau$-weight spaces with $\ang{\mu,\nu_\tau} =
  \ang{\mu,\lambda_\tau}$. Let $\sigma : k_w \to \Fpbar$ denote the embedding induced by~$\tau$. Since
  $N_{a_\sigma}$ has a weight space decomposition, as it is an algebraic subrepresentation of ${M_{\lambda_\tau}}_{/\fp}$, it
  follows that $\tau(\varpi)^{-\ang{\mu,\lambda_\tau}} \mu(\varpi)$ stabilises the subspace
  $N_{a_\sigma}(k_w) \otimes_{k_w,\sigma} \Fpbar$. This implies that part~(ii) holds.

  By comparing formulae~\eqref{eq:20} and~\eqref{eq:12} we see that part~(iii) holds provided $\o
  T_{\mu,\varpi}(\mu(\varpi)) \in \End_\Fpbar(F_a)$ is the linear projection onto the $\nu_\sigma$-weight spaces
  of $N_{a_\sigma}(k_w) \otimes_{k_w,\sigma} \Fpbar$
  such that $\ang{\mu,\nu_\sigma} = \ang{\mu,a_\sigma}$ for all $\sigma$.

  Fix an embedding $\sigma_0 : k_w \to \Fpbar$ and let $\sigma_i := \sigma_0^{p^i}$ for any integer $i$.  Let $f := [k_w :
  \fp]$. By Steinberg's tensor product theorem (Corollary II.3.17 in \cite{bib:Jan-reps}) we have $\bigotimes_{i=0}^{f-1}
  N_{a_{\sigma_i}}^{(i)} \cong N_{A}$, where $A := \sum_{i=0}^{f-1} p^i a_{\sigma_i}$ and where the superscript $(i)$
  denotes the $i$-th Frobenius twist (i.e.\ the composition with the $i$-th power of the Frobenius endomorphism of
  ${\GL_n}_{/\fp}$). Therefore
  \[ F_a \cong \bigotimes_\sigma \big(N_{a_\sigma}(k_w) \otimes_{k_w,\sigma} \Fpbar\big) \cong N_{A}(k_w)
  \otimes_{k_w,\sigma_0} \Fpbar. \] 
  Recall from \S\ref{subsec:hecke-operators} that $\o T_{\mu,\varpi}(\mu(\varpi))$ projects onto the subspace of
  $N_{A}$ that is invariant by the $k_w$-points of the unipotent radical of the parabolic of $\GL_n$ defined by 
  $-\mu \in X_*(T_n)_+$. It consists of the $\nu$-weight spaces of $N_{A}$ satisfying $\ang{\mu,\nu} = \ang{\mu,A}$ 
  by Lemma~2.3 in \cite{bib:herzig-classification}. We can
  write $\nu = \sum_{i=0}^{f-1} p^i\nu_{\sigma_i}$ where $\nu_{\sigma_i}$ is a weight of $N_{a_{\sigma_i}}$ and we see
  that $\ang{\mu,\nu} = \ang{\mu, A}$ if and only if $\ang{\mu,\nu_{\sigma_i}} = \ang{\mu,a_{\sigma_i}}$ for all $i$, 
  which is what we wanted to prove.
\end{proof}

\begin{coroll}\label{cor: compatibility of Hecke actions in char 0 and
    char p}
  The action of $(\prod_\tau \tau(\varpi)^{-\sum_{i=1}^j \lambda_{\tau,n+1-i}})T_j$ on
  $\Slift(W_\lambda)$ induces the action of $\o T_{j} \in \HH(F_{a})$ on
  $S(F_a) \subset \Slift(W_\lambda) \otimes_\Zpbar \Fpbar$. 
\end{coroll}

\subsection{A criterion for the vanishing of Hecke operators at $p$}
\label{subsec:niveau-3-implies}

In this subsection we explain how one can deduce information about
Hecke operators from information about the associated Galois
representation via Axiom~\ref{axiom: crystalline lifts}. In particular, we deduce that certain Hecke operators
in characteristic $p$ have all eigenvalues equal to $0$
in the situations that we consider in
this paper.

First we recall a variant of a lemma of Deligne and Serre (\cite{MR0379379}).

\begin{lm}\label{lem:Deligne--Serre}
Let $A$ be a commutative $\Zpbar$-algebra acting on a finite-rank free $\Zpbar$-module
$M$.  If $\n$ is a maximal ideal of $A$ such that $(M\otimes_{\Zpbar} \Fpbar)_{\n}
\neq 0$,
then there exists a homomorphism of
$\Zpbar$-algebras $\theta:A \to \Qpbar$ whose kernel is contained in $\n$,
such that the $\theta$-eigenspace of $M\otimes_{\Zpbar}\Qpbar$ is nonzero.
\end{lm}
\begin{proof}
Since $(M\otimes_{\Zpbar} \Fpbar)_{\n} \neq 0$,
we see that $M_{\n} \neq 0,$
and hence that $(M\otimes_{\Zpbar} \Qpbar)_{\n} \neq 0$.
Now $M\otimes_{\Zpbar} \Qpbar$ is a finite-dimensional $\Qpbar$-vector space,
of which the localisation $(M\otimes_{\Zpbar} \Qpbar)_{\n}$ 
is a subspace.  (More precisely, it is the maximal subspace
on which all of the elements $a \not\in \n$ act
invertibly.)  In particular, $(M\otimes_{\Zpbar} \Qpbar)_{\n}$ 
is again finite-dimensional.

If we let $A'$ denote the image of $A\otimes_{\Zpbar} \Qpbar$ in
$\End_\Qpbar\bigl((M\otimes_{\Zpbar} \Qpbar)_{\n}\bigr),$ then
$A'$ is a finite-dimensional commutative $\Qpbar$-algebra, and hence
admits a surjection $\theta':A' \to \Qpbar$. Since $A'$ acts faithfully
on $(M\otimes_{\Zpbar} \Qpbar)_{\n}$, the $\theta'$-eigenspace
of $(M\otimes_{\Zpbar} \Qpbar)_{\n}$ is then nonzero.  If we
let $\theta: A \to \Qpbar$ denote the homomorphism obtained by composing
$\theta'$ with the natural map $A \to A',$ then $\theta$ satisfies
the conditions of the lemma.
\end{proof}

In the remainder of this subsection we will follow the analysis of Section~2.7 of \cite{ger}, which
considers the ordinary situation. Let $E$ be an algebraic extension of
$\Qp$. Let $w|p$ be a place of $F$, and let
$\rho:G_{F_w}\to\GL_n(E)$ be a crystalline representation. Assume that
$E$ contains the images of all embeddings $F_w\into\overline{E}$. Let
$D=\Dcris(\rho)=(\Bcris\otimes_\Qp \rho)^{G_{F_w}}$. This is a weakly
admissible filtered $\varphi$-module; if $F_w^0$ is the maximal
absolutely unramified subfield of $F_w$, and $\phi^0$ is the
absolute Frobenius automorphism of $F_w^0$, then $D$ is a finite
free $F_w^0\otimes_\Qp E$-module of rank $n$ with
\begin{itemize}
\item a $\phi^0$-semilinear, $E$-linear automorphism $\varphi$, and
\item a separated and exhaustive decreasing filtration $(\Fil^i
  D_{F_w})_{i\in \Z}$ on the tensor product $D_{F_w}:=D\otimes_{F^0_w}F_w$ by
  $F_w\otimes_\Qp E$-submodules.
\end{itemize}
For each embedding $\sigma:F_w^0\into E$, we let
$D_\sigma:=D\otimes_{F_w^0\otimes E,\sigma\otimes 1}E$, so that
$D=\prod_{\sigma:F_w^0\into E}D_\sigma$. Similarly, for each embedding $\tau:F_w\into E$, we let
$D_{F_w,\tau}:=D_{F_w}\otimes_{F_w\otimes E,\tau\otimes 1}E$, so that
$D_{F_w}=\prod_{\tau:F_w\into E}D_{F_w,\tau}$. The
filtration on $D_{F_w}$ is induced by filtrations on each
$D_{F_w,\tau}$, so that we may write
$\Fil^i(D_{F_w})=\prod_{\tau:F_w\into E}\Fil^i(D_{F_w,\tau})$. Then
by definition $\HT_\tau(\rho)$ is the multiset of integers in which 
$i \in \Z$ occurs $\dim_E\gr^i(D_{F_w,\tau})$ times.

The map $\varphi^{[F_w^0:\Qp]}$ induces an isomorphism of $E$-vector
spaces $\varphi^{[F_w^0:\Qp]}:D_\sigma\isoto D_\sigma$ for each
$\sigma$. Assume now that $E$ is a finite extension of $\Qp$. Let $D'$
be a $\varphi$-stable free $F_w^0\otimes_\Qp E$-submodule of $D$. We define the Hodge and
Newton numbers of $D'$ in the usual way, by forgetting the $E$-structure, and
thinking of $D'$ as a finite-dimensional $F_w^0$-vector
space. Thus \[t_H(D')=\sum_{i\in\Z}(\dim_{F_w}\gr^iD'_{F_w})i,\] \[t_N(D')=\sum_{\alpha\in\Q}(\dim_{F^0_w}
D'_{\alpha})\alpha,\]where $D'_{\alpha}$ is the
slope-$\alpha$ part of $D'$. Since $D$ is weakly admissible, we
have $t_N(D')\ge t_H(D')$, with equality if and only if there is a
crystalline sub-$E$-representation $\rho'$ of $\rho$ with $D'=\Dcris(\rho')$.

In the following, we let $\val_p$ denote the valuation of $\qpb$,
normalised such that $\val_p(p) = 1$. We note that the next
Proposition does \emph{not} require the running assumption that the
extension $F_w/\Qp$ is unramified, so we phrase the statement and
proof without this assumption.

\begin{prop}\label{prop: newton = hodge implies reducible}
  Let $\lambda\in(\Z^n_+)^{\Hom(F_w,\Qpbar)}$, and assume that $\rbar|_{G_{F_w}}$
  has a crystalline lift $\rho : G_{F_w} \to \GL_n(\qpb)$
  such that for all $\tau : F_w \to \Qpbar$ we have $\HT_\tau(\rho) =
  -(\lambda_{\tau,1} + n-1,\dots, \lambda_{\tau,n})$, and that
  $\varphi^{[F_w:\qp]}$ acting on the $\qpb$-vector space $\Dcris(\rho)$
  has characteristic polynomial
  \[(X^n+\dots+(-1)^jq_w^{j(j-1)/2}t_j X^{n-j}+\dots+(-1)^nq_w^{n(n-1)/2}t_n)^{[F_w:\qp]},\] where $q_w = \# k_w$. 

If for some $1\le j\le n-1$ we
have
\begin{equation}
\val_p(t_j)=\frac{1}{[F_w:F_w^0]}\sum_{\tau:F_w\into \Qpbar}\sum_{i=1}^j\lambda_{\tau,n+1-i},\label{eq:8}
\end{equation}
then
$\rho|_{G_{F_w}}$ is reducible.
\end{prop}

Note that by the proof below (or in applications, by Corollary~\ref{cor: compatibility of Hecke actions in char 0 and
  char p}) we know that $\val_p(t_j)$ is always at least as big as the right-hand side.

\begin{proof}
By the usual Baire category argument, we may assume that
$\rho$ is valued in $\GL_n(E)$ for some finite extension
$E/\Qp$. Assume that $E$ has been chosen large enough that it contains
the images of all embeddings $F_w\into\Qpbar$ and that it contains all
the roots of $X^n+\dots+(-1)^jq_w^{j(j-1)/2}t_jX^{n-j}+\dots+(-1)^nq_w^{n(n-1)/2}t_n$. Assume that for some $1\le j\le n-1$ we
have \[\val_p(t_j)=\frac{1}{[F_w:F_w^0]}
\sum_{\tau:F_w\into \Qpbar}\sum_{i=1}^j\lambda_{\tau,n+1-i}.\] Let
$D=\Dcris(\rho)$, and let $\alpha_1,\dots,\alpha_n \in E$ denote
the roots
of \[X^n+\dots+(-1)^jq_w^{j(j-1)/2}t_jX^{n-j}+\dots+(-1)^nq_w^{n(n-1)/2}t_n\]
ordered so that $\val_p(\alpha_1)\le\dots\le\val_p(\alpha_n)$. 
For each embedding $\sigma : F_w^0 \into \qpb$ choose an $E$-subspace $D'_\sigma$ of $D_\sigma$
such that the eigenvalues of $\varphi^{[F_w^0:\Qp]}$ on it are $\alpha_1,\dots,\alpha_j$. Let
$D' = \prod_\sigma D'_\sigma$,
so that $D'$ is a free
$F_w^0\otimes_\Qp E$-submodule of $D$ of rank $j$. 
We can ensure that $D'$ is $\varphi$-stable by first picking $D'_{\sigma_0}$ for some fixed embedding $\sigma_0$
and then defining $D'_{\sigma\circ (\phi^0)^{-i}} := \varphi^i (D'_{\sigma_0})$ for all $i$.
Then we
have \begin{align*}t_N(D')&=\frac{\dim_{\qp}E}{\dim_{\qp}F_w^0}\sum_{i=1}^j\val_p(\alpha_i)\\&\le
\frac{\dim_{\qp}E}{\dim_{\qp}F_w^0}\bigl( \frac{j(j-1)}2[F^0_w:\Qp]+\val_p(t_j)\bigr). \end{align*}Since for each embedding
$\tau:F_w\into E$ the elements of $\HT_\tau(\rho)$ are
$\lambda_{\tau,n}<\lambda_{\tau,n-1}+1<\dots<\lambda_{\tau,1}+n-1$, we
also have
that \begin{align*}t_H(D')&\ge\frac{\dim_{\qp}E}{\dim_{\qp}F_w}\sum_{\tau:F_w\into \Qpbar}\sum_{i=1}^j(\lambda_{\tau,n+1-i}+i-1)\\&=
\frac{\dim_{\qp}E}{\dim_{\qp}F_w^0}[F_w:F_w^0]^{-1}\bigl(\frac{j(j-1)}2[F_w:\Qp]+[F_w:F_w^0]\val_p(t_j)\bigr)\\&\ge
  t_N(D'). \end{align*}Since $D$ is weakly admissible we have
$t_N(D')\ge t_H(D')$, so we have $t_H(D')=t_N(D')$, so that $\rho$ is reducible.
\end{proof}

\begin{corollary}\label{cor:supersingular-hecke-evals}
  Assume that $\rbar|_{G_{F_w}}$ is irreducible, and assume that
  Axioms~\emph{\ref{axiom: finite and ext of scalars}} and~\emph{\ref{axiom:
    crystalline lifts}} of Subsection {\ref{subsec:char zero
    framework}} hold. Then Axiom~\emph{\ref{axiom: Hecke is zero}} of Subsection {\ref{subsec:mod p framework}}
  holds.

\end{corollary}

\begin{proof}Let $V=F_a$ be a weight of $\GL_n(k_w)$ such that $S(V)_\m\ne
  0$, and let $\lambda\in(\Z^n_+)^{\Hom(F_w,\Qpbar)}$ be the lift of~$a$.
  As $S(F_a) \into \Slift(W_\lambda) \otimes_\Zpbar \Fpbar$, Axiom~\ref{axiom: finite and ext of scalars} implies that
  $S(F_a)$ is finite-dimensional.
  For $1 \le j \le n-1$ let us write $\wt T_j := (\prod_\tau \tau(\varpi)^{-\sum_{i=1}^j \lambda_{\tau,n+1-i}})T_j$. Suppose
  that $\alpha_j \in \fpb$ for $1 \le j \le n-1$ are eigenvalues of the $\wt T_j$ on $(\Slift(W_\lambda) \otimes_\Zpbar
  \Fpbar)_\m$. By Corollary \ref{cor: compatibility of Hecke actions in char 0 and char p} it suffices to show that $\alpha_j
  = 0$ for all $j$.  We apply Lemma~\ref{lem:Deligne--Serre} with $M = \Slift(W_\lambda)$, $A = \T[\wt T_1,\dots,\wt
  T_{n-1}]$, and the maximal ideal $\n$ that is the kernel of $\o\theta : A \to \T/\m = \fpb$ sending $\wt T_j$ to $\alpha_j$. The
  $\theta$-eigenspace provided by that lemma lies in $(\Slift(W_\lambda) \otimes_\Zpbar \Qpbar)_\m$ (as $\ker \theta \subset
  \ker \o\theta$) and $\theta(\wt T_j) \in \zpb$ lifts $\alpha_j$ for all $j$. If we had $\alpha_j \ne 0$ for some $j$, then the
  eigenvalue of $\wt T_j$ would be a unit, so the eigenvalue of $T_j$ would satisfy equation \eqref{eq:8}. Hence it
  would follow from Axioms~\ref{axiom: finite and ext of scalars}, \ref{axiom: crystalline lifts} and Proposition \ref{prop:
    newton = hodge implies reducible} that there is a lift $\rho$ of $\rbar|_{G_{F_w}}$ such that $\rho$ is reducible. Since
  $\rbar|_{G_{F_w}}$ is irreducible by assumption, this is impossible. The result follows.
\end{proof}
\begin{thm}\label{thm:char-0-axioms-imply-char-p-axioms}
Assume that $\rbar|_{G_{F_w}}$ is irreducible. If
Axioms~{\em 
\ref{axiom: finite and ext of scalars}--\ref{axiom: tame lifts}} of
Subsection {\ref{subsec:char zero framework}}
are satisfied by $\Slift$, and if $S$ is an $\Fpbar$-subspace of
$\Slift \otimes_{\Zpbar} \Fpbar$ which is stable under $\GL_n(F_w)$ and $\T$
and such that $\GL_n(F_w)$ acts smoothly on it,
then $S$ also satisfies
Axioms~{\em
\ref{axiom: Hecke is zero}--\ref{axiom: crystalline lifts of specified type}}
of Subsection {\ref{subsec:mod p framework}}.
\end{thm}
\begin{proof}
Corollary~\ref{cor:supersingular-hecke-evals}
shows that~\ref{axiom: finite and ext of scalars}
and~\ref{axiom: crystalline lifts} together imply~\ref{axiom: Hecke is
zero}. 
If $F_w = \qp$, and $a = (a_1,\dots,a_n)$ is restricted, let $\lambda \in \Z^n_+$ denote the lift of $a$. From the
embedding~\eqref{eq:11} we see that $(\Slift(W_\lambda) \otimes_\Zpbar \Fpbar)_\m \ne 0$, hence
$(\Slift(W_\lambda) \otimes_{\zpb} \Qpbar)_\m \ne 0$. 
From~\ref{axiom: finite and ext of scalars} and~\ref{axiom: crystalline
  lifts} we deduce that the first part of~\ref{axiom: crystalline lifts of specified type} holds.
For the second and third parts of~\ref{axiom: crystalline lifts of specified type} we argue similarly,
letting a $\GL_n(k_w)$-stable lattice $R_0$ in the cuspidal
(resp.\ principal series) representation $R$ such that $R_0 \otimes_\Zpbar \Fpbar$ contains
$F_a$ as a subrepresentation play the role of $W_\lambda$ and
using axioms \ref{axiom: finite and ext of scalars} and~\ref{axiom: tame lifts}.
\end{proof}

\section{Weight elimination}
\label{sec:elimination}

We preserve all the notation of the preceding section, and assume that
Axiom~\ref{axiom: crystalline lifts of
  specified type} holds.  We will however specialise to the case when
$n = 3$ and our fixed place $w|p$ of $F$ is split, that is,
$F_w = \qp$. The Serre-type conjecture in
\cite{bib:herzig-thesis} lets us describe a set $W_w^?(\rbar)$ of
weights in which $\rbar$ should be modular, and which is conjectured
to coincide with the set of all ``regular'' weights of~$\rbar$.
Our goal in this section is to prove ``weight elimination'', i.e.\ that if $\rbar$ is modular
of some sufficiently generic weight,
then this weight does in fact lie in the set $W_w^?(\rbar)$.

The construction of $W_w^?(\rbar)$ is recalled in Subsection~\ref{subsec:weights},
while elimination itself is proved in Subsection~\ref{subsec:elimination}.  

\subsection{The definition of $W_w^?(\rbar)$.}
\label{subsec:weights}
The set $W_w^?(\rbar)$ has a purely local definition depending only on
$\rbar|_{G_{F_w}}$, and we will sometimes refer to it as
$W^?_w(\rbar|_{G_{F_w}})$ in order to emphasise the local nature of
the definition. In Section~6 of~\cite{bib:herzig-thesis} we naturally
associated to a tame Galois representation $\rhobar : G_\qp \to
\GL_3(\fpb)$ a representation $V(\rhobar)$ of $\GL_3(\fp)$ over
$\qpb$. (It is a Deligne--Lusztig representation up to sign and only
depends on $\rhobar|_{I_\qp}$.) We also defined a certain operator $\cR$
on the set of irreducible $\fpb$-representations of $\GL_3(\fp)$. Let
$\JH(\o{V(\rhobar)})$ denote the set of Jordan--H\"older constituents of
the reduction mod~$p$ of a lattice inside $V(\rhobar)$. We let
\begin{equation*}
  W^?_w(\rbar|_{G_{F_w}}) = \cR(\JH(\o{V(\rbar|_{G_{F_w}})})).
\end{equation*}

\subsection{Elimination}
\label{subsec:elimination}
Before we get to our main result on weight elimination we need some
preliminary lemmas. 
We recall from~\S6.4 in~\cite{bib:herzig-thesis} the shorthand
notation for inertial Galois representations $\tau : I_\qp \to
\GL_3(\fpb)$ that can be extended to $G_\qp$, at least in the special
case we need. Suppose that $\xi \in S_3$ is a permutation of the set $\{1,2,3\}$ of order~3 and that $\mu = (\mu_1,\mu_2,\mu_3)
\in \Z^3$. Then we define $\tau(\xi,\mu) := \psi \oplus \psi^p \oplus
\psi^{p^2}$, where
\[\psi = \omega_3^{\sum_{i=0}^2 \mu_{\xi^i(1)} p^i}. \]
(Recall that the tame fundamental character $\omega_3$ was defined in Subsection~\ref{subsec:breuil-modules-two}.)

\begin{lm}\label{lm:gen-gl3-weights}
  Suppose that $0 \le x-y \le p-3$ and $0 \le y-z \le p-3$ and that $\rhobar : G_{F_w} \to \GL_3(\fpb)$ is irreducible.
  Then $F(x,y,z) \in W^?_w(\rhobar)$ if and only if
  \begin{equation*}
    \text{$\rhobar|_{I_{F_w}} \cong \tau(\xi,(x+2,y+1,z))$, for some $\xi \in S_3$ of order $3$,}
  \end{equation*}
  or $x-z > p-2$ and
  \begin{equation*}
    \text{$\rhobar|_{I_{F_w}} \cong \tau(\xi,(z+p,y+1,x-p+2))$, for some $\xi \in S_3$ of order $3$.}
  \end{equation*}
\end{lm}

\begin{proof}
  This follows from Proposition~7.4 of \cite{bib:herzig-thesis}.
\end{proof}

We need a lemma to help us distinguish niveau~$3$ inertial Galois representations of the form $\tau(\xi,(a,b,c))$ 
with $\xi \in S_3$ of order~$3$.

\begin{lm}\label{lm:niv3-exponents}
  Suppose $n \in \Z$. Then $n$ is not divisible by $p^2+p+1$ if and only if one of the following happens:
  \begin{enumerate}
  \item $n = x+py+p^2z$ with $x > y \ge z$ and $x-z \le p$.
  \item $n = p^2x+py+z$ with $x \ge y > z$ and $x-z \le p$.
  \end{enumerate}
  The two cases are disjoint. In either case, $(x,y,z)$ is uniquely determined by $n$.
\end{lm}

\begin{proof}
  Without loss of generality, we can assume that $0 \le n < p^2+p+1$. (Adding $p^2+p+1$ to $n$ corresponds to adding $1$ to
  each of $x$, $y$, $z$.)

  If $n \ne 0$, in base $p+1$ we have $n-1 = \alpha(p+1)+\beta$ with $0 \le \alpha \le p-1$, $0 \le \beta \le p$.
  If $\alpha+\beta \le p-1$, we take $(x,y,z) = (\alpha+\beta+1,\alpha,0)$. If $\alpha+\beta \ge p$,
  we take $(x,y,z) = (1,\alpha+2-p,\alpha+\beta+1-2p)$. This establishes the ``only if'' direction.

  Now note that $0 \le n < p^2+p+1$ is equivalent to $z = 0$ if $n$ is as in (i) (resp.\ to $x = 1$ if $n$ is as in (ii)).
  Each of these cases happens $\binom{p+1}2 = (p^2+p)/2$ times. This implies the remaining claims.
\end{proof}

\begin{lm}\label{lm:comparing-niv3-tau}
  Suppose $a > b > c$, $a-c \le p$, suppose that $x > y > z$, $x-z \le p$, and suppose that $a+b+c = x+y+z$.
  If
  \begin{equation*}
    \tau(\xi,(a,b,c)) \cong \tau(\xi',(x,y,z))
  \end{equation*}
  for some $\xi$, $\xi' \in S_3$ of order~$3$, then $\xi = \xi'$ and $(a,b,c) = (x,y,z)$.
\end{lm}

\begin{proof}
  Note that $\tau(\xi,(a,b,c))^\vee \cong \tau(\xi^{-1},(-c,-b,-a))$. Thus, by dualising if necessary, we may assume that $\xi =
  (1\;2\;3)$.  It follows that one of $p^2x +py+z$, $x+py+p^2z$ has to coincide with one of
  \begin{gather*}
    a+pb+p^2c, \\
    b+pc+p^2a \equiv (b+1)+pc+p^2(a-p), \\
    c+pa+p^2b \equiv (c+p)+p(a-1)+p^2b
  \end{gather*}
  modulo~$p^3-1$. Our conditions on $(a,b,c)$ and $(x,y,z)$ imply that the first of these five expressions is as in case (ii)
  and the others are as in case (i) of Lemma~\ref{lm:niv3-exponents}. As those two cases are disjoint, we see that $\xi' =
  (1\;2\;3)$. Lemma~\ref{lm:niv3-exponents} shows that $(x,y,z)$ has to agree with one of $(a,b,c)$, $(b+1,c,a-p)$,
  $(c+p,a-1,b)$ modulo $(p-1,p-1,p-1)\Z$. Since $a+b+c=x+y+z$ it follows that $(a,b,c) = (x,y,z)$.
\end{proof}

\begin{lm}\label{lm:comparing-weights}
  If $(x,y,z)$ and $(a,b,c)$ in $\Z^3_+$ are restricted weights, then
  $F(x,y,z) \cong F(a,b,c)$ if and only if $x-y = a-b$, $y-z = b-c$, and $x+y+z \equiv a+b+c \pmod{3(p-1)}$.
\end{lm}

\begin{proof}
  This follows immediately from the discussion in \S\ref{subsubsec:weights-gln}.
\end{proof}

We can now state and prove our main result on weight elimination.
(Recall from Definition~\ref{df:modular} that if $\rbar: G_F \to \GL_3(\Fbar_p)$ is continuous and irreducible,
then we let $W_w(\rbar)$ denote the
set of weights $V$ for which $\rbar$ is modular of weight $V$. This set depends on the axiomatic setup of Section~\ref{sec: abstract framework}.)

\begin{thm}\label{thm:elimination}
  Suppose that $\rbar : G_F \to \GL_3(\fpb)$, and let $w|p$ be a place of $F$ such that
  $\rbar|_{G_{F_w}}$ is irreducible and such that Axiom~\emph{\ref{axiom: crystalline lifts of specified type}}
  is satisfied.
  Continue to assume that $F_w = \qp$. Suppose that the weight
  $F(x,y,z)$ is such that either
  \begin{equation}\label{eq:C_0-elim-bound}
    x-z < p-3
  \end{equation}
  or that
  \begin{equation}\label{eq:C_1-elim-bound}
    x-y < p-5,\ y-z < p-5,\text{\ and\ } x-z > p+1.
  \end{equation}
  Then $F(x,y,z) \in W_w(\rbar)$ implies $F(x,y,z) \in W^?_w(\rbar|_{G_{F_w}})$.
\end{thm}

\begin{proof} Suppose that $F(x,y,z) \in W_w(\rbar)$.

  We first suppose that \eqref{eq:C_1-elim-bound} holds.
  Let $a = y$, $b = x-(p-1)$, $c = z$. By Lemma~\ref{lm:serre-wts-invts} and Frobenius reciprocity we see that $F(x,y,z)$ is
  a constituent of the reduction mod $p$ of the principal series $\Ind_{B_3(k_w)}^{\GL_3(k_w)}(\chi^a \otimes \chi^{b} \otimes
  \chi^c)$, where $B_3 \subset \GL_3$ is the Borel subgroup of upper-triangular matrices and $\chi : k_w\s = \fp\s \to \qpb\s$
  is the Teichm\"uller lift. (Note that the principal series is unchanged, up to semisimplification, when the characters are permuted.)
  By Axiom~\ref{axiom: crystalline lifts of specified type} there is a lift of $\rbar|_{G_{F_w}}$ to a potentially
  semistable representation $\rho:G_{F_w}\to\GL_3(\Qpbar)$ with Hodge--Tate weights $-2$, $-1$, $0$ and
  $\WD(\rho)|_{I_{F_w}}\cong \wt{\omega}^{a}\oplus\wt{\omega}^{b}\oplus\wt{\omega}^{c}$. By Theorem~\ref{thm:
    explicit list of possible characters for specific descent data}, we have
  $\rbar|_{I_{F_w}}\cong\psi\oplus\psi^p\oplus\psi^{p^2}$, where either
  \begin{equation}\label{eq:1}
    \psi=\omega_3^{(y+a_0)+p(z+a_2)+p^2(x-p+1+a_1)}
  \end{equation}
  or
  \begin{equation}\label{eq:2}
    \psi=\omega_3^{(y+2-a_2)+p(x-p+3-a_1)+p^2(z+2-a_0)}
  \end{equation}
  with $(a_0,a_1,a_2) \in \{(1,1,1),(1,2,0),(2,1,0)\}$.
  (Here we use that $x-y < p-3$, $y-z < p-3$, $x-z > p+1$.)

  Next we use that $F(x,y,z)$ is a constituent of the reduction of two cuspidal representations.
  Choose a $\qp$-linear embedding of $F_{w,3}$ into $\qpb$
  and let $\chi_3 : k_{w,3}\s \to F_{w,3}\s \to \qpb\s$ denote the Teichm\"uller lift.  By formula~(7.8) in~\cite{bib:herzig-thesis}
  we see that $F(x,y,z)$ is a constituent of the reduction mod~$p$ of the cuspidal representation
  $R_\TT^\theta$ in the notation of~\S\ref{subsec:inertial-ll}, where $\theta = \chi_3^{a+pb+p^2c}$ and $a = y+1$, $b =
  x-(p-1)$, $c = z-1$.  (Note that all highest weights in that formula are restricted, due to our bounds on $(x,y,z)$.)
  By Axiom~\ref{axiom: crystalline lifts of specified type} there is a lift of $\rbar|_{G_{F_w}}$ to a potentially
  semistable representation $\rho:G_{F_w}\to\GL_3(\Qpbar)$ with Hodge--Tate weights $-2$, $-1$, $0$ and
  $\WD(\rho)|_{I_{F_w}}\cong\wt{\omega}_3^{a+pb+p^2c}\oplus\wt{\omega}_3^{b+pc+p^2a}\oplus\wt{\omega}_3^{c+pa+p^2b}$.
  By Theorem~\ref{thm: explicit list of possible characters for specific descent data}, we have
  $\rbar|_{I_{F_w}}\cong\psi\oplus\psi^p\oplus\psi^{p^2}$, where either
  \begin{equation}\label{eq:3}
    \psi=\omega_3^{(y+1+a_0)+p(z-1+a_2)+p^2(x-p+1+a_1)}
  \end{equation}
  with $(a_0,a_1,a_2) \in \{(0,2,1),(1,1,1),(1,2,0)\}$
  or 
  \begin{equation}\label{eq:4}
    \psi=\omega_3^{(y+1+a_0)+p(x-p+1+a_2)+p^2(z-1+a_1)}
  \end{equation}
  with $a_0$, $a_1$, $a_2\in[0,2]$ with $a_1+a_2+a_3=3$.
  (Here we use that $x-y < p-2$, $y-z < p-5$, $x-z > p$.)

  Similarly we claim that $F(x,y,z)$ is a constituent of the reduction mod~$p$ of the cuspidal representation $R_\TT^{\theta}$ where
  $\theta = \chi_3^{c+pb+p^2a}$ and $a = x+1$, $b = z+p-1$, $c = y-1$. Equivalently, $F(x,y,z)\dual \cong F(-z,-y,-x)$
  is a constituent of the reduction of $(R_\TT^{\theta})\dual \cong R_\TT^{\theta^{-1}}$. (For the last isomorphism
  see p.~136 in \cite{bib:DL}.) This is again true by formula~(7.8) in~\cite{bib:herzig-thesis}. 
  By Axiom~\ref{axiom: crystalline lifts of specified type} there is a lift
  of $\rbar|_{G_{F_w}}$ to a potentially semistable representation $\rho:G_{F_w}\to\GL_3(\Qpbar)$ with Hodge--Tate weights
  $-2$, $-1$, $0$ and
  $\WD(\rho)|_{I_{F_w}}\cong\wt{\omega}_3^{c+pb+p^2a}\oplus\wt{\omega}_3^{b+pa+p^2c}\oplus\wt{\omega}_3^{a+pc+p^2b}$.
  By Theorem~\ref{thm: explicit list of possible characters for specific descent data}, we have
  $\rbar|_{I_{F_w}}\cong\psi\oplus\psi^p\oplus\psi^{p^2}$, where either
  \begin{equation}\label{eq:5}
    \psi=\omega_3^{(y+2-a_0)+p(z+1-a_2)+p^2(x-p+4-a_1)}
  \end{equation}
  with $a_0$, $a_1$, $a_2\in[0,2]$ with $a_1+a_2+a_3=3$, or
  \begin{equation}\label{eq:6}
    \psi=\omega_3^{(y+2-a_0)+p(x-p+3-a_2)+p^2(z+2-a_1)}
  \end{equation}
  with $(a_0,a_1,a_2) \in \{(0,2,1),(1,1,1),(1,2,0)\}$.
  (Here we use that $x-y < p-5$, $y-z < p-2$, $x-z > p$.)

  Lemma~\ref{lm:comparing-niv3-tau} allows us to distinguish between the various inertial Galois representations.  (Note that
  whenever Theorem~\ref{thm: explicit list of possible characters for specific descent data} applies, it shows that
  $\rbar|_{I_{F_w}} \cong \tau(\xi,(a+b_0,b+b_1,c+b_2))$ for some $\xi \in S_3$ of order~$3$ and some $b_i \in [0,2]$ with $\sum
  b_i = 3$. We have $a+b_0 > b+b_1 > c+b_2$ and $(a+b_0)-(c+b_2) \le p$ due to the assumptions placed on $(a,b,c)$.
  Moreover, whenever we applied that theorem, we had $a+b+c=x+y+z+p-1$.) In particular, we see that the possibilities for
  $\rbar|_{I_{F_w}}$ described in \eqref{eq:1}, \eqref{eq:3}, \eqref{eq:5} are distinct from the possibilities described in
  \eqref{eq:2}, \eqref{eq:4}, \eqref{eq:6}.

  The above results thus show that either
  \begin{equation}\label{eq:7}\rbar|_{I_{F_w}} \cong \tau((1\;3\;2), (y+b_0, x-p+1+b_1, z+b_2))\end{equation}
  with $(b_0,b_1,b_2) \in \{(1,2,0),(2,1,0)\}$ (use \eqref{eq:1} and \eqref{eq:3}), or
  \begin{equation}\label{eq:9}\rbar|_{I_{F_w}} \cong \tau((1\;2\;3),(y+b_0,x-p+1+b_1,z+b_2))\end{equation}
  with $(b_0,b_1,b_2) \in \{(1,1,1), (2,1,0)\}$ (use \eqref{eq:2} and \eqref{eq:6}).

  On the other hand, Lemma~\ref{lm:gen-gl3-weights} shows that $F(x,y,z) \in W^?_w(\rbar|_{G_{F_w}})$ is equivalent to
  $\rbar|_{I_{F_w}}$ being isomorphic to $\tau(\xi,(x+2,y+1,z))$ or $\tau(\xi,(z+p,y+1,x-p+2))$ for some $\xi \in S_3$ of order
  $3$, but it is not hard to see that this is equivalent to \eqref{eq:7} or \eqref{eq:9}.
  It follows that $F(x,y,z) \in W^?_w(\rbar|_{G_{F_w}})$.

  Now suppose that \eqref{eq:C_0-elim-bound} holds.
  By Axiom~\ref{axiom: crystalline lifts of specified type} there is a lift of $\rbar|_{G_{F_w}}$ to a crystalline representation
  $\rho:G_{F_w}\to\GL_3(\Qpbar)$ with Hodge--Tate weights $-(x+2)$, $-(y+1)$, $-z$. Since $\rbar|_{G_{F_w}}$ is irreducible
  and $(x+2)-z \le p-2$, 
  a standard calculation in Fontaine--Laffaille theory (cf.\ Theorem~1.2' of \cite{MR2441079}) shows that $\rbar|_{I_{F_w}}
  \cong \tau(\xi,(x+2,y+1,z))$ for some $\xi \in S_3$ of order~$3$ and thus $F(x,y,z) \in W^?_w(\rbar|_{G_{F_w}})$ by
  Lemma~\ref{lm:gen-gl3-weights}. (Alternatively we could proceed as in the previous case, using for example the two
  principal series representations whose reduction contains $F(x,y,z)$. But then we would need to assume that $x-y > 2$, $y-z
  > 2$, and $x-z < p-5$.)
\end{proof}

\section{Weight cycling for $\mathrm{GL}_3$}
\label{subsec:weight-cycling-gl_3}
In this section we further develop the general
weight cycling formalism of Subsection~\ref{subsec:weight-cycling}
in the particular case of the group $\GL_3(\mathbb Q_p)$, and then use it to prove our main theorem.
More precisely, in Subsection~\ref{subsec:local results} we prove Proposition~\ref{prop:implications},
which makes Corollary~\ref{cor:weight-cycling} completely explicit
in this context, while in Subsection~\ref{subsec:main-theorem}, we prove Theorem~\ref{thm:main},
which is the main theorem of the paper in the axiomatic setting of Section~\ref{sec: abstract framework}; as explained in the introduction, it establishes
that if $\rbar:G_F \to \GL_3(\Fbar_p)$ is modular of some strongly generic weight,
then the set $W_w^?(\rbar)$ of conjectured  weights consists precisely
of those generic  weights for which $\rbar$ is modular.
(See Definition~\ref{defn: strongly generic weight} for the definition
of generic and strongly generic weights.)

\subsection{Local results}
\label{subsec:local results}

We start with a basic but important lemma on parabolically induced representations for the group $\GL_3(\fp)$. We introduce
the following notation for standard parabolic subgroups of $\GL_3$: let $P_1 := P_{(0,0,-1)} = \Big(\begin{smallmatrix} * & *
& * \\ * & * & * \\ & & * \end{smallmatrix}\Big)$ and $P_2 := P_{(0,-1,-1)} = \Big(\begin{smallmatrix} * & * & * \\ & * & * \\
& * & * \end{smallmatrix}\Big)$.

Suppose that $V$ is an irreducible $\fpb$-representation of $\GL_3(\fp)$
We say that $V$ \emph{lies in the lower alcove} if $V \cong F(x,y,z)$ with $x-y \ge 0$, $y-z \ge 0$, and $x-z < p-2$.
We say that $V$ \emph{lies in the closure of the lower alcove} if the same conditions hold, except that $x-z \le p-2$.
We say that $V$ \emph{lies in the upper alcove} if $V \cong F(x,y,z)$ with $x-y < p-1$, $y-z < p-1$, and $x-z > p-2$.
(These notions do not depend on the choice of $(x,y,z)$ by Lemma~\ref{lm:comparing-weights}.)

\begin{lm}\label{lm:constit-parab-ind}
  Suppose that $(a,b,c) \in \Z^3$ with $a-b > 0$, $b-c > 0$, and $a-c < p-1$.
  Then the induced representation $\Ind_{P_2(\fp)}^{\GL_3(\fp)} \big(F(a) \otimes F(b,c)\big)$
  is of length three with irreducible constituents given by 
  \begin{equation*}
    F(b,c,a-p+1), F(b+p-1,a,c), F(a,b,c).
  \end{equation*}
  The induced representation $\Ind_{P_2(\fp)}^{\GL_3(\fp)} \big(F(a) \otimes F(c,b-p+1)\big)$
  is of length six with irreducible constituents given by 
  \begin{gather*}
    F(c+p-1,b,a-p+1), F(c+p-1,a,b), F(c+p-2,a,b+1),\\
    F(a-1,b,c+1), F(b-1,c,a-p+2), F(a,c,b-p+1).
  \end{gather*}
\end{lm}

These nine irreducible representations of $\GL_3(\fp)$ are non-isomorphic, by Lemma~\ref{lm:comparing-weights}.  
(It helps to note that the second, the fourth, and the ninth weight lie in the upper alcove, and the others lie
in the closure of the lower alcove.)

\begin{proof}
  Recall that we denote by $B_n \subset \GL_n$ the Borel subgroup of upper-triangular matrices and
  by $T_n \subset B_n$ the diagonal maximal torus.
  We have a short exact sequence
  \begin{equation*}
    0 \to F(c,b-p+1) \to \Ind_{B_2(\fp)}^{\GL_2(\fp)} \big(F(b)\otimes F(c)\big) \to F(b,c) \to 0.
  \end{equation*}
  (We get the two maps using Frobenius reciprocity together with Lemma~2.3
  in~\cite{bib:herzig-classification}.
  For dimension reasons
  exactness follows, since the representations at the ends are non-isomorphic.) Tensoring with $F(a)$ and parabolically
  inducing to $\GL_3(\fp)$ we see that the principal series representation 
$$\Pi := \Ind_{B_3(\fp)}^{\GL_3(\fp)} \big(F(a) \otimes F(b)\otimes F(c)\big)$$
is an extension between the two induced representations in the statement of this lemma. But
  $\Pi = \overline{R_1(a,b,c)}$ in the Grothendieck group, for example by Lemma~4.7 in~\cite{bib:herzig-thesis} (note that a
  factor $(-1)^{n-r}$ is missing on the left-hand side). Here, $R_\xi(\mu)$ for $\xi \in S_3$ and $\mu \in X^*(T_3)$ denotes
  a Deligne--Lusztig representation of $\GL_3(\fp)$ over $\qpb$. We can thus determine the irreducible constituents of $\Pi$
  from the formula in the proof of Proposition~7.4 in~\cite{bib:herzig-thesis}.  Decomposing dual Weyl modules using Proposition~3.18
  in~\cite{bib:herzig-thesis}, we obtain the nine irreducible representations listed above. Since they are distinct,
  for any irreducible constituent $V$ of $\Pi$, the principal series $\Pi$ has a unique quotient with socle~$V$.  We compute
  the socle of $\Ind_{P_2(\fp)}^{\GL_3(\fp)} \big(F(a) \otimes F(b,c)\big)$ using Lemma~2.3
  in~\cite{bib:herzig-classification}: if $V = F(x,y,z)$ is in the socle, then by Frobenius reciprocity we have a
  $\GL_1(\fp)\times \GL_2(\fp)$-linear map $F(z) \otimes F(x,y) \to F(a) \otimes F(b,c)$, so by our bounds on $(a,b,c)$ we
  see that $V \cong F(b,c,a-p+1)$.  Thus $\Ind_{P_2(\fp)}^{\GL_3(\fp)} \big(F(a) \otimes F(b,c)\big)$ is the unique quotient
  of~$\Pi$ with socle $F(b,c,a-p+1)$.

  Let $\alpha_2$ be the simple root $(0,1,-1) \in X^*(T_3)$, whose corresponding simple reflection is $(2\;3)$ under the natural
  identification of the Weyl group with $S_3$. Jantzen \cite{bib:Jan-princ-series} considers a natural intertwining
  map
  \begin{equation*}
    \o T_{\!\alpha_2} : \Ind_{B_3(\fp)}^{\GL_3(\fp)} \big(F(a) \otimes F(b)\otimes F(c)\big) \to
    \Ind_{B_3(\fp)}^{\GL_3(\fp)} \big(F(a) \otimes F(c)\otimes F(b)\big)
  \end{equation*}
  and gives in Satz~4.1 the following formula for its image, in the Grothendieck group of $\GL_3(\fp)$-representations over $\fpb$:
  \begin{equation*}
    \im(\o T_{\!\alpha_2}) = \frac 12 \sum_{0 \le l \le b-c} (\o{R_1(a,b-l,c+l)}-\o{R_{(2\;3)}(a,b-l,c+l)}).
  \end{equation*}
  Note that Jantzen computes the restriction of this representation to $\SL_3(\fp)$ in~\S5.2
  of~\cite{bib:Jan-princ-series}, but it is not completely straightforward to deduce the result for $\GL_3(\fp)$ due to the
  fact that the restriction of $\Pi$ to $\SL_3(\fp)$ may have repeated Jordan--H\"older factors. Instead, we redo his
  calculation for $\GL_3(\fp)$, using his result on the reduction mod~$p$ of Deligne--Lusztig representations. From the proof
  of Proposition~7.4 in~\cite{bib:herzig-thesis} we find that for all $(i,j,k) \in \Z^3$,
  $\o{R_1(i,j,k)}-\o{R_{(2\;3)}(i,j,k)}$ equals
  \begin{gather*}
    \big(W(j+p-1,i,k)-W(j+p-2,i,k+1)\big) \\
    + \big(W(i,j,k) - W(i,j-1,k+1)\big) \\
    + \big(W(j,k,i-p+1)    - W(j-1,k+1,i-p+1)\big) \\
    + \big(W(k+p-1,i,j)    - W(k+p-2,i,j+1)\big),
  \end{gather*}
  where $W(\mu)$ denotes the dual Weyl module of highest weight $\mu$. 
  Thus the sum over $l$ above is telescoping and we obtain eight terms initially. But two of them are dual Weyl modules
  with non-dominant highest weights, which we can make dominant using formula~(3.5) in~\cite{bib:herzig-thesis}.
  After simplifying we find that $\im(\o T_{\!\alpha_2})$ equals
  \begin{equation*}
    W(b+p-1,a,c) + W(a,b,c) + W(b,c,a-p+1) - W(c+p-2,a,b+1).
  \end{equation*}
  By Proposition~3.18 in~\cite{bib:herzig-thesis} this equals
  \begin{equation*}
    F(b+p-1,a,c) + F(a,b,c) + F(b,c,a-p+1).
  \end{equation*}

  Finally, $\im(\o T_{\!\alpha_2})$ is isomorphic to $\Ind_{P_2(\fp)}^{\GL_3(\fp)} \big(F(a) \otimes F(b,c)\big)$, since both
  are quotients of~$\Pi$ with socle $F(b,c,a-p+1)$. (Note that even the codomain of $\o T_{\!\alpha_2}$ has this socle.) This
  completes the proof of the first part of the lemma. The second part follows from our computation of the union of the
  constituents above.
\end{proof}

\begin{rk}
  With some more work it is possible to compute the complete submodule structure of $\Ind_{B_3(\fp)}^{\GL_3(\fp)} \big(F(a)
  \otimes F(b)\otimes F(c)\big)$. This uses also Jantzen's 3-step filtration with semisimple graded pieces of this representation
  \cite{bib:Jan-princ-series}.
\end{rk}

Now we go back to the setting of Section~\ref{sec:repr-theory},
in the special case of the algebraic
group $G = \GL_3$ with $F$ taken to be $\qp$.
(We could equally well assume that $F/\qp$ is totally ramified, but we will
not need it.) We have $K = \GL_3(\zp)$.
Suppose that $V$ is a weight of the group $\GL_3(\fp)$, i.e.\ an
irreducible $\fpb$-representation of $\GL_3(\fp)$. 
Recall that we defined Hecke operators $\o T_1 = \o T_{(0,0,1),\varpi}$ and $\o T_2 = \o T_{(0,1,1),\varpi}$ in $\HH_G(V)$ in
Subsection~\ref{subsec:mod p framework}.

\begin{prop}\label{prop:implications}
  Suppose that $V = F(x,y,z)$, where $(x,y,z)$ is a restricted weight and that $\pi$ is a smooth
  $\GL_3(\qp)$-representation over $\fpb$.
  \begin{enumerate}
  \item Suppose $x - y > 0$, $y - z > 0$, and $x - z < p-1$.
    If $\o T_1$ fails to be injective on $(V \otimes_{\fpb} \pi)^K$, then
    $(V' \otimes_{\fpb} \pi)^K \ne 0$ where $V'$ is one of
    \begin{equation*}
      F(z+p-1,x,y),\ F(x,z,y-p+1).
    \end{equation*}
    If $\o T_2$ fails to be injective on $(V \otimes_{\fpb} \pi)^K$, then
    $(V' \otimes_{\fpb} \pi)^K \ne 0$ where $V'$ is one of
    \begin{equation*}
      F(y,z,x-p+1),\ F(y+p-1,x,z).
    \end{equation*}
  \item Suppose $x - y < p-1$, $y-z < p-1$, and $x - z > p-1$.
    If $\o T_1$ fails to be injective on $(V \otimes_{\fpb} \pi)^K$, then
    $(V' \otimes_{\fpb} \pi)^K \ne 0$ where $V'$ is one of
    \begin{gather*}
      F(x,z+p-1,y),\ F(x-1,z+p-1,y+1),\\ F(y-1,x-p+1,z+1),\ F(z+p-2,y,x-p+2),\\ F(z+2p-2,x,y).
    \end{gather*}
    If $\o T_2$ fails to be injective on $(V \otimes_{\fpb} \pi)^K$, then
    $(V' \otimes_{\fpb} \pi)^K \ne 0$ where $V'$ is one of
    \begin{gather*}
      F(y,x-p+1,z),\ F(y-1,x-p+1,z+1),\\ F(x-1,z+p-1,y+1),\ F(z+p-2,y,x-p+2),\\ F(y,z,x-2p+2).
    \end{gather*}
  \end{enumerate}
  In each case, the last weight lies in the upper alcove and the other weights lie in the closure of the lower alcove.
\end{prop}

\begin{proof}
  We apply Corollary~\ref{cor:weight-cycling} with $\lambda = (0,0,1)$, respectively $\lambda = (0,1,1)$.
  In these two cases it is clear that $K \cap {}^{\lambda(\varpi)} K$ contains
  $\ker(K \to \GL_3(\fp))$, so by Lemma~\ref{lm:buildings-lemma} we see that $\pi|_K$ contains an irreducible constituent of
  \begin{equation*}
    \ker(\Ind_{P_{-\lambda}(\fp)}^{\GL_3(\fp)} V^{N_{-\lambda}(\fp)} \to V).
  \end{equation*}
  By Lemma~2.3 in~\cite{bib:herzig-classification} the induced representation equals
  \begin{equation*}
    I_1 := \Ind_{P_1(\fp)}^{\GL_3(\fp)} \big(F(x,y) \otimes F(z)\big), \text{\ resp.\ } 
    I_2 := \Ind_{P_2(\fp)}^{\GL_3(\fp)} \big(F(x) \otimes F(y,z)\big).
  \end{equation*}
  
  To compute the constituents of $I_2$, 
  we apply Lemma~\ref{lm:constit-parab-ind} with $a = x$, $b = y$, $c = z$ if $x-z < p-1$ and with $a = x$, $b = z+p-1$,
  $c = y$ if $x-z > p-1$. This completes the case when $\o T_2$ fails to be injective.
  
  We will compute the constituents of $I_1$ from those of $I_2$. The outer automorphism $\zeta : g \mapsto
  \Big(\begin{smallmatrix} &&1\\&1\\1 \end{smallmatrix}\Big)\cdot {}^t g^{-1}\cdot \Big(\begin{smallmatrix}
    &&1\\&1\\1 \end{smallmatrix}\Big)$ of $\GL_3(\fp)$ maps every irreducible $\fpb$-representation $F(x,y,z)$ to its dual
  $F(-z,-y,-x)$. (Consider formal characters on the level of algebraic group representations.)  Thus, when we apply $\zeta$
  to $I_2$ we obtain $\Ind_{P_1(\fp)}^{\GL_3(\fp)} \big(F(-z,-y) \otimes F(-x)\big)$ and its constituents are dual to
  the ones of $I_2$. It remains to relabel $(-z,-y,-x)$ as $(x,y,z)$ to obtain $I_1$. Concretely,
  we obtain the constituents of~$I_1$ from those of~$I_2$ by first dualising and then relabelling $(-z,-y,-x)$ as $(x,y,z)$.
  It is now easy to verify that we get the required constituents, so
  this completes the case when $\o T_1$ fails to be injective.
\end{proof}

\subsection{The main theorem}
\label{subsec:main-theorem}
We preserve all the notation of Section~\ref{sec: abstract framework}, and assume that Axioms~\ref{axiom: Hecke is zero}
and~\ref{axiom: crystalline lifts of specified type} holds.  As in Section~\ref{sec:elimination} we assume furthermore
that $n = 3$ and that $F_w = \qp$.
In this subsection we prove our main theorem
in the axiomatic setting. 
However, before doing so, we have to make precise the notion of
a weight being (strongly) generic.  To this end we introduce the following definitions.

\begin{df}
  Suppose that $\delta \in \Z_{\ge 0}$.  We say that an irreducible $\fpb$-representation of $\GL_3(\fp)$ is
  \emph{$\delta$-generic} if it isomorphic to $F(x,y,z)$ for some $(x,y,z) \in \Z^3_+$ such that
  \begin{gather*}
    -1+\delta < x-y < p-1-\delta, \\
    -1+\delta < y-z < p-1-\delta, \\
    |x-z-(p-2)| > \delta.
  \end{gather*}
\end{df}

Note that by Lemma~\ref{lm:comparing-weights} this definition does not depend on the choice of $(x,y,z)$.

\begin{df}\label{defn: strongly generic weight}
  We say that a weight $V$ is \emph{generic} (resp.\ \emph{strongly generic}) if $V$ is $4$-generic (resp.\ $6$-generic).
\end{df}

If $\rbar:G_F \to \GL_3(\Fbar_3)$ is continuous and irreducible, then we recall that
$W_w(\rbar)$ denotes the set of weights $V$ for which $\rbar$ is modular of weight $V$,
and we let $W\genw(\rbar) := \{ V \in W_w(\rbar) : \text{$V$ is generic} \}$.
With this, we are finally ready to state and prove our main theorem.

\begin{thm}\label{thm:main}
  Suppose that $\rbar: G_F \to \GL_3(\fpb)$ is a continuous
  representation, and let $w|p$ be a place of $F$ such that $\rbar|_{G_{F_w}}$ is irreducible and such that 
  Axioms~\emph{\ref{axiom: Hecke is zero}} and~\emph{\ref{axiom: crystalline lifts of specified type}} are 
  satisfied. Continue to assume that $F_w = \qp$.  If $\rbar$ is
  modular of some strongly generic weight, then $W\genw(\rbar) =
  W_w^?(\rbar)$. 
\end{thm}

\begin{proof}
  By Theorem~\ref{thm:elimination} we know that $W\genw(\rbar) \subset
  W_w^?(\rbar)$. To establish the reverse inclusion we will show that all weights in $W_w^?(\rbar)$ are generic
  and that $W^?_{w}(\rbar) \subset W_w(\rbar)$.

  \newcounter{step}\setcounter{step}{0}

  \addtocounter{step}{1}\newcounter{two-conclusions}\setcounter{two-conclusions}{\value{step}}
  \emph{Step \arabic{step}:} We show that either (a) $\rbar|_{I_{F_w}}$ or (b) $\rbar|_{I_{F_{w}}}\dual \otimes \varepsilonbar^2$ is
  isomorphic to $\tau((1\;2\;3),(a+2,b+1,c))$ for some $(a,b,c) \in \Z^3$ such that $a-b > 5$, $b-c > 4$, $a-c < p-7$. 
  By assumption we can pick a 6-generic weight $F(x,y,z) \in W_w(\rbar)$.  Thus $F(x,y,z) \in W^?_w(\rbar)$ and we get from
  Lemma~\ref{lm:gen-gl3-weights} that either $\rbar|_{I_{F_w}} \cong \tau(\xi, (x+2,y+1,z))$, or $x-z > p-2$ and
  $\rbar|_{I_{F_w}} \cong \tau(\xi, (z+p,y+1,x-p+2))$, where $\xi \in S_3$ is of order~3.  If $x-z < p-2$ and $\xi =
  (1\;2\;3)$ then conclusion (a) holds with $(a,b,c) = (x,y,z)$. If $x-z > p-2$ and $\rbar|_{I_{F_w}}$ is isomorphic to
  $\tau((1\;2\;3), (z+p,y+1,x-p+2))$ or to $\tau((1\;3\;2), (x+2,y+1,z)) \cong \tau((1\;2\;3), (z+p+1,y,x-p+2))$, then it is
  easy to check that conclusion (a) holds (for the obvious choices of $(a,b,c)$).  For the three remaining cases, note that either $\rbar|\dual_{G_{I_{w}}} \otimes
  \varepsilonbar^2 \cong \tau((1\;3)\xi(1\;3), (-z+2,-y+1,-x))$ or $x-z > p-2$ and $\rbar|\dual_{G_{I_{w}}} \otimes
  \varepsilonbar^2 \cong \tau((1\;3)\xi(1\;3), (-x+p,-y+1,-z-p+2))$. Compared to the first three cases, $(x,y,z)$
  is interchanged with $(-z,-y,-x)$ and $\xi$ with $(1\;3)\xi(1\;3)$. We thus find that conclusion (b) holds.

  From now on until the last step we suppose that conclusion (a) holds in Step~\arabic{two-conclusions}.

  \addtocounter{step}{1}\newcounter{nine-weights}\setcounter{nine-weights}{\value{step}}
  \emph{Step \arabic{step}:} We analyse the set $W^?_w(\rbar|_{G_{F_{w}}})$. It follows from Proposition 7.4 and Lemma~7.6 in \cite{bib:herzig-thesis}
  that it consists of the ``obvious weights in the lower alcove'',
  \begin{equation*}
    F(a,b,c), F(c+p-2,a,b+1), F(b,c-1,a-p+2),
  \end{equation*}
  the ``obvious weights in the upper alcove'',
  \begin{equation*}
    F(c+p-2,b+1,a-p+1), F(b+p-1,a+1,c-1), F(a,c,b-p+1),
  \end{equation*}
  and the ``shadow weights in the upper alcove'',
  \begin{equation*}
    F(c+p-2,b,a-p+2), F(b+p-1,a,c), F(a,c-1,b-p+2).
  \end{equation*}
  (In general we say that a weight $F(x,y,z)$ that lies in the upper alcove is the \emph{shadow} of the weight 
  $F(z+p-2,y,x-p+2)$ in the lower alcove. The two highest weights are related by a simple reflection in the affine
  Weyl group.)
  Note that the map $\theta : (a,b,c) \mapsto (c+p-2,a,b+1)$ preserves $\tau((1\;2\;3),(a+2,b+1,c))$ and the bounds $a-b >
  5$, $b-c > 4$, $a-c < p-7$. In particular, it induces a permutation of $W^?_w(\rbar|_{G_{F_{w}}})$. It cyclically permutes
  each of the three kinds of weights above. (Note that it preserves these weights, though not the particular choice of highest
  weights used in the list above. See Lemma~\ref{lm:comparing-weights}.)

  We now verify that all nine weights above are 4-generic.
  For $F(a,b,c)$ we have
\begin{gather*}
 a-b = (a-b-5)+5 > 3, \\ b-c = (b-c-4)+4 > 3,  \\ a-c = p-7-(p-7-a+c) < p-6,
\end{gather*}
  and this also covers its shadow $F(c+p-2,b,a-p+2)$.
  For $F(c+p-2,b+1,a-p+1)$ we have
  \begin{gather*}
    (c+p-2)-(b+1)=p-7-(b-c-4) < p-5, \\ (b+1)-(a-p+1)=p-5-(a-b-5) < p-5, \\ (c+p-2)-(a-p+1)=(p-7-a+c)+p+4 > p+2.
  \end{gather*}
  We have written the inequalities in such a way that they immediately
  carry over for the other six weights using $\theta$. (Note that $\theta$ permutes $\{a-b-5, b-c-4, p-7-a+c\}$.)

  In the next three steps we consider an arbitrary element $V \in W_w(\rbar)$ and analyse the effect of weight cycling.

  \addtocounter{step}{1}\newcounter{lower-alcove}\setcounter{lower-alcove}{\value{step}}
  \emph{Step \arabic{step}:} Suppose that $V$ is an obvious weight in the lower alcove.  Using $\theta$, we may
  assume without loss of generality that $V \cong F(a,b,c)$.
  By Lemma~\ref{lm:commuting-lemma} we have an $\HH_G(V)$-equivariant isomorphism
  \begin{equation*}
    S(V)_{\mf m} \cong (V\otimes_{\fpb} S_{\mf m})^{\GL_3(\O_{F_w})}.
  \end{equation*}
  Axiom~\ref{axiom: Hecke is zero} 
  implies that both $\o T_{1}$ and $\o T_{2}$ in
  $\HH_G(V)$ act nilpotently on the finite-dimensional vector space $S(V)_{\mf m}$.
  We can thus apply weight cycling for both Hecke operators.

  In the case of $\o T_{1}$ Proposition~\ref{prop:implications} shows that $S(V')_{\mf m} \ne 0$ where $V'$ is
  isomorphic to either $F(c+p-1,a,b)$ or $F(a,c,b-p+1)$.  Suppose the first possibility holds. Then $F(c+p-1,a,b)\in
  W_w(\rbar) \subset W^?_w(\rbar|_{G_{F_{w}}})$ by Theorem~\ref{thm:elimination}. Since $(c+p-1)-b < p-2$ it would have
  to equal one of the obvious weights in the lower alcove. But Lemma~\ref{lm:comparing-weights} implies that this is not the
  case, so $F(c+p-1,a,b) \not \in W_w(\rbar)$. Thus $F(a,c,b-p+1) \in W_w(\rbar)$.

  The case of $\o T_{2}$ is very similar. We eliminate $F(b,c,a-p+1)$ and deduce
  that $F(b+p-1,a,c) \in W_w(\rbar)$.

  \addtocounter{step}{1}\newcounter{upper-alcove-obvious}\setcounter{upper-alcove-obvious}{\value{step}}
  \emph{Step \arabic{step}:} Suppose that $V$ is an obvious weight in the upper alcove.  Using $\theta$, we may
  assume without loss of generality that $V \cong F(a,c,b-p+1)$. We apply weight cycling with $\o
  T_{2}$ and deduce from Proposition~\ref{prop:implications} that $S(V')_{\mf m} \ne 0$ where $V'$ is
  isomorphic to one of $F(c+p-1,a,b)$, $F(c+p-2,a,b+1)$, $F(a-1,b,c+1)$, $F(b-1,c,a-p+2)$, and $F(c+p-1,b,a-p+1)$.  As in the
  previous step we can eliminate all these possibilities except for the second.  Let us give some detail for the last
  weight, the only one that lies in the upper alcove. We have
  \begin{gather*}
    (c+p-1)-b = p-5-(b-c-4) < p-5, \\
    b-(a-p+1) = p-6-(a-b-5) < p-5, \\
    (c+p-1)-(a-p+1) = (p-7-a+c)+p+5 > p+1,
  \end{gather*}
  so we may indeed apply Theorem~\ref{thm:elimination} to conclude
  that $F(c+p-1,b,a-p+1) \in W_w^?(\rbar|_{G_{F_w}})$, a contradiction.
  Thus $F(c+p-2,a,b+1) \in W_w(\rbar)$.

  \addtocounter{step}{1}\newcounter{upper-alcove-shadow}\setcounter{upper-alcove-shadow}{\value{step}}
  \emph{Step \arabic{step}:} Suppose that $V$ is a shadow weight in the upper alcove.  Using $\theta$, we may assume
  without loss of generality that $V \cong F(b+p-1,a,c)$. We apply weight cycling with $\o
  T_{1}$ and deduce that $S(V')_{\mf m} \ne 0$ where $V'$ is
  isomorphic to one of $F(b,c,a-p+1)$, $F(b-1,c,a-p+2)$, $F(a-1,b,c+1)$, $F(c+p-2,a,b+1)$, and
  $F(c+p-1,b,a-p+1)$. In the previous two steps we eliminated all these possibilities except for the third.
  So $F(c+p-2,a,b+1) \in W_w(\rbar)$.
  
  \addtocounter{step}{1}
  \emph{Step \arabic{step}:} We now put the previous three steps together and take into account the symmetry $\theta$. Schematically
  the cycling process does the following.
  
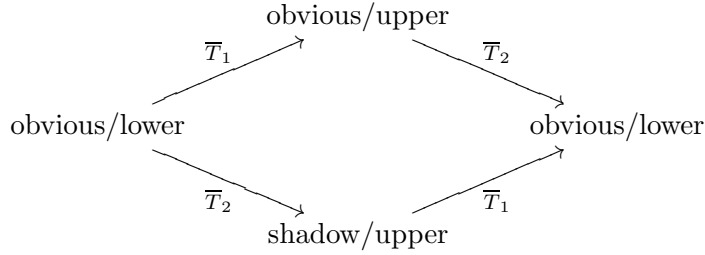
\begin{equation*}
    \xymatrix{& \text{obvious/upper}\ar[dr]^{\o T_{2}} \\ 
      \text{obvious/lower}\ar[ur]^{\o T_{1}}\ar[dr]_{\o T_{2}} && \text{obvious/lower} \\
      & \text{shadow/upper}\ar[ur]_{\o T_{1}}}
  \end{equation*}
  It is clear that it passes through all nine weights in $W^?_w(\rbar|_{G_{F_{w}}})$. (The weights on the left and on the
  right are related by $\theta$.)  This completes the proof of the theorem, in case conclusion (a) in Step~\arabic{two-conclusions} holds.

  \addtocounter{step}{1}
  \emph{Step \arabic{step}:} Finally, we deal with the case when conclusion (b) holds in Step~\arabic{two-conclusions}.
  The situation is dual to the one we considered: since
  $W^?_w(\rbar|_{G_{F_{w}}})\dual = W^?_w(\rbar|_{G_{F_{w}}}\dual \otimes \varepsilonbar^2)$ by Proposition 6.23 
  in~\cite{bib:herzig-thesis}, the weights in Step~\arabic{nine-weights} become
  dualised. A weight $V$ lies in the same alcove as its dual $V\dual$ and if one of them is generic, so is the
  other. The argument in the proof of Proposition~\ref{prop:implications} shows that for $i = 1$, $2$, the implied weights
  for $V$ if $\o T_{i}$ fails to be injective are dual to the implied weights for $V\dual$ if $\o T_{3-i}$ fails to
  be injective. (By ``implied weights'' we mean the weights $V'$ listed in Proposition~\ref{prop:implications} for a
  given weight $V$ and Hecke operator $\o T_i$.)
\end{proof}

\begin{rk}
  In Steps~\arabic{upper-alcove-obvious} and~\arabic{upper-alcove-shadow} we only applied weight cycling with one of the
  two Hecke operators. In either situation it is not hard to
  see that, if we apply weight cycling with the other Hecke operator we can only eliminate two out of the five possibilities.
\end{rk}

\section{Global applications}
\label{sec:global-applications}
In this section
we show that the axioms of Section~\ref{sec: abstract framework}
are satisfied for certain spaces of automorphic forms on definite unitary
groups, and we establish our main theorems for these groups.
Subsection~\ref{subsec:basic-facts} introduces notation and
terminology related to unitary groups and automorphic forms on them.
Subsection~\ref{subsec:galo-repr-hecke}
presents various results relating Galois
representations, automorphic forms, and Serre weights; we signal in
particular Definition~\ref{defn: modular of some Serre weight}, in
which we give precise meaning to the notion of a global residual
Galois representation being modular of a given Serre
weight. Subsection~\ref{sec:verify-char-zero} verifies that the axioms
of Section~\ref{sec: abstract framework} are satisfied, and the main
theorems are deduced in Subsection~\ref{subsection:main theorem for
  unitary groups}.
Finally, Subsection~\ref{examples:automorphic induction} uses automorphic induction of
Hecke characters to produce specific examples to which our main theorems
apply.

\subsection{Automorphic forms on unitary groups}
\label{subsec:basic-facts}
Throughout this subsection, and those that follow,
we fix an imaginary CM field $F$ with
maximal totally real subfield $F^+$,
and let $c$ denote the non-trivial element of $\Gal(F/F^+)$.
We assume that every place of
$F^+$ above $p$ splits in~$F$. 
(Note that in our main theorems, we will further assume that 
$p$ splits completely in $F$.) We also assume that $F^+\ne \Q$. 
We now proceed to introduce the various objects that
will be required for our global theorems.

\subsubsection{Unitary groups}
\eqninc
We fix a positive integer $n$ (which
will eventually be specialised to equal~$3$), 
and let $G$ be a reductive group over $F^+$ which is an outer form of
$\GL_n$, which splits over $F$, and which has the property
that for every infinite
place $v$ of~$F^+$,  $G(F^+_v)$ is compact.
(We could relax the assumption that $G$ splits over $F$ to $G$ becoming an inner
form over $F$ that splits at all places dividing $p$, but we refrain from doing so
in order to keep notation simple.)

It will be convenient to fix an integral model for $G$.  We may
find a nonzero $N\in \cO_{F^+} $ that is prime to $p$ such that $G$ admits a reductive (in particular, smooth)
model $\cG$ over $\cO_{F^+}[1/N]$, which
furthermore admits an isomorphism
$\iota:\cG_{/\cO_F[1/N]} \isoto (\GL_n)_{/\cO_F[1/N]}.$ (To see that
such a model exists, we may argue as follows. Since $G$ splits over
$F$, there is an involution ${\#}$ of the second kind on
$M:=M_n(F)$ such that for any $F^+$-algebra $A$, we have \[G(A)=\{g\in
(M\otimes_{F^+}A)^\times|gg^{\#}=1\}.\] Choose any $\O_F$-order
$M_1$ in $M$, and let $M_0:=M_1\cap M_1^\#$, so that $M_0$ is a  $\#$-stable
$\cO_F$-order in $M$. Then we have a model $G_0$ over $\O_{F^+}$ for $G$
defined by \[G_0(A)=\{g\in
(M_0\otimes_{\O_{F^+}}A)^\times|gg^{\#}=1\}\] for any $\O_{F^+}$-algebra
$A$. This model is connected reductive outside the finite set of places at which
$M_0$ is not maximal, so by localising we obtain a reductive model over some $\cO_{F^+}[1/N]$.

In order to see that we may assume that $p\nmid N$, recall that we are
free to modify an order at any finite set of places, so as places above $p$
split in $F/F^+$, we can assume that $M_0$ is maximal at each place dividing $p$. In order to
obtain the isomorphism $\iota$, it is enough to show that over some
$\cO_{F}[1/N]$ (with $p\nmid N$), $M_0$ becomes isomorphic to
$M_n(\cO_F)$. 
By weak approximation, we can choose $g \in\GL_n(F)$ so that at each place $w|p$ of $F$, we have the equality
$g_w(M_0\otimes\O_{F_w})g_w^{-1}=M_n(\O_{F_w})$. Then $g M_0 g^{-1}$ becomes
equal to $M_n(\cO_F)$ after localising at some $N$ that is prime to $p$.)


If $v$ is any finite place of $F^+$ which splits in $F$,
and if $w$ is a place of $F$ lying above $v$,
then $\iota$ gives rise to an isomorphism
$$\iota_w:G(F_v^+) \isoto \GL_n(F_w),$$
which restricts to an isomorphism
$\cG(\cO_{F_v^+}) \isoto \GL_n(\cO_{F_w})$ if $v$ does not divide $N$.
If $w^c$ is the other prime of $F$ lying over $v$,
then $c$ induces an isomorphism $c:\GL_n(F_w)\isoto \GL_n(F_{w^c}),$
and our assumption that $G$ is an outer form of $\GL_n$
ensures that $c\circ \iota_w$ is conjugate (by an element
of $\GL_n(\O_{F_{w^c}})$, if $v\nmid N$) to the inverse transpose of~$\iota_{w^c}$.

\subsubsection{Modular forms and Hecke operators away from $p$}
\eqninc
If $W$ is a $\Zpbar$-module with an action of $\cG(\cO_{F^+,p})$ (where $\cO_{F^+,p} := \prod_{v|p}\cO_{F^+_v}$), and if
$U$ is a compact open subgroup of $G(\A_{F^+}^{\infty,p})\times \cG(\cO_{F^+,p})
\subset G(\A_{F^+}^\infty)$,
then we let $S(U,W)$ denote
the space of  modular forms on $G$ of level $U$ with coefficients
$W$, i.e.\ the space of functions \[f:G(F^+)\backslash
G(\A_{F^+}^\infty)\to W\] such that $f(gu)=u_p^{-1}f(g)$ for all $u\in
U$. (Note that when $W=W_\lambda$ with $W_\lambda$ as
in~\S\ref{subsubsec:modular forms with weights},
these give examples of spaces of algebraic modular
forms in the sense of \cite{MR1729443}.)

For any $U$ as above, we may write
$G(\A_{F^+}^\infty)=\coprod_i G(F^+)t_i U$ for some finite set
$\{t_i\}$. Then there is an isomorphism \[S(U,W)\to\bigoplus_i W^{U\cap
  t_i^{-1}G(F^+)t_i}\]given by $f\mapsto (f(t_i))_i$. 
We say that
$U$ is \emph{sufficiently small} if for some finite place $v$ of $F^+$
the projection of $U$ to $G(F^+_v)$ contains no element of finite
order other than the identity. If $U$ is sufficiently
small, then for each $i$ as above we have
$U\cap t_i^{-1}G(F^+)t_i=\{1\}$,
so for any $W$ as above and any $\Zpbar$-algebra $A$,
we have a natural isomorphism
\begin{equation}
\label{eqn:general formula for extending scalars}
S(U,W)\otimes_\Zpbar A \isoto S(U,W\otimes_{\Zbar_p} A).
\end{equation}

If $U$ is any compact open subgroup as above, 
then we will say that $U$ is \emph{unramified} at
a place $v$ of $F^+$ which splits in $F$ and does not divide $N$ if $U = \cG(\cO_{F^+_v})U^v$
for a compact open subgroup $U^v$ of~$G(\A_{F^+}^{\infty,v}).$
(Recall that $\cG(\cO_{F^+_v}) \cong \GL_n(\O_{F_w})$ via
$\iota_{w}$, for either of the places $w$ of $F$ lying over $v$.)
We will say that $U$ is \emph{unramified at $p$}
if it is unramified at all places~$v|p$.
We note that $U$ {\em is} unramified at all but finitely many places of~$F^+$ that split in~$F$.

Given $U$, we let
$\cP_U$ denote the set of finite places of $F$ which lie above split places of $F^+$
at which $U$ is unramified and
which do not divide $pN$, and we let $\cP \subset \cP_U$ denote a subset with finite complement.
We then define
$\mathbb{T}^{\cP}$ as in
Subsection~\ref{subsec:mod p framework}. 
The algebra
$\mathbb{T}^{\cP}$ acts on $S(U,W)$ (for any $\cG(\O_{F^+,p})$-module~$W$) via the Hecke
operators
  \[ T_{w}^{(j)}:=  \iota_{w}^{-1} \left[ \GL_n(\mc{O}_{F_w}) \left( \begin{matrix}
      \varpi_{w}1_j & 0 \cr 0 & 1_{n-j} \end{matrix} \right)
\GL_n(\mc{O}_{F_w}) \right] 
\] for $w\in \cP$ and $\varpi_w$ a uniformiser in
$\mc{O}_{F_w}$. 
Note that $T_{w^c}\su = T_{w}\su[n-j](T_w\su[n])^{-1}$ on
$S(U,W)$.

\ssinc
\subsubsection{Modular forms of weight $\lambda$}
\label{subsubsec:modular forms with weights}
\eqninc
Let $(\Z^n_+)^{\Hom(F,\Qpbar)}_0$ be the elements $\lambda$ of
$(\Z^n_+)^{\Hom(F,\Qpbar)}$ with \[\lambda_{\sigma,i}+\lambda_{\sigma
  c,n+1-i}=0\] for all $\sigma\in\Hom(F,\Qpbar)$ and $1\le i\le
n$. For the remainder of this subsection fix $\lambda \in (\Z^n_+)^{\Hom(F,\Qpbar)}_0$.
For any place $w|p$ of $F$ let $\lambda_w$ denote the projection
of $\lambda$ to
$(\Z^n_+)^{\Hom(F_w,\qpb)}$. In~\S\ref{subsubsec:dual-weyl-modules} we
defined a finite free $\zpb$-module $W_{\lambda_w}$ with an action of
$\GL_n(\O_{F_w})$. We give this an action of $\cG(\cO_{F^+_v})$ via
$\iota_w$. This depends only on $v:=w|_{F^+}$, as $\iota_{w^c}$ is
conjugate to the inverse-transpose of $c\circ\iota_w$ (consider the
formal characters of the $M_{\lambda_\tau}$), so we will also denote
it by $W_{\lambda_v}$.
Let $W_\lambda$ be the finite free $\Zpbar$-module with an action of $\cG(\cO_{F^+,p})$ defined
by \[W_\lambda:=\bigotimes_{v|p} W_{\lambda_v}.\]

If $A$ is a $\Zpbar$-module and $U$ is a compact open subgroup
of $G(\A_{F^+}^{\infty,p})\times \cG(\cO_{F^+,p})$,
then we let \[S_\lambda(U,A):=S(U,W_\lambda\otimes_\Zpbar A).\]
From~(\ref{eqn:general formula for extending scalars}), we see that
if $U$ is sufficiently small, then there is a natural isomorphism
\begin{equation}
\label{eqn:extending scalars in algebraic setting}
S_\lambda(U,\Zpbar)\otimes_\Zpbar A \isoto S_{\lambda}(U,A).
\end{equation}

We now recall the relationship between our spaces of (algebraic) modular forms and
the space of (classical) automorphic forms on $G$. Write $S_\lambda(\Qpbar)$ for
the direct limit of the spaces $S_\lambda(U,\Qpbar)$ over compact open
subgroups $U$ of $G(\A_{F^+}^{\infty,p})\times \cG(\cO_{F^+,p})$
(with the transition maps being the inclusions
$S_\lambda(U,\Qpbar)\subset S_\lambda(V,\Qpbar)$ whenever $V\subset
U$). Concretely, $S_\lambda(\Qpbar)$ is the set of
functions \[f:G(F^+)\backslash G(\A_{F^+}^\infty)\to
W_\lambda\otimes_\Zpbar\Qpbar\] such that there is a compact open
subgroup $U$ of $G(\A_{F^+}^{\infty,p})\times \cG(\cO_{F^+,p})$
with \[f(gu)=u_p^{-1}f(g)\] for all $u\in U$, $g\in G(\A_{F^+}^{\infty})$. This
space has a natural left action of $G(\A_{F^+}^\infty)$ via \[(g\cdot
f)(h):=g_pf(hg).\]

Recall that we have fixed an isomorphism $\imath:\Qpbar\isoto\C$.
Fix a set of embeddings $\{\wt\tau : F \to \qpb\}$ such that $\Hom(F,\qpb) = \{ \wt\tau \} \sqcup \{ \wt\tau \circ c \}$.  Let
$\sigma_\lambda$ denote the representation of $G(F_\infty^+)$ given by $\bigotimes_{\wt\tau} M_{\lambda_{\wt\tau}}(\C)$,
where $G(F_\infty^+)$ acts on $M_{\lambda_{\wt\tau}}(\C)$ via the homomorphism $G(F_\infty^+) \to G(\C) \xrightarrow{\iota}
\GL_n(\C)$, where both maps are induced by $\imath \circ \wt\tau : F \to \C$. As usual, $\sigma_\lambda$ does not depend on the
choice of the $\wt\tau$ up to isomorphism.
Let $\cA$ denote the space of automorphic forms on $G(F^+)\backslash
G(\A_{F^+})$. As in the proof of Proposition 3.3.2 of \cite{cht} one
easily obtains the following lemma.

\begin{lem}
  \label{lem: relationship of algebraic automorphic forms to classical
    automorphic forms}There is an isomorphism of
  $G(\A_{F^+}^\infty)$-modules \[S_\lambda(\Qpbar)\otimes_{\Qpbar,\imath}\C\isoto\Hom_{G(F^+_\infty)}(\sigma_\lambda^\vee,\cA).\]
\end{lem}In particular, we note that $S_\lambda(\Qpbar)$ is a semisimple admissible $G(\A_{F^+}^\infty)$-module.

\subsection{Galois representations attached to modular forms}
\label{subsec:galo-repr-hecke}

In the following theorem we recall the results we will need regarding the existence
of Galois representations attached to modular forms on the unitary group~$G$. In 
the statement of the theorem,
we let $|.|^{(1-n)/2} : F\s \to \qpb\s$ denote the unique square root of $|.|^{1-n}$ 
whose composite with $\imath : \qpb \congto \C$ takes positive values. We will write $\rec_w$ for $\rec_{F_w}$.

\begin{thm}
  \label{thm: existence of Galois reps attached to algebraic modular forms}
If $\pi$ is an irreducible subrepresentation of the $G(\A_{F^+}^{\infty})$-representation
$S_\lambda(\Qpbar)$, then there is
a continuous semisimple
representation \[r_\pi:G_F\to\GL_n(\Qpbar)\] such that
\begin{enumerate}
\item $r_\pi^c\cong
  r_\pi^\vee\otimes\varepsilon^{1-n}$.
\item The representation $r_\pi$ is de Rham, and is
  crystalline if $\pi$ has level prime to $p$. If $\tau:F\into\Qpbar$
  then \[\HT_\tau(r_\pi)=\{\lambda_{\tau,1}+n-1,\dots,\lambda_{\tau,n}\}.\]
\item
  If 
  $v\nmid p$ is a place of $F^+$ which splits as $v=ww^c$ in $F$, then \[\WD(r_\pi|_{G_{F_w}})\Fss\cong
  \rec_w((\pi_v\circ\iota_w^{-1})\otimes|\cdot|^{(1-n)/2}).\]
\item If $w|p$ is a place of $F$, write $v=w|_{F^+}$. Then \[\WD(r_\pi|_{G_{F_w}})\ss\cong
 \rec_w((\pi_v\circ\iota_w^{-1})\otimes|\cdot|^{(1-n)/2})\ss.\]
\end{enumerate}

\end{thm}
\begin{proof}
  This follows from Lemma \ref{lem: relationship of algebraic
  automorphic forms to classical automorphic forms}, 
  Corollaire 5.3 of~\cite{labesse} and the main results of \cite{shin},
  \cite{chenevierharris}, \cite{ana} and \cite{blggt2}. More
  specifically, under a mild hypothesis on the weight $\lambda$, a
  representation $r_\pi$ satisfying (i), (ii) and (iii) is constructed
  in \cite{shin}, using the cohomology of unitary Shimura
  varieties. In \cite{chenevierharris} this construction is extended
  to general $\lambda$ via congruences, and properties (i) and (ii)
  are checked, and (iii) is checked up to semisimplification. That
  (iii) holds without semisimplification is the main result of
  \cite{ana}, and (iv) is the main result of \cite{blggt2}.
\end{proof}

\subsection{Serre weights}
\label{subsec:Serre weights}

We now wish to define what it means for an irreducible
representation $\rbar:G_F\to\GL_n(\Fpbar)$ to be modular of some
weight.  The basic notion required for this is that of a {\em Serre 
weight}, which we now define.

  \begin{df}\label{df:serre-weight}
    A \emph{Serre weight} is an isomorphism class of irreducible smooth $\fpb$-representations of $\cG(\cO_{F^+,p})$.
  \end{df}

\begin{remark}
\label{rem:Serre wt abuse}
We will engage in the same abuse of language with regard to Serre weights
as the one that was signalled in the case of (general)
weights in Remark~\ref{rem:abuse}. 
\end{remark}

The kernel of the natural map $\cG(\cO_{F^+,p}) \to \prod_{v | p}
 \cG(k_v)$, where $k_v$ denotes the residue field of $F^+_v$, is a normal pro-$p$ group, and so a Serre weight
may equally well be regarded as an isomorphism class of
irreducible representations of the product $\prod_{v | p} \cG(k_v)$.
If we fix a prime $\tv$ of $F$ lying over $v$ for each $v$ lying over $p$,
then the isomorphisms $\iota_{\tv}$ induce an isomorphism
\begin{equation}
\label{eqn:group isomorphism}
\prod_{v | p} \cG(k_v) \isoto \prod_{v | p} \GL_n(k_{\tv}),
\end{equation}
and hence an irreducible representation of this product
is nothing but a tensor product of irreducible representations
of the groups $\GL_n(k_{\tv})$.  This relates the
concept of Serre weight to the general notion of weight
discussed in Subsections~\ref{subsec:intro-to-weights} and
in~\S\ref{subsubsec:weights-gln}.
This allows us to
  give an ``explicit'' description of all Serre weights, as we now explain.

First note that if $w | p$, then $c$ induces an isomorphism
$k_w \isoto k_{w^c}$, and hence a bijection
$\Hom(k_{w^c},\Fbar_p) \isoto \Hom(k_w,\Fbar_p)$ via $\sigma \mapsto
\sigma c$. 
In this way we obtain an action of $c$ on 
  $(\Z^n_+)^{\coprod_{w|p}\Hom(k_w,\Fpbar)}$.
We let $(\Z^n_+)_0^{\coprod_{w|p}\Hom(k_w,\Fpbar)}$ denote
  the subset of $(\Z^n_+)^{\coprod_{w|p}\Hom(k_w,\Fpbar)}$ consisting
  of elements $a$ such that for each $w|p$, if
  $\sigma\in\Hom(k_w,\Fpbar)$ and $1\le i\le n$
  then \[a_{\sigma,i}+a_{\sigma c,n+1-i}=0.\]
We say that an element
  $a\in(\Z^n_+)_0^{\coprod_{w|p}\Hom(k_w,\Fpbar)}$ is a \emph{restricted
    weight} if, for each $w|p$,
the projection $a_w \in (\Z^n_+)^{\Hom(k_w,\Fpbar)}$ is 
  restricted in the sense of~\S\ref{subsubsec:weights-gln};
that is, if, for each $w|p$ and each $\sigma\in\Hom(k_w,\Fpbar)$, we
have \[a_{\sigma,i}-a_{\sigma,i+1}\le p-1\] for all $1\le i\le n-1$.

Let $a\in(\Z^n_+)_0^{\coprod_{w|p}\Hom(k_w,\Fpbar)}$ be a restricted weight.
For any place $w|p$ of $F$ we defined in~\S\ref{subsubsec:weights-gln} a finite free $\Fpbar$-module
$F_{a_w}$ with an action of $\GL_n(k_w)$.   We give this an action of $\cG(k_v)$ via
$\iota_w$. Just as in the case of $W_\lambda$ in \S\ref{subsubsec:modular forms with weights}, this
depends only on $v:=w|_{F^+}$, so we will also denote
it by $F_{a_v}$.

We then define the Serre
weight $F_a$ to be the tensor product
\[F_a :=\bigotimes_{v|p} F_{a_v},\]
thought of as an irreducible representation of $\prod_{v|p} \cG(k_v)$,
or equivalently, of $\cG(\cO_{F^+,p})$.

\begin{lem}
\label{lem:Serre weight classficiaton}
Every Serre weight admits a representative of the form $F_a$,
for some restricted weight
$a \in (\Z^n_+)_0^{\coprod_{w|p} \Hom(k_w,\Fbar_p)}$.
\end{lem}
\begin{proof}
This follows from the isomorphism~(\ref{eqn:group isomorphism}),
together with the fact,
noted in \S\ref{subsubsec:weights-gln},
that every irreducible representation of $\GL_n(k_w)$ is
of the form $F_{a_w}$ for some restricted weight $a_w \in (\Z^n_+)^{\Hom(k_w,\Fbar_p)}.$
\end{proof}

By virtue of this lemma (and taking into account Remark~\ref{rem:Serre
wt abuse}),
from now on we will typically denote a Serre weight
by $F_a$, for some restricted weight
$a \in (\Z^n_+)_0^{\coprod_{w|p} \Hom(k_w,\Fbar_p)}$.

\medskip

If $p$ is unramified in $F$, then we may introduce the notion
of the {\em lift} of a restricted weight.
Namely, if $p$ is unramified in $F$,
so that $\Hom(F,\Qpbar) = \coprod_{w|p}\Hom(F_w,\Qpbar)$ and $\Hom(\cO_F,\Fpbar) = \coprod_{w|p}\Hom(k_w,\Fpbar)$ are in natural bijection,
then we say that a weight $\lambda\in (\Z^n_+)_0^{\Hom(F,\Qpbar)}$ is the
\emph{lift}
of a restricted weight $a$ if for each $w|p$,
we have that
$\lambda_w$ is the lift of $a_w$ in the sense of \S\ref{subsubsec:weights-gln}; i.e.\ if for each $w|p$
and each $\tau\in\Hom(F_w,\Qpbar)$, we have
$\lambda_{\tau}=a_{\o\tau}$.

\medskip

Suppose now that $\rbar:G_F\to\GL_n(\Fpbar)$ is a continuous
irreducible representation.  If $\cP$ is a cofinite subset of the
finite places of $F$ which lie above split places of $F^+$ and do not lie over $p$ with the property that
$\rbar$ is unramified at all the places lying in $\cP$, then we may
associate a maximal ideal $\mathfrak m$ of $\mathbb{T}^{\cP}$ to
$\rbar$, as in Subsection~\ref{subsec:mod p framework}.  (We remark
that, conversely, $\m$ determines $\rbar$ by the Chebotarev density
theorem.)

\begin{defn}\label{defn: modular of some Serre weight}
If $\rbar: G_F \to \GL_n(\Fbar_p)$ is continuous and irreducible,
then we say that $\rbar$ 
{\em modular of Serre weight~$F_a$}, if
  there is a compact open subset $U$ of $G(\A_{F^+}^{\infty,p})\times \cG(\cO_{F^+,p})$
which is unramified at $p$, a subset $\cP \subset \cP_U$ with finite complement such that $\rbar$ is unramified at all
the places lying in $\cP$, and such that,
if $\mathfrak m$ is the maximal ideal in $\mathbb{T}^{\cP}$ associated to $\rbar$,
then $S(U,F_a)_{\mf{m}}\neq~0$. 

  We say that $\rbar$ is \emph{modular} if it is modular of
  some Serre weight.  
\end{defn}

  We write $W(\rbar) := \{ F_a \, | \, \rbar \text{ is modular of Serre weight } F_a\},$
and refer to this as the set of Serre weights of $\rbar$.
 (Thus $\rbar$ is modular if and only if $W(\rbar) \neq \varnothing$.)

\begin{rem}
\label{rem:shrinking the level}
If, in the context of the preceding definition, we have that $S(U,F_a)_{\mf{m}}\neq 0$,
then certainly $S(U',F_a)_{\mf{m}} \neq 0$ for any compact open subgroup $U'$ 
of $U$.  (Here we have assumed, without loss of generality, that $\cP \subset \cP_{U'} (\subset \cP_U)$,
so that $\T^\cP$ acts on both spaces of modular forms.)
Given $U$ we can choose $U' \subset U$ that is sufficiently small and unramified at $p$.
Thus, when working in the context of this definition,
we may assume that the level
  $U$ that appears is sufficiently small, and we will
  assume this from now on.
\end{rem}

\begin{remark}
  Note that if $\rbar:G_F\to\GL_n(\Fpbar)$ is modular then
  $\rbar^c\cong\rbar^\vee\varepsilonbar^{n-1}$. (Note that in
  Subsection \ref{subsec:mod p framework} we defined the maximal ideal
  $\m$ via the characteristic polynomials of $\rbar^\vee(\Frob_v)$,
  which is why we have $\varepsilonbar^{n-1}$ rather than
  $\varepsilonbar^{1-n}$ here.)
\end{remark}

\subsection{Verifying the characteristic zero axioms}
\label{sec:verify-char-zero}

We now
apply the theory developed in the earlier sections of this paper to
the specific situation of the unitary groups considered earlier in
this section. 
To this end, we will define appropriate spaces $S$ and $\Slift$ which satisfy
the axioms of Section~\ref{sec: abstract framework}.
We maintain the notation established earlier in this section,
and furthermore assume that $p$ is unramified in $F$.

Assume that $\rbar$ is modular of Serre weight $F_a$, in the sense of Definition~\ref{defn: modular of some Serre weight}, 
for some restricted weight $a\in(\Z^n_+)_0^{\coprod_{w|p}\Hom(k_{w},\Fpbar)}$.  Then, by definition, there is a
compact open subgroup $U$ of $G(\A_{F^+}^{\infty,p})\times \cG(\cO_{F^+,p})$
that is unramified at $p$ and a subset $\cP \subset \cP_U$ with finite complement
such that $S(U,F_a)_\m\ne 0$, where
$\m$ is the maximal ideal of $\T^{\cP}$ corresponding to $\rbar$. 
Shrinking $\cP$ if necessary, by Remark~\ref{rem:shrinking the level} we may (and do) assume that $U$ is sufficiently small.

We fix a place $w|p$ of $F$, and assume that $\rbar|_{G_{F_w}}$ is irreducible. We write $v=w|_{F^+}$,
and write $U = U_v \times U^v$, where $U^v \subset G(\A_{F^+}^{\infty,v})$ is compact open.

We now set up
some further notation in preparation for defining the spaces $S$ and $\Slift$.

If $W$ is any $\zpb$-module with an action of $\cG(\O_{F^+,p})$ such that $\cG(\O_{F^+_v})$ acts trivially, we define $S(U^v, W)$
to be the space of functions \[f : G(F^+)\backslash G(\A_{F^+}^\infty) \to W\] such that
$f(gu) = u_p^{-1}f(g)$ for all $g \in G(\A_{F^+}^\infty)$ and all $u \in U^v$.
It is a $\zpb$-module with an action of $G(F_v^+)$, given by $(\gamma \cdot f)(g) = f(g\gamma)$ for $\gamma \in G(F_v^+)$ and
$g \in G(\A_{F^+}^\infty)$, and a commuting action of $\mathbb{T}^{\cP}$. Let $S\sm(U^v, W)$ be the largest smooth
$G(F_v^+)$-subrepresentation of $S(U^v, W)$, i.e.\ \[ S\sm(U^v, W) = \ilim_\text{$U_v' \le U_v$ open} S(U^v, W)^{U_v'}.\]

If $W$ is a finite free $\zpb$-module with an action of $\cG(\O_{F^+,p})$ such that $\cG(\O_{F^+_v})$ acts
trivially, 
we define $S\lalg(U^v, W)$ to be the space of functions \[f : G(F^+)\backslash G(\A_{F^+}^\infty) \to W\] such that $f(gu) =
u_p^{-1}f(g)$ for all $g \in G(\A_{F^+}^\infty)$ and all $u \in U^v$ and such that for all $g \in G(\A_{F^+}^\infty)$ 
and all $\eta \in W\dual$ the map $G(F_v^+) \to \zpb$, $\gamma \mapsto \langle f(g\gamma), \eta\rangle$ is locally algebraic.
Note that $S\lalg(U^v, W)$ has natural, commuting
actions of $G(F_v^+)$ and $\mathbb{T}^{\cP}$.

We let $S=S\sm(U^v,\otimes_{v' \ne v} F_{a_{v'}})$ and $\Slift = S\lalg(U^v, \otimes_{v' \ne v} W_{\lambda_{v'}})$,
where $\lambda \in (\Z^n_+)_0^{\Hom(F,\Qpbar)}$ denotes the lift of $a$. Note that these spaces do not depend on $a_w$ and
$\lambda_w$. We identify $G(F^+_v)$ with $\GL_n(F_w)$ via
$\iota_w$. 
Then $S$ and $\Slift$ have commuting actions of $\T$ and $\GL_n(F_w)$, and the action of $\GL_n(F_w)$ on $S$ is smooth.

\begin{lm}\label{lem:reducing-S-tilde-mod-p}
  Recall that we are assuming that $U$ is sufficiently small.
  We have a natural isomorphism
  \[ \Slift \otimes_{\zpb} \fpb \congto S\sm(U^v, \otimes_{v' \ne v} (W_{\lambda_{v'}} \otimes_{\zpb} \fpb)) \]
  that is compatible with the $\T$- and $\GL_n(F_w)$-actions.
\end{lm}

\begin{proof}
  Note first that $G(F^+)\bs G(\A^\infty_{F^+})$ is compact, as $G$ is anisotropic at infinity. This implies in particular
  that the above map is well defined.
  The above map is surjective, since any element of the right-hand side takes only finitely many values 
  so that we can lift it to a smooth element of $\Slift$. (Here we use that $U$ is sufficiently small.)

  Let $f \in \Slift$ be such that $f \otimes 1$ is in the kernel of the map. Let $e_1\dual, \dots, e_d\dual$ denote a basis
  of $(\otimes_{v' \ne v} W_{\lambda_{v'}})\dual$.  By compactness of $G(F^+)\bs G(\A^\infty_{F^+})$ and by definition of
  $\Slift$ we can find a compact open subgroup $U_v' \subset \GL_n(\O_{F_w})$ with the following property: writing
  $G(\A^\infty_{F^+}) = \coprod_{i = 1}^r G(F^+) g_i U_v' U^v$ there are algebraic functions $f_{ij} : U_v' \to \zpb$ such
  that $\langle f(g_i u_v'), e_j\dual\rangle = f_{ij}(u_v')$ for all $i$ and all $u_v' \in U_v'$. Since $\Res_{F_w/\qp}
  \GL_n$ is defined over $\qp$, there is a finite extension $K$ of $\qp$ with uniformiser $\varpi_K$ such that all $f_{ij}$
  take values in $\O_K$. Since $f \otimes 1$ is in the kernel, the $f_{ij}$ take values in $\varpi_K\O_K$. It follows that
  $\varpi_K^{-1} f \in \Slift$, so $f \otimes 1 = 0$.
\end{proof}

By the lemma we get an embedding $S \into \Slift \otimes_{\zpb} \fpb$ that is compatible with the $\T$- and
$\GL_n(F_w)$-actions.

\begin{lm}\label{lem:mapping-in-locally-algebraic-V}
  If $\wt V$ is a free finite-rank $\zpb$-module equipped with a locally algebraic action of $U_v = \GL_n(\O_{F_w})$
  and $M$ is a finite free $\zpb$-module with an action of $\prod_{v' \ne v} \cG(\O_{F_v^+})$,
  then there is a natural isomorphism of $\T$-modules
  \begin{equation*}
    (\wt V \otimes_{\zpb} S\lalg(U^v, M))^{\GL_n(\O_{F_w})} \congto S(U, \wt V \otimes_{\zpb} M).
  \end{equation*}
\end{lm}

\begin{proof}
  Note that we have a $\T$- and $\GL_n(F_w)$-equivariant map $\eta : \wt V \otimes S(U^v, M) \to S(U^v, \wt V \otimes M)$
  given by $x \otimes f_0 \mapsto f(g) := x \otimes f_0(g)$. It is an isomorphism since $\wt V$ is finite free as
  $\zpb$-module.  By restricting to locally algebraic functions on the left-hand side and by passing to
  $\GL_n(\O_{F_w})$-invariants, we get a map $\eta'$ as in the statement of the lemma that is moreover injective. Pick a
  basis $e_1, \dots, e_d$ of $\wt V$.  Given any element $f \in S(U, \wt V \otimes_{\zpb} M)$, write $\eta^{-1}(f) = \sum e_i
  \otimes f_i$ for unique $f_i \in S(U^v, M)$.  To show that $\eta'$ is surjective, we need to prove that $f_i \in
  S\lalg(U^v, M)$ for all $i$.  But a simple computation, using the fact that $f$ is $\GL_n(\O_{F_w})$-invariant, shows that
  for fixed $g \in G(\A^\infty_{F^+})$ and $m\dual \in M\dual$, the function $\gamma \mapsto \langle f_i(g\gamma), m\dual \rangle$ on
  $\GL_n(\O_{F_w})$ is a linear combination of matrix coefficients of $\wt V$, hence locally algebraic.
\end{proof}

A similar but easier argument establishes the following lemma.

\begin{lm}\label{lem:mapping-in-smooth-V}
  If $V$ is a free finite-rank $\fpb$-module equipped with a smooth action of $U_v = \GL_n(\O_{F_w})$
  and $M$ is an $\fpb$-module with an action of $\prod_{v' \ne v} \cG(\O_{F_v^+})$,
  then there is a natural isomorphism of $\T$-modules
  \begin{equation*}
    (V \otimes_{\fpb} S\sm(U^v, M))^{\GL_n(\O_{F_w})} \congto S(U, V \otimes_{\fpb} M).
  \end{equation*}
\end{lm}

\begin{prop}
  \label{prop:the axioms hold for unitary groups}In the present
  setting, Axioms~\emph{\ref{axiom: finite and ext of scalars}}--\emph{\ref{axiom: tame lifts}} hold.
\end{prop}
\begin{proof}
To verify Axiom~\ref{axiom: finite and ext of scalars}, first note that $\Slift(\wt V)$ is finite free by 
Lemma~\ref{lem:mapping-in-locally-algebraic-V} (as $U$ is sufficiently small). It remains to show that
  \begin{equation*}
    (\wt V \otimes_{\zpb} \Slift)^{\GL_n(\O_{F_w})} \otimes_{\zpb} A \congto
    (\wt V \otimes_{\zpb} \Slift \otimes_{\zpb} A)^{\GL_n(\O_{F_w})}
  \end{equation*}
  when $A = \fpb$ or $A = \qpb$. In the first case we can rewrite the left-hand side using
  Lemma~\ref{lem:mapping-in-locally-algebraic-V}, (\ref{eqn:general formula for extending scalars}),
Lemma~\ref{lem:mapping-in-smooth-V}, and Lemma~\ref{lem:reducing-S-tilde-mod-p} (in this order):
  \begin{multline*}
     (\wt V \otimes_{\zpb} \Slift)^{\GL_n(\O_{F_w})} \otimes_{\zpb}
     \Fpbar \\
\cong S(U, \wt V \otimes_{\zpb} \bigotimes_{v' \ne v} W_{\lambda_{v'}})
    \otimes_{\zpb} \fpb
\qquad \qquad \qquad \qquad \qquad \qquad \qquad \\
    \cong S(U, \wt V \otimes_{\zpb} \bigotimes_{v' \ne v} W_{\lambda_{v'}} \otimes_{\zpb} \fpb)
\qquad \qquad \qquad \qquad \qquad \qquad \quad \\
    \cong \big((\wt V \otimes_{\zpb} \fpb) \otimes_{\fpb} S\sm(U^v, \otimes_{v' \ne v} (W_{\lambda_{v'}} \otimes_{\zpb}
    \fpb))\big)^{\GL_n(\O_{F_w})}\\ 
 \cong   (\wt V \otimes_{\zpb} \Slift \otimes_{\zpb} \Fpbar)^{\GL_n(\O_{F_w})}.
  \end{multline*}
  In the second case note that the map
  is injective because $\qpb$ is $\zpb$-flat. Any element of the right-hand side can be written as $m \otimes \lambda$,
  where $m \in \wt V \otimes_{\zpb} \Slift$ and $\lambda \in \qpb\s$. Since $\wt V \otimes_{\zpb} \Slift$ is $\zpb$-torsion
  free it follows that $m$ is $\GL_n(\O_{F_w})$-stable, as required.
  
  To verify Axiom~\ref{axiom: crystalline lifts}, let $\lambda_w \in
  (\Z^n_+)^{\Hom(F_w,\Qpbar)}$. Our definition of $\Slift$ only
  depended on $\lambda$ away from $w$, $w^c$, so we are free to modify $\lambda$ at $w$ and $w^c$ so that it agrees with 
  the given $\lambda_w$ at $w$.
  By Lemma~\ref{lem:mapping-in-locally-algebraic-V} and Axiom~\ref{axiom: finite and ext of scalars}
  we have $\Slift(W_{\lambda_w} \otimes_{\zpb} \qpb) = S(U,W_\lambda \otimes_{\zpb} \qpb)$, so by assumption
  $S_\lambda(U,\Qpbar)_\mf{m}\ne 0$. It follows from Lemma~\ref{lem:Deligne--Serre}, as well as Theorem~\ref{thm:
    existence of Galois reps attached to algebraic modular forms}(iii)
  together with Corollary 3.1.2 of \cite{cht} and the Chebotarev density theorem
  that there is an
  irreducible $G(\A_{F^+}^\infty)$-subrepresentation $\pi$ of
  $S_\lambda(\Qpbar)$ with $\rbar_\pi\cong\rbar^\vee$. By Theorem
  \ref{thm: existence of Galois reps attached to algebraic modular
    forms}, $r_\pi|_{G_{F_w}}$ is crystalline with $\HT_\tau(r_\pi|_{G_{F_w}}) =
  (\lambda_{\tau,1}+n-1,\dots,\lambda_{\tau,n})$ for all $\tau : F_w \to \Qpbar$, so $r^\vee_\pi|_{G_{F_w}}$
  provides the required lift. The claim about the characteristic
  polynomial of $\varphi^{[F_w:\qp]}$ follows from Theorem~\ref{thm:
    existence of Galois reps attached to algebraic modular forms}(iv)
  together with Corollary 3.1.2 of \cite{cht}.

  It remains to verify Axiom~\ref{axiom: tame lifts}.  Suppose that $R$ is a representation of $\GL_n(k_w)$ as in Axiom~\ref{axiom: tame lifts}.
  As our definition of $\Slift$ only
  depended on $\lambda$ away from $w$, $w^c$, we are free to modify $\lambda$ at $w$ and $w^c$ so that $\lambda_w = 0$ at~$w$.
  By choosing a $\GL_n(k_w)$-stable lattice $R_0$ in $R$
  and applying Axiom~\ref{axiom: finite and ext of scalars} 
  and Lemma~\ref{lem:mapping-in-locally-algebraic-V} to the smooth representation $\wt V = R_0$
  we have $\Slift(R) \cong S(U, R \otimes_{\qpb} (W_{\lambda} \otimes_{\zpb} \qpb))$.
  Therefore
  \begin{align*}
    \Slift(R)_\m &\cong S(U, R \otimes_{\qpb} (W_{\lambda} \otimes_{\zpb} \qpb))_\m \\
    & \cong S(U', R \otimes_{\qpb} (W_{\lambda} \otimes_{\zpb} \qpb))^{U/U'}_\m \\
    & \cong \Hom_{U_v}\big(R\dual, S(U',W_\lambda \otimes_{\zpb} \qpb)_\m\big),
  \end{align*}
  where $U'$ is the kernel of the map $U \to U_v \xrightarrow{\iota_w} \GL_n(\O_{F_w}) \to \GL_n(k_w)$.
  (In the second line, note that $u \in U$ acts on a function $f : G(F^+)\backslash G(\A_{F^+}^\infty) \to R \otimes_{\qpb} (W_{\lambda} \otimes_{\zpb} \qpb)$
  by $(u\cdot f)(g) = u_p f(gu)$.)
  As in the preceding paragraph (the verification of Axiom~\ref{axiom: crystalline lifts}) there is thus an
  irreducible subrepresentation $\pi$ of $S_\lambda(\Qpbar)$ such that $\pi^{U'} \cap S(U',W_\lambda \otimes_{\zpb} \qpb)_\m \ne 0$
  and $\pi|_{\cG(\cO_{F_v^+})}$ contains $R\dual$.
  In particular, $\rbar_\pi\cong\rbar^\vee$. 
  Since $\lambda_w=0$, it follows from Theorem \ref{thm:
    existence of Galois reps attached to algebraic modular
    forms}
  and
  Proposition \ref{prop:inertial-llc} that $r_\pi^\vee|_{G_{F_w}}$ provides
  the required lift.
\end{proof}

\subsection{The main theorem for definite unitary
  groups}\label{subsection:main theorem for unitary groups}

From now on we assume that $n = 3$ and that $p$ splits in $F$.

The set $W^?(\rbar)$ of predicted Serre weights for $\rbar$ is defined in terms of the set
of predicted weights $W^?(\rbar|_{G_{F_w}})$ for the
representations $\rbar|_{G_{F_w}}$ for $w|p$ a place of $F$ (see Subsection~\ref{subsec:weights}).
Namely, if $v$ is a place of $F^+$ lying over~$p$,
if $w$ is a place of $F$ lying over $v$, 
if $V_v$ is an irreducible representation of $\cG(k_v)$ over $\Fbar_p$,
and if $w$ is a place of $F$ lying over $v$, 
then we introduce the following condition:
\begin{equation}
\label{eqn:condition}
V_v \circ \iota_w^{-1} \in W^?_w(\rbar|_{G_{F_w}}).
\end{equation}
\begin{df}
\label{df:Serre weights}
We define the set of predicted Serre weights for $\rbar$ to be
\begin{equation*}
  W^?(\rbar) := \{ V = \otimes_{v|p} V_{v} : V_{v} \text{ satisfies condition~(\ref{eqn:condition})
 for all } w|v|p \}.
\end{equation*}
\end{df}

\begin{remark}
In Definition~\ref{df:Serre weights},
it suffices to fix a place $\tv$ over each place $v|p$ of $F^+$,
and then to impose the condition~(\ref{eqn:condition})
just at these chosen places $\tv$.
This follows from the facts
that
  $\rbar^c\cong\rbar^\vee \otimes\epsilonbar^{2}$,
that 
$W^?_w(\rbar|_{G_{F_{w}}})\dual = W^?_w(\rbar|_{G_{F_{w}}}\dual \otimes \varepsilonbar^2)$
(by Proposition 6.23 in~\cite{bib:herzig-thesis})
and that $\iota_{w^c}$ is conjugate to the inverse
transpose of $c\circ\iota_w$.   Indeed, taken together,
these imply that imposing condition~(\ref{eqn:condition})
at $\tv$ is equivalent to imposing it at $\tv^c$. 
\end{remark}

\begin{defn}\label{defn: generic global Serre weight}
  We say that a Serre weight $V=\otimes_{v|p}V_v$ is \emph{generic}
  \(resp.\ \emph{strongly generic}\) if $V_v\circ\iota_w^{-1}$ is
  a generic \(resp.\ strongly generic\) weight in the sense of Definition~\ref{defn: strongly generic weight}
  for all places $w|p$ of $F$. We let
  $W\gen(\rbar)$ denote the set of generic Serre weights in $W(\rbar)$.
\end{defn}

\begin{thm}\label{thm:elimination unitary global} Let $F$ be an imaginary CM
  field in which the prime $p$ splits completely, and suppose that
  $\rbar: G_F \to \GL_3(\fpb)$ is a continuous representation.
  For all places $w | p$ of $F$ suppose that $\rbar|_{G_{F_w}}$ is irreducible.
  Then $W\gen(\rbar) \subset W^?(\rbar)$.
\end{thm}
\begin{proof}
  Fix a Serre weight $V =\otimes_{v|p}V_v\in W\gen(\rbar)$. Since the definition of
  $W^?(\rbar)$ is local, we need only prove that if $w|p$ is a place
  of $F$ with $w|_{F^+}=v$, then $V_v \circ \iota_w^{-1}\in W^?_w(\rbar|_{G_{F_w}})$. 
  We set up the axiomatic framework as in Subsection~\ref{sec:verify-char-zero} using the place $w$ and the Serre weight $V$.
  Then we see by Lemma~\ref{lem:mapping-in-smooth-V} that $V_v \circ \iota_w^{-1}\in W_w(\rbar)$.
  We can now conclude by Proposition \ref{prop:the axioms hold for unitary
    groups} and Theorems \ref{thm:char-0-axioms-imply-char-p-axioms}, \ref{thm:elimination}.
\end{proof}

\begin{thm}\label{thm:main unitary global}
  Let $F$ be an imaginary CM field in which the prime $p$ splits
  completely, and suppose that $\rbar: G_F \to \GL_3(\fpb)$ is a
  continuous representation.
  For all places $w | p$ of $F$ suppose that $\rbar|_{G_{F_w}}$ is irreducible.
  If $\rbar$ is modular of some strongly generic Serre weight, then $W\gen(\rbar) = W^?(\rbar)$.
\end{thm}
\begin{proof}
  By Theorem \ref{thm:elimination unitary global} we have
  $W\gen(\rbar) \subset W^?(\rbar)$. By assumption there is some $V=\otimes_{v|p}V_v\in
  W\gen(\rbar)$ such that $\rbar$ is modular of Serre weight $V$. Fix some
  place $w|p$ of $F$. Arguing with each such $w$ in turn, one over each place $v|p$ of $F^+$, and iterating, we see that
  it suffices to show that if $V'_v\circ\iota_w^{-1}\in
  W^?_w(\rbar|_{G_{F_w}})$ then $(\otimes_{v'\ne v}V_{v'})\otimes
  V'_v\in W(\rbar)$.  Setting up the axiomatic framework as in Subsection~\ref{sec:verify-char-zero} using the place $w$
  and the Serre weight $V$, we see that this follows from Proposition \ref{prop:the
    axioms hold for unitary groups} and Theorems \ref{thm:char-0-axioms-imply-char-p-axioms}, \ref{thm:main}.
\end{proof}

\subsection{Automorphic induction}\label{examples:automorphic induction}
In this brief section
we show that one can construct many examples in
which our main theorem (Theorem~\ref{thm:main}) applies. Firstly, choose a CM field $F$ with
maximal totally real subfield $F^+$ such that $F/F^+$ is unramified at
all finite places and such that $p$ splits completely in $F$. (For
example, choose any CM field $F$ in which $p$ splits completely, and use Lemma~4.1.2
of~\cite{cht} to find a Galois totally real field $K^+$ in which $p$
splits completely such that the
extension $FK^+/F^+K^+$ is everywhere unramified. The condition that
$F/F^+$ is unramified at all finite places is used below in order to
apply Theorem 5.4 of \cite{labesse}.) 
Fix a weight $\lambda\in (\Z^3_+)^{\Hom(F,\Qpbar)}_0$. 

By Theorem~6 of chapter~10 of \cite{MR2467155} (and its proof, which
shows that we may assume that $p$ is unramified in the extension), we may choose a cyclic
totally real extension $L^+/F^+$ of degree 3 in which every place of
$F^+$ lying over $p$ is inert. Let $L=FL^+$. By Lemma 4.1.5 of \cite{cht}, we may choose
a crystalline character $\psi:G_L\to\qpb\s$ such that
\begin{itemize}
\item $\psi\psi^{c}=\varepsilon^{-2}$,
\item if $w|p$ is a place of $F$ and $\tau:F_w\into\Qpbar$, then 
  \[ \bigcup \HT_{\wt\tau}(\psi|_{G_{L_w}}) =\{\lambda_{\tau,1}+2,\lambda_{\tau,2}+1,\lambda_{\tau,3}\},\]
  where $\wt\tau$ runs through all embeddings $L_w \into \Qpbar$ such that $\wt\tau|_{F_w} = \tau$.
\end{itemize}
Assume furthermore that for each $\tau$ we have
$\lambda_{\tau,1}-\lambda_{\tau,3}\le p-2$. Set
$r_\psi:=\Ind_{G_L}^{G_F}\psi$. Then for each place $w|p$ of $F$ the
representation $\rbar_\psi^\vee|_{G_{F_w}}$ is irreducible (use the
bounds on the $\lambda_{\tau,i}$ and the assumption that each place of
$F^+$ lying over $p$ is inert in $L^+$), so in particular the global
representation $\rbar_\psi^\vee$ is irreducible.  Then by Theorems 4.2
and 5.1 of \cite{MR1007299}, together with Theorem 5.4 of
\cite{labesse} and Lemma~\ref{lem: relationship of algebraic
  automorphic forms to classical automorphic forms}
above, it follows that $\rbar_\psi^\vee$ is modular of
Serre weight $F_a$,
  where $a\in(\Z^3_+)_0^{\coprod_{w|p}\Hom(k_w,\Fpbar)}$ is the unique restricted weight of which $\lambda$ is the
  lift. (Note that $W_a \otimes_{\zpb} \fpb \cong F_a$ due to our bounds on the $\lambda_{\tau,i}$,
  for example by Proposition~3.18 in~\cite{bib:herzig-thesis}. Note
  also that there is a unitary group satisfying our assumptions which
  is furthermore quasisplit at all finite places, cf.\ Section 2.1 of \cite{ger}.)

  Thus if we choose any restricted weight
  $a\in(\Z^3_+)_0^{\coprod_{w|p}\Hom(k_w,\Fpbar)}$ such that $F_{a_w}$
  is strongly generic and in the lower alcove for each $w|p$, the representation $\rbar_\psi\dual$
  constructed above (with $\lambda$ the lift of $a$) satisfies the
  assumptions of Theorem \ref{thm:main}. (Note that the condition that
  $\lambda_{\tau,1}-\lambda_{\tau,3}\le p-2$ is automatic from the
  assumption that all $F_{a_w}$ are strongly generic.)

\bibliography{egh}

\bibliographystyle{amsalpha} 

\end{document}